\documentclass[leqno,12pt,a4paper]{article}
\usepackage[margin=0.97in]{geometry}
\usepackage[authoryear,compress]{natbib}
\usepackage{etoolbox} \usepackage[amsthm,thmmarks]{ntheorem}
\usepackage{thmtools}
\usepackage{amsfonts,amssymb,amsmath}
\usepackage{xypic,mathtools}
\usepackage{times}
\usepackage{microtype}
\usepackage{hyphenat}
\usepackage[usenames,dvipsnames]{color}
\usepackage{rotating}
\usepackage{longtable}
\usepackage{caption}

\usepackage{epigraph}
\usepackage{enumitem}
\setlist[enumerate]{listparindent=0.5in}
\usepackage{mathrsfs}
\DeclareMathAlphabet{\mathscrbf}{OMS}{mdugm}{b}{n} \newcommand{\be}{\begin{equation}}
\newcommand{\ee}{\end{equation}}
\newcommand{\bes}{\begin{equation*}}
\newcommand{\ees}{\end{equation*}}
\newcommand{\bea}{\begin{eqnarray}}
\newcommand{\eea}{\end{eqnarray}}
\newcommand{\beas}{\begin{eqnarray}}
\newcommand{\eeas}{\end{eqnarray}}
\newcommand{\ben}{\begin{note}}
\newcommand{\een}{\end{note}}
\newcommand{\bexl}{\vskip0.1em\noindent\hrulefill\vskip1em\begin{ExerciseList}}
\newcommand{\eexl}{\end{ExerciseList}\hrulefill}

\newcommand{\bthm}{\begin{theorem}}
\newcommand{\ethm}{\end{theorem}}
\newcommand{\bpro}{\begin{prop}}
\newcommand{\epro}{\end{prop}}
\newcommand{\bcor}{\begin{corollary}}
\newcommand{\ecor}{\end{corollary}}
\newcommand{\bcon}{\begin{conjecture}}
\newcommand{\econ}{\end{conjecture}}
\newcommand{\bp}{\begin{proof}}
\newcommand{\ep}{\end{proof}}
\newcommand{\blem}{\begin{lemma}}
\newcommand{\elem}{\end{lemma}}
\newcommand{\bn}{\begin{note}}
\newcommand{\en}{\end{note}}
\newcommand{\benum}{\begin{enumerate}}
\newcommand{\eenum}{\end{enumerate}}
\newcommand{\bed}{\begin{defn}}
\newcommand{\eed}{\end{defn}}
\newcommand{\brem}{\begin{remark}}
\newcommand{\erem}{\end{remark}}

\newcommand{\btik}{\begin{tikzpicture}\begin{axis}[scale=0.5,axis y line=center, axis x line=middle]}
\newcommand{\etik}{\end{axis}\end{tikzpicture}}

\let\into=\hookrightarrow
\let\mapsto=\longmapsto

\newcommand{\upperRomannumeral}[1]{\uppercase\expandafter{\romannumeral#1}}

 \usepackage{tikz-cd}
\usepackage{stmaryrd}
\let\cite=\citep
\usepackage{fancyhdr}
\usepackage{titlesec}
\usepackage{fancyhdr}
\usepackage[hypertexnames=false,colorlinks=true,hyperindex, citecolor=blue, urlcolor=cyan]{hyperref}
\usepackage{zref-clever}
	\let\Cref=\zcref

\theoremstyle{theorem}
\newtheorem{theorem}[equation]{Theorem}      \newtheorem{lemma}[equation]{Lemma}          \newtheorem{corollary}[equation]{Corollary}  \newtheorem{proposition}[equation]{Proposition}
\newtheorem{prop}[equation]{Proposition}

\theoremstyle{definition}
\newtheorem{conj}[equation]{Conjecture}
\newtheorem{conjecture}[equation]{Conjecture}
\newtheorem{example}[equation]{Example}
\newtheorem{question}[equation]{Question}

\theoremstyle{definition}
\newtheorem{defn}[equation]{Definition}
\theoremsymbol{\ensuremath{\bullet}}

\theoremsymbol{\ensuremath{\bullet}}

\theoremstyle{remark}

\theoremsymbol{\ensuremath{\bullet}}

\theoremstyle{definition}
\newtheorem{remark}[equation]{Remark}
\theoremsymbol{\ensuremath{\bullet}}

\numberwithin{equation}{section}

\newcommand{\bdefn}{\begin{defn}}
\newcommand{\edefn}{\end{defn}}

\let\into=\hookrightarrow
\let\isom=\simeq

\let\tensor=\otimes
\newcommand{\A}{\mathcal{A}}
\newcommand{\abs}[1]{\left\vert#1\right\vert}

\newcommand{\bF}{{\bar{F}}}
\newcommand{\bQ}{{\bar{\Q}}}

\newcommand{\C}{{\mathbb C}}

\newcommand{\End}{\rm{End}}

\newcommand{\Ext}{{\rm Ext}\,}
\newcommand{\F}{{\mathbb F}}

\newcommand{\gal}{{\rm Gal}}
\newcommand{\Gm}{\mathbb{G}_m}

\newcommand{\mydot}{{\small{\bullet}}}
\newcommand{\N}{\mathcal{N}}
\newcommand{\norm}[1]{\left\Vert#1\right\Vert}

\newcommand{\Q}{{\mathbb Q}}
\newcommand{\R}{{\mathbb R}}

\newcommand{\Spec}{{\rm Spec}}

\newcommand{\Z}{{\mathbb Z}}

\renewcommand{\O}{{\mathcal O}}
\renewcommand{\P}{{\mathbb P}}
\renewcommand{\wp}{{\mathfrak p}}

\newcommand{\invlim}{\varprojlim}

\newcommand{\mapright}[1]{{\xymatrix{{}\ar[r]^{#1}&{}}}}

\newcommand{\anabmapright}[1]{\stackrel{#1}{\leftrightsquigarrow}}

\renewcommand{\bpro}{\begin{proposition}}
	\renewcommand{\epro}{\end{proposition}}
\renewcommand{\bcon}{\begin{conj}}
	\renewcommand{\econ}{\end{conj}}

\let\mathcal=\mathscr

\title{On Mochizuki's idea of Anabelomorphy and its applications
}
\author{Kirti Joshi}

\begin{document}
\maketitle

\lhead{}

\tableofcontents 
\newpage
\epigraphwidth0.45\textwidth
\epigraph{\emph{Jean-Marc Fontaine}}{In Memoriam}

\newcommand{\act}{\curvearrowright}
\newcommand{\lmp}{{\Pi\act\Ot}}
\newcommand{\lmpi}{{\lmp}_{\int}}
\newcommand{\lmpf}{\lmp_F}
\newcommand{\Om}{\O^{\times\mu}}
\newcommand{\Omf}{\O^{\times\mu}_{\bF}}
\renewcommand{\N}{\mathbb{N}}
\newcommand{\yoga}{Yoga}
\newcommand{\gl}[1]{{\rm GL}(#1)}
\newcommand{\bK}{\bar{K}}
\newcommand{\bL}{\bar{L}}
\newcommand{\reptrip}{\rho:G_K\to\gl V}
\newcommand{\reptripp}[1]{\rho\circ\alpha:G_{\ifstrempty{#1}{K}{{#1}}}\to\gl V}
\newcommand{\benumlab}{\begin{enumerate}[label={{\bf(\arabic{*})}}]}
\newcommand{\ord}{\mathop{\rm ord}\nolimits}	
\newcommand{\kcs}{K^\circledast}
\newcommand{\lcs}{L^\circledast}
\renewcommand{\A}{\mathbb{A}}
\newcommand{\bfq}{\bar{\mathbb{F}}_q}
\newcommand{\tripod}{\P^1-\{0,1728,\infty\}}

\newcommand{\vseq}[2]{{#1}_1,\ldots,{#1}_{#2}}
\newcommand{\anab}[4]{\left({#1},\{#3 \}\right)\anabelmap\left({#2},\{#4 \}\right)}

\newcommand{\GL}{\rm GL}
\newcommand{\gln}[1]{\GL_n(#1)}
\newcommand{\glo}[1]{\GL_1(#1)}
\newcommand{\glt}[1]{\GL_2(#1)}

\newcommand{\iut}{\cite{mochizuki-iut1, mochizuki-iut2, mochizuki-iut3,mochizuki-iut4}}
\newcommand{\topics}{\cite{mochizuki-topics1,mochizuki-topics2,mochizuki-topics3}}
\newcommand{\iutfour}{\cite{mochizuki-iut1-4}}
\let\iut=\iutfour
\newcommand{\linv}{\mathfrak{L}}
\newcommand{\bedef}{\begin{defn}}
\newcommand{\eedef}{\end{defn}}

\titleformat{\subsection}[runin]{\normalfont\bfseries}{\S\ \thesubsection}{.5em}{}[{\ \ }]
\titlespacing{\subsection}{0pt}{1.5ex plus .1ex minus .2ex}{0pt}

\zcsetup{cap} \AddToHook{env/theorem/begin}{\zcsetup{countertype={equation=theorem}}}
\AddToHook{env/corollary/begin}{\zcsetup{countertype={equation=corollary}}}
\AddToHook{env/lemma/begin}{\zcsetup{countertype={equation=lemma}}}
\AddToHook{env/proposition/begin}{\zcsetup{countertype={equation=proposition}}}
\AddToHook{env/prop/begin}{\zcsetup{countertype={equation=prop}}}
\AddToHook{env/theoremdef/begin}{\zcsetup{countertype={equation=theoremdef}}}
\AddToHook{env/defn/begin}{\zcsetup{countertype={equation=defn}}}
\AddToHook{env/definition/begin}{\zcsetup{countertype={equation=definition}}}
\AddToHook{env/remark/begin}{\zcsetup{countertype={equation=remark}}}
\AddToHook{env/rem/begin}{\zcsetup{countertype={equation=rem}}}
\AddToHook{env/conj/begin}{\zcsetup{countertype={equation=conj}}}
\AddToHook{env/conjecture/begin}{\zcsetup{countertype={equation=conjecture}}}
\AddToHook{env/question/begin}{\zcsetup{countertype={equation=question}}}
\AddToHook{env/example/begin}{\zcsetup{countertype={equation=example}}}

\zcRefTypeSetup{equation}{
	Name-sg = eq. ,
	name-sg = eq. ,
	Name-pl = eqns. ,
	name-pl = eqns. ,
}

\zcRefTypeSetup{theorem}{
	Name-sg = Theorem ,
	name-sg = theorem ,
	Name-pl = Theorems ,
	name-pl = theorems ,
}

\zcRefTypeSetup{theoremdef}{
	Name-sg = Theorem-Definition ,
	name-sg = theorem-definition ,
	Name-pl = Theorem-Definitions ,
	name-pl = theorem-definitions ,
}

\zcRefTypeSetup{corollary}{
	Name-sg = Corollary ,
	name-sg = corollary ,
	Name-pl = Corollaries ,
	name-pl = corollaries ,
}

\zcRefTypeSetup{proposition}{
	Name-sg = Proposition ,
	name-sg = proposition ,
	Name-pl = Propositions ,
	name-pl = propositions ,
}

\zcRefTypeSetup{prop}{
	Name-sg = Proposition ,
	name-sg = proposition ,
	Name-pl = Propositions ,
	name-pl = propositions ,
}

\zcRefTypeSetup{lemma}{
	Name-sg = Lemma ,
	name-sg = lemma ,
	Name-pl = Lemmas ,
	name-pl = lemmas ,
}

\zcRefTypeSetup{remark}{
	Name-sg = Remark ,
	name-sg = remark ,
	Name-pl = Remarks ,
	name-pl = remarks ,
}

\zcRefTypeSetup{question}{
	Name-sg = Question ,
	name-sg = question ,
	Name-pl = Questions ,
	name-pl = questions ,
}

\zcRefTypeSetup{example}{
	Name-sg = Example ,
	name-sg = example ,
	Name-pl = Examples ,
	name-pl = examples ,
}

\zcRefTypeSetup{defn}{
	Name-sg = Definition ,
	name-sg = definition ,
	Name-pl = Definitions ,
	name-pl = definitions ,
}

\zcRefTypeSetup{conj}{
	Name-sg = Conjecture ,
	name-sg = conjecture ,
	Name-pl = Conjectures ,
	name-pl = conjectures ,
}

\zcRefTypeSetup{conjecture}{
	Name-sg = Conjecture ,
	name-sg = conjecture ,
	Name-pl = Conjectures ,
	name-pl = conjectures ,
}

\newcommand{\fixnumberwithin}[1]{
\numberwithin{equation}{#1}
	\numberwithin{theorem}{#1}
	\numberwithin{conj}{#1}
	\numberwithin{conjecture}{#1}
	\numberwithin{lemma}{#1}
	\numberwithin{proposition}{#1}
	\numberwithin{prop}{#1}
	\numberwithin{corollary}{#1}
	\numberwithin{defn}{#1}
	\numberwithin{definition}{#1}
	\numberwithin{remark}{#1}
	\numberwithin{rem}{#1}
	\numberwithin{question}{#1}
}

\newcommand{\nws}{\fixnumberwithin{section}}
\newcommand{\nwss}{\fixnumberwithin{subsection}}
\newcommand{\nwsss}{\fixnumberwithin{subsubsection}}

\newcommand{\ENDDOCUMENT}{\bibliography{../../master/masterofallbibs.bib}\end{document}}

\newcommand{\arithspaces}{\cite{joshi-teich,joshi-untilts,joshi-teich-estimates,joshi-teich-def,joshi-teich-rosetta}}

\newcommand{\ssep}{\S\,}
\newcommand{\phisen}{\Phi_{\rm Sen}}
\numberwithin{example}{subsection}
\section{Introduction}\label{se:intro}
\newcommand{\anabelmap}{\leftrightsquigarrow}
\subsection{What is Anabelomorphy?}\label{ss:intro1.1}
The term \emph{anabelomorphy}  (pronunciation guide \textit{anabel-o-morphy}; the root of this term is in Alexander Grothendieck's Anabelian Program) is coined (by the Author) and introduced here as a concise way of expressing Shinichi Mochizuki's notion of an anabelian way of changing base fields or base rings.  Roughly speaking, one may understand anabelomorphy as the branch of arithmetic in which one studies arithmetic by fixing the absolute Galois group of a field rather than the field itself. 
For $p$-adic fields, the idea of anabelomorphy is firmly grounded in two theorems (1) the well-known theorem of Jarden-Ritter (\Cref{th:fourth-fun-anab}) which provides a necessary and sufficient condition for the absolute Galois groups of two $p$-adic fields to be topologically isomorphic (for explicit examples of such non-isomorphic $p$-adic fields, see \Cref{le:basic-example}), and (2) a well-known theorem of Mochizuki  (\Cref{th:second-fun-anab}) which asserts that a $p$-adic field is  determined by the absolute Galois group equipped with the upper numbering ramification filtration but not determined by the absolute Galois group.

This leads to the definition (\Cref{def:anabelomorphy-fields}) of anabelomorphic $p$-adic fields: two $p$-adic fields $K,L$ are anabelomorphic (denoted $K\anabmapright{\alpha} L$) if there exists a topological isomorphism $G_K\overset{\alpha}{\isom} G_L$ of their absolute Galois groups (for some choice of algebraic closures of $K$ and $L$). Anabelomorphism of $p$-adic fields is an equivalence relation. Anabelomorphic $p$-adic fields have naturally isomorphic multiplicative group and additive groups i.e. the multiplicative structure and the additive structure  of the field remains separately fixed within the  anabelomorphism class, but the fields may not be isomorphic (\Cref{th:third-fun-anab}). This means that the intertwining between the additive and multiplicative structure of the  field deforms or wiggles   around within the anabelomorphism class, or as Mochizuki views it,  the multiplicative structure remains fixed but the additive structure deforms. Since the absolute Galois group of a field is its \'etale fundamental group,  this behavior is quite analogous to the case of (compact, connected) Riemann surfaces of a fixed genus (and hence isomorphic fundamental groups) but with possibly non-isomorphic (sheaves of) rings of complex functions. So the term anabelomorphy can be applied to Riemann surfaces: for example, any pair of points of any Teichm\"uller space gives rise to anabelomorphic Riemann surfaces.

The first key realization which emerges  from (1) and (2), and upon which this paper is founded, is that the upper numbering ramification filtration is the Galois-theoretic stand-in for the  intertwining between the additive and multiplicative structures of a $p$-adic field and therefore the field structure manifests itself wherever one encounters the upper numbering ramification filtration.

The second key realization of this paper is this:  anabelomorphic $p$-adic fields $K,L$ have topologically isomorphic absolute galois groups, hence representations of $G_K$ can be viewed as representations of $G_L$ and vice versa.

The third key realization is this: suppose $K$ is a $p$-adic field and $\bK$ is an algebraic closure of $K$. Then there are at most finitely many fields $L\subset \bK$ with an anabelomorphism $K\anabelmap L$. Because  all these finitely many fields $L$ have isomorphic absolute Galois groups i.e. isomorphic \'etale fundamental groups, this set of fields can be considered as a zero-dimensional or discrete Teichm\"uller space.

The fourth key realization is that  anabelomorphy and the phenomena that go with it appear in many diverse contexts. Here is a surprising example: anabelomorphy is a key feature of the theory of perfectoid fields and spaces, and like the $p$-adic case, (suitably defined) multiplicative structures remains fixed.  Moreover, theorems arising in anabelomorphy of $p$-adic fields have parallels in the perfectoid setting: for example (\Cref{th:third-fun-anab}{\bf(3)} and \Cref{th:perfect-anab}) and (\Cref{th:anab-proj-spaces} and \Cref{th:anab-proj-spaces-perf}). 

These four realizations make it clear that anabelomorphy and the phenomena that go with it have consequences for number theory and algebraic geometry, and anabelomorphy provides a common umbrella under which they may be studied. For completeness, \ssep\ref{se:five-fundamental} recalls several classical results of anabelomorphy. A few of the themes explored in this paper are  \ssep\ref{se:anabel-galois-reps} (Galois representations),  \ssep\ref{se:p-adic-hodge} ($p$-adic Hodge Theory), \ssep\ref{se:anab-langlands} (local Langlands Correspondence).   The archimedean case (monodromy and mixed Hodge Theory) is discussed in \ssep\ref{se:anab-arch}. In \ssep\ref{se:open-question}, readers will find  several open questions which can serve as starting points for new investigations.  These questions show that anabelomorphy has applications to other topics of current interest.

The rest of this Introduction provides a detailed discussion of the results.

\subsection{Amphoricity and anabelomorphy of $p$-adic fields}\label{ss:intro1.2}
Anabelomorphy leads us naturally to introduce the following notion. 

A quantity (resp. a property, an algebraic structure) associated with a $p$-adic field is said to be \emph{amphoric} (\Cref{def:amphoric-var}) if two $p$-adic fields $K,L$ in the same anabelomorphism class have the same quantity (resp. same property, have naturally isomorphic algebraic structures). [The word \textit{amphoric} has its root in the word \textit{amphora} which was a storage jar used in Ancient Greece and Rome.] For examples of classically known amphoric quantities, properties and structures see \Cref{th:third-fun-anab}.
[For a more geometric, topos theoretic view of amphoricity, see \Cref{re:topos-of-gk}.]

Classical results (recalled in \ssep\ref{se:five-fundamental}) of anabelomorphy of $p$-adic fields  have focused on amphoric quantities, properties and structures. However, it is important to understand quantities, properties and structures which are not amphoric. Such quantities are quite common, even in the theory of Riemann surfaces (\Cref{th:hodge3}) and their existence signals the existence of Teichm\"uller Theory (\Cref{pr:linv-inf-unamph}--its $p$-adic analog is \Cref{th:l-inv-unamphoric}). It was in this context, that the Author discovered \Cref{th:discriminant-is-unamphoric} which shows that the absolute discriminant and the absolute different of a $p$-adic field are not amphoric. This observation and computations of \ssep\ref{ss:weak-anab-elliptic} were the starting point of this paper.

\subsection{Anabelomorphy and   Galois representations}\label{ss:intro1.3}
One has the following results: the category of ordinary ($\ell$-adic and $p$-adic) representations is amphoric (\Cref{th:ordinary-amphoric}), the maximal tamely ramified and maximal unramified extensions of anabelomorphic $p$-adic fields are anabelomorphic (\Cref{le:anab-knr}), the property of being peu-ramifi\'ee or tres ramifi\'ee are not amphoric (\Cref{th:peu-tres}), Frobenius elements are amphoric (\Cref{th:frob-amphoric}) and hence  $L$-functions of $p$-adic Galois representations are amphoric (\Cref{th:char-poly-amphoric}). The Artin and Swan conductors are not amphoric (\Cref{th:artin-swan-unamphoric} and also \Cref{th:discriminant-is-unamphoric2}). \Cref{pr:iwasawa-coh} records the amphoricity of the Iwasawa cohomology.

Section \ssep\ref{se:p-adic-hodge} considers anabelomorphy and $p$-adic Hodge Theory and proves the following:
the property of being a crystalline representation (and  hence  of being a Hodge-Tate representation) is not amphoric in general (\Cref{th:crystalline-not-amphoric}), but the property of being a Hodge-Tate representation of pure weight is amphoric (\Cref{th:sen-thm-and-anabelomorphy}). Anabelomorphic $p$-adic fields have anabelomorphic cyclotomic fields of norms 
 and the  category of  \'etale $(\varphi,\Gamma)$-modules is amphoric (\Cref{th:phi-gamma}). The property of the Sen-invariant $\phisen(\rho,V)$ being ``semisimple and has integer eigenvalues'' is not amphoric. The Fontaine subspace $H^1_f(G_K,\Q_p(1))$ (of ordinary crystalline two-dimensional representations of $G_K$) and the subspace $H^1_e(G_K,\Q_p(1))$ are amphoric (\Cref{th:anab-galois-h1}). The
 $\linv$-invariant of an ordinary two-dimensional $p$-adic representation is not amphoric (\Cref{th:l-inv-unamphoric}). Consequences for deformation theory of Galois representations are given in \Cref{th:anab-deform-thry}.

\newcommand{\sS}[1]{\mathscr{S}(#1)}
\newcommand{\dmu}[2]{d\mu_{#1}({#2})}
\newcommand{\dmus}[2]{d\mu_{#1}^*({#2})}

\subsection{Anabelomorphy and the local Langlands Correspondence}\label{ss:intro1.4}
\Cref{le:weil-deligne-groups} establishes an isomorphism between the Weil and the Weil-Deligne groups of anabelomorphic $p$-adic fields. Since
the local Langlands Correspondence \cite{henniart1988} establishes a bijection between semisimple representations of the Weil-Deligne group of a $p$-adic field $K$ and the set of irreducible admissible representations of $\gln K$,  the natural question which arises is this: 
\numberwithin{question}{subsection}
\begin{question}\label{qu:basic-langlands}
If $K\anabelmap L$ are anabelomorphic $p$-adic fields, then how are the irreducible admissible representations  of $\gln{K}$ and $\gln L$ related? 
\end{question}

This question is taken up in \ssep\ref{se:anab-langlands}.  This leads to the following results:  suppose one has anabelomorphic $p$-adic fields $K\anabelmap L$. The Schwartz spaces $\sS{K}$ and $\sS{K^*}$ are amphoric (\Cref{th:meas-thry}). \Cref{th:anab-glo-reps} establishes the $\GL_1$ case.   \Cref{th:automorphic-ordinary-synchronization}, shows that  one has a natural bijection between  principal series representations of $\gln{K}$ and $\gln L$ (the Galois analog of this is \Cref{th:ordinary-amphoric}). For $(p,n)=1$,  one has a natural bijection between all irreducible supercuspidal representations of $\gln{K}$ and $\gln{L}$ (\Cref{th:syn-super-cusp-n}). \Cref{pr:central-simple} sets up a natural bijection between central division algebras over $K$ and $L$. For $p>2$, one can also synchronize Weil representations of $\glt K$ and $\glt L$ (\Cref{th:weil-rep-gl2}).  So for example,  for any odd prime $p$, one has a complete correspondence (compatible with the local Langlands correspondence) between irreducible admissible representations  of $\glt K$ and $\glt{L}$ (\Cref{th:anab-local-l-gl2-odd}). The situation for $\GL_n$  with $p| n$ needs substantial clarification and remains open.

\subsection{Constructions of varieties via anabelomorphy}
\Cref{th:anab-affine-spaces}, \Cref{th:anab-proj-spaces}, \Cref{th:anab-toric} establish the relationships between affine, projective spaces and smooth, projective toric varieties over anabelomorphic $p$-adic fields (as is noted in \ssep\ref{se:perfectoids}, the perfectoid analogs of these are due to \cite{scholze12-perfectoid-ihes}). \Cref{th:anabelomorphy-tate-curves} and \Cref{th:anabelomorphy-abelian-var} provide constructions of Tate elliptic curves and  abelian varieties with split multiplicative reduction from such varieties over anabelomorphic fields. In all of these cases the (constructed) varieties are themselves anabelomorphic. In \Cref{th:grp-ord-p}, \Cref{th:group-schemes} this is carried out for group schemes of order $p$ and $\F_q$-vector space schemes of rank one (where $q$ is the common cardinality of the residue fields of anabelomorphic $p$-adic fields).

\subsection{Local anabelomorphy as Galois Theoretic Surgery on Number Fields}\label{ss:intro1.6}
The validity of Grothendieck's Anabelian Conjecture for number fields means that a number field $M$ is anabelomorphically rigid  (\Cref{th:first-fun-anab}).  So the question of incorporating local changes of arithmetic into global arithmetic geometry  is quite a subtle one. But examples of such local (anabelomorphic) changes  already occur in many results  related to automorphic forms and Galois representations--for example \cite{taylor02}.

This idea leads to  the notion of \textit{anabelomorphically connected number fields} in which two number fields have isomorphic local (absolute) Galois groups at respective finite sets of primes (see \Cref{def:ef:anab-connected} and \Cref{ex:basic-ex2}).  This may be thought of as  \emph{Galois-theoretic surgery on number fields}; and one is interested in transferring, via local anabelomorphy at the relevant finite set of primes, objects of arithmetic and geometric interest (say automorphic forms, algebraic varieties) from one number field to the other. A basic existence for such number fields is (\Cref{th:anab-connectivity-thm2}).
\newcommand{\Out}{{\rm Out}}

Grothendieck's Section Conjecture, suggests a conjectural anabelomorphic version of Moret-Bailly's Theorem about density of global points in $p$-adic topologies for anabelomorphically connected number fields (\Cref{th:general-density-thm}). In the simplest cases, one can establish this unconditionally in \Cref{th:anabelomorphic-density} (for $\P^1-\{0,1,\infty\}$), \Cref{cor:anab-moret-bailly-proj} for projective and affine spaces. A trivial arithmetic application of \Cref{th:anabelomorphic-density} is \Cref{th:anabelomorphic-connectivty-theorem-elliptic}, but more sophisticated applications should exist (see \Cref{qu:deform-galois}).

\subsection{Weak anabelomorphy}
Since we have suggested that anabelomorphy should be roughly understood as providing an anabelian way of base-change, it is interesting to study the behavior of (say) a variety over $\Q_p$ when viewed over two anabelomorphic extensions of $\Q_p$. This leads to the notion of weak anabelomorphy (\Cref{def:weakly-anab}) and weak amphoricity (\Cref{def:weakly-amphoric}). In Theorem~\ref{th:kodaira-sym-unamphoric}(3), we show that for an elliptic curve $E$ over a $p$-adic field, all the four quantities:  the exponent of the discriminant, the exponent of the conductor,  the Kodaira Symbol and the Tamagawa Number are not weakly amphoric. 
This phenomenon arises due to the fact that elliptic curves (and curves of higher genus) may acquire potentially good reduction over wildly ramified extensions and is indicated by the presence of the Swan conductor (i.e. presence of wild ramification) in the Grothendieck-Ogg-Shafarevich formula  (also see \Cref{th:artin-swan-unamphoric},  \Cref{th:discriminant-is-unamphoric2}). 

\subsection{Anabelomorphy in perfectoid spaces, the archimedean case, and differential equations} In  \ssep\ref{se:perfectoids}, we show that anabelomorphy also appears non-trivially  in the theory of perfectoid fields (\Cref{th:perfect-anab}) and  perfectoid spaces (\Cref{th:perfect-anab3}) considered in \cite{scholze12-perfectoid-ihes}.

Anabelomorphy in the archimedean case is treated in \ssep\ref{se:anab-arch}. Our ideas were shaped by results in Hodge Theory, but the definition (\Cref{def:anab-arch}) of anabelomorphy in this case  is broader than what is considered in literature (for example , \cite{hain1987}, \cite{deligne-zeta-values}) and allows one to work with quasi-conformal mappings and their analog in higher dimensions. \Cref{th:anab-arch-diff-eq} deals with gluing differential equations by their monodromy. \Cref{pr:hodge1,pr:hodge2} deal with unipotent variations of mixed Hodge structures and should be thought of as archimedean analogs of \Cref{th:ordinary-amphoric}. This leads to \Cref{th:hodge3}, which shows that there is a natural equivalence between the categories of unipotent mixed Hodge structures on  anabelomorphic complex, quasi-projective varieties (which glues differential equations by their monodromy) and that such varieties have naturally quasi-equivalent categories of commutative, differential graded $\Q$-algebras (this last assertion shows that our definition of anabelomorphic complex quasi-projective varieties is better suited in the archimedean case).

These archimedean results and \Cref{qu:basic-langlands} motivate the following results on $p$-adic differential equations. In \ssep\ref{se:anab-diff-eq}, one shows that rank one $p$-adic differential equations  (in the sense of \cite{andre-book}) on a geometrically connected, smooth, quasi-projective and anabelomorphic  varieties can also be synchronized under anabelomorphy. One expects the higher rank case of this result to hold (\Cref{conj:p-adic-diff-eq}). \Cref{conj:irreg-unamphoric} is the natural analog of \Cref{th:artin-swan-unamphoric}.

\subsection{A picturesque way of thinking about anabelomorphy}\label{ss:picturesque-anab}
One could think of anabelomorphy in the following way:

\indent One has two parallel universes (in the sense of physics) of geometry/arithmetic over $p$-adic fields $K$ and $L$ respectively. If $K,L$ are anabelomorphic (i.e. $K\anabelmap L$) then there is a worm-hole or a conduit through which one can funnel arithmetic/geometric information in the $K$-universe to the $L$-universe through the choice of an isomorphism of Galois groups $G_K\isom G_L$, which serves as a wormhole. Information is transferred by means of amphoric quantities, properties and algebraic structures. The $K$ and $L$ universes themselves follow different laws (of algebra) as addition and multiplication has different meaning in the two anabelomorphic fields $K,L$ (in general). As one might expect, some information appears unscathed on the other side, while some is altered by its passage through the wormhole. [This is the amphoric/not amphoric dichotomy.] Readers will find  ample evidence of such phenomena throughout this paper.

\subsection{Acknowledgments}  I met  Jean-Marc Fontaine in  1994--1995 at the Tata Institute (Mumbai) where he taught a course on $p$-adic Hodge theory. I was fortunate enough to learn $p$-adic Hodge theory  directly from him. In the coming years, Fontaine arranged my stays in Paris (1996, 1997, and 2003) which provided me an opportunity to further my understanding of $p$-adic Hodge Theory from him. Influence of Fontaine's ideas on this paper and my work on Arithmetic Teichm\"uller Spaces detailed in \cite{joshi-teich,joshi-teich-estimates,joshi-teich-def,joshi-teich-rosetta,joshi-teich-abc}  should be obvious. I dedicate this paper to the memory of Jean-Marc Fontaine. 

Some of the reflections recorded herein began during my sabbatical stay at RIMS (Kyoto, Spring 2018). Support and hospitality from RIMS (Kyoto) is gratefully acknowledged. I thank Shinichi Mochizuki for his invitation and for many conversations around his work.

I thank: Yuichiro Hoshi for some conversations about \cite{hoshi-mono};  Yu Yang for some conversations; Machiel van Frankenhuijsen for many conversations.   For providing comments to the 2020 version of this paper I thank: Taylor Dupuy, Anton Hilado, Tim Holzschuh, S. Mochizuki, P. Scholze.

I thank Dinesh Thakur for encouragement and conversations, spanning many years, surrounding this paper and my work on related topics.

I take this opportunity to express my sincere gratitude to all the referees for a careful reading, several corrections, and  suggesting many improvements, which have vastly increased the readability of this paper.  I thank the editors for their kind consideration and their patience. 

\section{Anabelomorphy, amphoric quantities and amphoras}\label{se:anabelomorphy}\nwss
Let $p$ be a fixed prime number and write $\ell$ for an arbitrary prime number not equal to $p$. By a $p$-adic field, we mean a finite extension of $\Q_p$. Let $K$ be a field and let $X/K$ be a geometrically connected, smooth quasi-projective variety over $K$. The case $X=\Spec(K)$ is perfectly reasonable for understanding the definitions given below. By and large, we will assume that $K$ is either a $p$-adic field or a number field, but the ideas presented here can be used in wider contexts.

For a field $K$, let $\bK$ be a separable closure of $K$ (\textcolor{red}{note} the conflation of standard notation $K^{sep}$ and $\bK$), $G_K=\gal(\bK/K)$ be its absolute Galois group considered as a topological group. If $K$ is a $p$-adic field, let $I_K\subset G_K$ (resp. $P_K\subset G_K$) be the inertia (resp. wild inertia) subgroup of $G_K$. 

\subsection{Definitions}\label{ss:anab-defs} 
\newcommand{\sP}{\mathscr{P}}
\bdefn\label{def:anabelomorphy-fields}\ 
Let $K,L$ be two $p$-adic fields or two number fields.
\benumlab
\item  We  say that  $K,L$ are \emph{anabelomorphic} or \emph{anabelomorphs} (or anabelomorphs of each other) if one has a topological isomorphism $G_K\isom G_L$ of their absolute Galois groups.  We write $K\anabelmap L$ if $K,L$ are anabelomorphic and $\alpha:K\anabelmap L$ will mean a specific topological isomorphism $\alpha:G_K\mapright{\isom} G_L$ of topological groups. 
\item Obviously, if $L\anabelmap L'$ and $L'\anabelmap L''$, then $L\anabelmap L''$. \emph{So anabelomorphism (or anabelomorphy) is an equivalence relation on $p$-adic fields or number fields.}
\item The collection of all fields $L$ which are anabelomorphic to $K$ will be called the \emph{anabelomorphism class} of $K$. 
\item We say that \emph{$K$ is strictly anabelomorphic to $L$} or that $K\anabelmap L$ is a \emph{strict anabelomorphism} if $K\anabelmap L$ but $K$ is not isomorphic to $L$ (note that any abstract isomorphism between two $p$-adic fields $L,K$ is an isomorphism of valued fields).
\eenum
\edefn

\brem
These definitions extend  to include a broader class of fields. For example, the extension to perfectoid fields is discussed in \ssep\ref{se:perfectoids}.  For extension to quasi-projective varieties, see \ssep\ref{ss:anab-varieties}, for extension to higher dimensional local fields see Questions~\ref{qu:hdim1}--\ref{qu:hdim7}.
\erem

The phrase \textit{structure} (or an \textit{algebraic structure}) in the next definition generally refers to \cite[Chapitre I]{bourbaki_algebra1_4} and its variations (for example a group, a topological space etc. are structures). To keep the discussion precise,  an algebraic structure $A_K$, associated to a $p$-adic  field or a number field $K$,  one means a structure in this sense. However, readers should beware that in the context of anabelomorphy, there may be no homomorphisms between the fields, but a choice of an anabelomorphism between the fields can lead to isomorphisms between the associated algebraic structures (for examples see \Cref{th:third-fun-anab}, \Cref{le:basic-example}).   

\bdefn\label{def:amphoric-field}
Let $K,L$ be two $p$-adic fields or two number fields. A quantity $Q_K$ or an algebraic structure $A_K$ or a property $\sP$ of $K$ is said to be an \emph{amphoric quantity} (resp. \emph{amphoric algebraic structure}, \emph{amphoric property}) if this quantity (resp. algebraic structure or property) depends only on the anabelomorphism class of $K$. More precisely, if $\alpha:K\anabelmap L$ is an anabelomorphism of $p$-adic fields or number fields, then $Q_K=Q_L$, and one has an isomorphism $\alpha:A_K\isom A_L$ of algebraic structures which is induced by the anabelomorphism $\alpha:G_K\isom G_L$ and which is functorial in the anabelomorphism $\alpha$ (note the abuse of notation); and $L$  has property $\sP$ if and only if $K$ has property $\sP$. The \emph{amphora} of $G_K$ is the collection of all quantities, properties, algebraic structures associated with $K$ which are amphoric. 
\edefn

In \ssep\ref{se:five-fundamental}, especially \Cref{th:third-fun-anab}, the reader will find examples illustrating the non-triviality of these definitions. The next remark helps clarify the functorial aspect of the above definition.

\brem\label{re:topos-of-gk}
Let $K$ be a $p$-adic field, $\bK$ an algebraic closure of $K$ and $G_K$ be the absolute Galois group of $K$ computed using $\bK$. Then there is  a natural Grothendieck topos associated to the topological group $G_K$ \cite[Chapter III, \ssep 9, Theorem 1]{macLane_moerdijk1992} and anabelomorphic $p$-adic fields give rise to isomorphic topoi and a specific anabelomorphism gives rise to a specific isomorphism of topoi. An amphoric quantity should be considered to be an invariant of this topos, amphoric properties to be properties of the topos. Further, one expects amphoric algebraic structures to be sheaves in this topos. For example, the rule which to an open subgroup $H\subset G_K$, assigns the abelian group $K_H^*$, where $K_H=\bK^H$ is the fixed field of $H$, is expected to be a sheaf of abelian groups in this topos. These assertions will be taken up in detail in a separate paper. 
\erem
\bdefn\label{def:anabelian-rigidity-fld}
A field $K$ is \textit{anabelomorphically rigid} if, whenever one has an anabelomorphism $K\anabelmap L$ (with both $K,L$ being  $p$-adic fields or  number fields), one has an isomorphism of fields $K\isom L$.
\edefn
\brem 
Recall that any field  isomorphism $K\isom L$ of $p$-adic fields is an isomorphism of valued fields.
\erem
\subsection{Anabelomorphy of quasi-projective varieties}\label{ss:anab-varieties}\nwss
The definition of anabelomorphy of fields readily extends to smooth varieties of higher dimension. If $X/K$ is a geometrically connected, smooth quasi-projective variety over $K$, then write $\Pi_{X}$ for the \'etale  fundamental group of $X/K$ (computed for a suitable choice of a geometric base-point). If $X=\Spec(K)$, then this group coincides with $G_K$.

\bdefn\label{def:anabelomorphy-var}\ 
\benumlab
\item If $X/K$ and $Y/L$ are two geometrically connected, smooth, quasi-projective varieties, then \textit{$X/K$ is anabelomorphic to $Y/L$} (denoted $X/K\anabelmap Y/L$) if one has a topological isomorphism of the \'etale fundamental groups  
$$\Pi_{X} \isom \Pi_{Y}.$$ 
(Note that this isomorphism is not required to be compatible with the passage to the quotient $G_K$ (resp. $G_L$) on either side.) Evidently isomorphic varieties are anabelomorphic. 
\item We write $X/K\anabelmap Y/L$ if $X/K,Y/L$ are anabelomorphic, and write $\alpha:X/K\anabelmap Y/L$ if one is given  specific  isomorphism $$\alpha:\Pi_{X} \mapright{\isom}  \Pi_{Y}$$ of topological groups.  
\item An anabelomorphism $X/K\anabelmap Y/L$ is a \emph{strict anabelomorphism} or that $X/K, Y/L$ are \emph{strictly anabelomorphic} if $X/K\anabelmap Y/L$, but $X$ and $Y$ are not isomorphic as $\Z$-schemes. 
\item Anabelomorphy is an equivalence relation: if $X/L\anabelmap X'/L'$ and $X'/L'\anabelmap X''/L''$ then $X/L\anabelmap X''/L''$.  Hence one can speak of the anabelomorphism class of $X/K$.
\item If $X=\Spec(K)$ and $Y=\Spec(L)$, then $X/K\anabelmap Y/L$, if  $K\anabelmap L$   i.e.  if their absolute Galois groups are topologically isomorphic $$G_K\isom G_L.$$ 
\eenum 
\edefn

\brem\ 
\benumlab
\item The hypothesis on $X/K$ in \Cref{def:anabelomorphy-var}{\bf(1)} imply that $K$ is the integral closure of $\Q$ in $\Gamma(X,\O_X)$ in the number field case (resp. $\Q_p$ in $\Gamma(X,\O_X)$ in case $K$ is a $p$-adic field). Hence, the $K$-scheme structure  can be recovered from the absolute scheme structure of $X$. Secondly, strictly anabelomorphic varieties exist (\Cref{re:anab-Pn}).
\item One may also extend \Cref{def:anabelomorphy-var}  to other fundamental group functors. For example, one may define 
\textit{`tame anabelomorphy'} (resp. \textit{`tempered anabelomorphy'}) using the tame fundamental group (resp. the tempered fundamental group) and so on. These variants will not be used in this paper, but the tempered variant is used in \cite{joshi-teich} (and its sequels) and in \iut. Birational anabelomorphy, using the absolute Galois groups of the function fields $K(X)$ of $X/K$ and $L(Y)$ of $Y/L$ respectively, appeared in the classic works \cite{uchida-function}, \cite{pop94}.
\eenum
\erem

The following is fundamental in understanding anabelomorphy of varieties:
\bpro\label{pr:anabelomorphy-of-base} Suppose $K,L$ are finite fields, $p$-adic fields or number fields.  Any anabelomorphism $X/K\anabmapright{\alpha} Y/L$ of geometrically connected, smooth, quasi-projective varieties induces an anabelomorphism $$K\anabmapright{\alpha} L.$$
\epro
\bp 
This is \cite[Corollary 2.8(ii)]{mochizuki-topics1}.
\ep

\brem\label{re:anab-Pn}
If $K,L$ are number fields, then the above proposition together with \Cref{th:first-fun-anab} implies that $K\isom L$. Note that if $K\anabelmap L$ are anabelomorphic $p$-adic fields, then $\P^n/K\anabelmap \P^n/L$. This is a strict anabelomorphism in general (see \Cref{th:anab-toric}). 
\erem

\bdefn\label{def:anabelian-rigidity-var}
Let $X/K$ be a geometrically connected, smooth, quasi-projective variety over a field $K$. We say that $X/K$ is \textit{anabelomorphically rigid} if any anabelomorphism $\alpha:X/K\anabelmap Y/L$ (with $Y/L$ of the same sort as $X/K$), one has an isomorphism of $\Z$-schemes $X\isom Y$.
\edefn

\bdefn\label{def:amphoric-var}
Let $X/K$ be a geometrically connected, smooth, quasi-projective variety over a $p$-adic field $K$. A quantity $Q_{X/K}$ or an algebraic structure $A_{X/K}$ or a property of $\sP_{X/K}$ associated to $X/K$ is said to be an \emph{amphoric quantity} (resp. \emph{amphoric algebraic structure}  (with functoriality in the sense of \Cref{def:amphoric-field}), \emph{amphoric property}) if this quantity (resp. algebraic structure or property) depends only on the anabelomorphism class of $X/K$ i.e. it depends only on the isomorphism class of the topological group  $\Pi_{X}$. More precisely: if $\alpha:\Pi_{X}\isom \Pi_{Y}$ is an isomorphism of topological groups, then $\alpha$ takes the quantity $Q_{X/K}$ (resp. algebraic structure $A_{X/K}$, property $\sP_{X/K}$) for $X/K$ to the corresponding quantity (resp. algebraic structure, property) of $Y/L$ (with functoriality for algebraic structures as given in \Cref{def:amphoric-field}). Otherwise, a quantity (resp. algebraic structure, property) of $X/K$ will not be an amphoric quantity (resp. algebraic structure, property).
\edefn

\Cref{re:topos-of-gk} is also relevant in the context of the above definition. For examples of amphoric quantities which have been known prior to this paper, see \ssep\ref{se:five-fundamental}.

\section{Classical theorems of anabelomorphy of number fields and $p$-adic fields}\label{se:five-fundamental}\nwss
For the reader's convenience, here are the five fundamental theorems of anabelian geometry upon which anabelomorphy of $p$-adic fields and number fields rests. This list is organized logically (as opposed to a chronologically). 
\subsection{First fundamental theorem of anabelomorphy}
\bthm[First fundamental theorem of anabelomorphy]\label{th:first-fun-anab} Number fields are anabelomorphically rigid i.e.
if $K,L$ are number fields then $K$ is anabelomorphic to $L$ if and only if $K$ is isomorphic to $L$ i.e. 
$$K\anabelmap L \iff K\isom L.$$
\ethm
\bp This is a classical result due to Neukirch and Uchida \cite[Theorem 12.2.1]{neukirch-coh-book}. \ep

\subsection{The Grothendieck-Mochizuki-Tamagawa Theorem} 
For completeness, we provide the following reformulation of the celebrated theorem  conjectured by A. Grothendieck (this conjecture is also known as the Absolute Grothendieck Conjecture for smooth hyperbolic curves):\bthm[Grothendieck-Mochizuki-Tamagawa Theorem] 
Let $L,L'$ be number fields and suppose that $X/L$ (resp. $Y/L'$) are geometrically connected, smooth, hyperbolic curves over $L$ (resp. $L'$). Then the following assertions are equivalent:
\benumlab
\item There exists an anabelomorphism $X/L\anabelmap Y/L'$ of schemes.
\item There exists an isomorphism $L\isom L'$ of fields and  an isomorphism $X\isom Y$ of $\Z$-schemes. 
\eenum
In particular, geometrically connected, smooth hyperbolic curves over number fields are anabelomorphically rigid.
\ethm
\bp 
As $X/L,Y/L'$ are geometrically connected,  the integral closure of $\Q\subseteq \Gamma(X,\O_X)$ (resp. $\Q\subseteq\Gamma(Y,\O_Y)$) is $L$ (resp. $L'$). Thus, if {\bf(2)} holds, then $X$ and $Y$ are isomorphic as $L$-schemes and hence  $X/L\anabelmap Y/L'$ holds. Hence {\bf(2)}$\implies${\bf(1)}. So the non-trivial part of the assertion is to prove that {\bf(1)}$\implies${\bf(2)}. The assertion $L\isom L'$ is immediate from the hypothesis of {\bf(1)},   \Cref{pr:anabelomorphy-of-base} and \Cref{th:first-fun-anab}. The isomorphism $X\isom Y$ of $\Z$-schemes follows from {\bf(1)} by \cite[Corollary 1.3.5]{mochizuki04}. The last assertion is immediate from the equivalence of {\bf(1), (2)} and \Cref{def:anabelian-rigidity-var}.
\ep

\brem 
The history of this theorem is as follows: the special case of punctured projective lines was established by \cite{nakamura90}; the affine (absolute) case was proved in \cite[Theorem 0.4]{tamagawa97-gconj}, the proper (but relative) case established in \cite{mochizuki96-gconj}; the above formulation refers to the absolute version \cite[Corollary 1.3.5]{mochizuki04}.
\erem

\subsection{Second fundamental theorem of anabelomorphy}
\bthm[Second fundamental theorem of anabelomorphy]\label{th:second-fun-anab}
If $K,L$ are $p$-adic fields then $K\isom L$ if and only if there is a topological isomorphism of their Galois groups equipped with the respective (upper numbering) inertia filtration i.e. $(G_K,G_K^\mydot)\isom (G_L,G_L^\mydot)$
\ethm
\bp 
This is the main theorem of \cite{mochizuki-local-gro}.
\ep
\brem\label{re:standin}
Thus, one sees from \Cref{th:second-fun-anab} that the upper numbering ramification filtration is a Galois-theoretic stand-in for the field structure.
\erem

\subsection{Third fundamental theorem of anabelomorphy}
The following theorem is a combination of many different results proved by (Neukirch, Uchida, Jarden-Ritter, Mochizuki) in different time periods.

\bthm[Third fundamental theorem of anabelomorphy]\label{th:third-fun-anab} 
Let $K$ be a $p$-adic field. Then
\benumlab
\item The residue characteristic $p$  of $K$ is amphoric.
\item The degree $[K:\Q_p]$ and  the absolute ramification index $e_K$ are amphoric.
\item The topological groups $K^*$, $\O_K^*$ and $(K,+)$ (viewed as topological groups) are amphoric.
\item The $p$-adic cyclotomic character $\chi_p:G_K\to \Z_p^*$ is amphoric.
\eenum
\ethm

\bp 
For proofs of {\bf(1), (2)} see \cite[Proposition 3.6]{hoshi-mono}; for the proof of {\bf(3)} see \cite[Proposition 3.11 and Lemma 3.12]{hoshi-mono}. For the proof of {\bf(4)}, see \cite[Proposition 1.1]{mochizuki-local-gro}. 
\ep

\brem
See \cite[Summary 3.15]{hoshi-mono} for a longer list of amphoric quantities, properties and algebraic structures.
\erem
\subsection{Fourth fundamental theorem of anabelomorphy}
The next assertion was first proved in \cite{jarden79} with the additional hypothesis $\zeta_p\in K$ ($\zeta_4\in K$ if $p=2$). The hypothesis  $\zeta_p\in K$ was removed for $p\neq2$ in \cite{ritter1978}. This important theorem provides a way of deciding if two $p$-adic fields are anabelomorphic or not.
\bthm[Fourth fundamental theorem of anabelomorphy]\label{th:fourth-fun-anab} 
Let $K,L$ be  $p$-adic fields with  $\zeta_4\in K$ if $p=2$. Write $K\supseteq K^0\supseteq \Q_p$ (resp. $L\supseteq L^0\supseteq\Q_p$) be the maximal abelian subfield contained in $K$ (resp. $L$).
Then the following are equivalent:
\benumlab
\item $K\anabelmap L$
\item $[K:\Q_p]=[L:\Q_p]$ and $K^0=L^0$.
\eenum
\ethm
\bp 
See \cite{jarden79}, \cite{ritter1978}. 
\ep
\subsection{Fifth fundamental theorem of anabelomorphy}
\bthm[Fifth Fundamental theorem of anabelomorphy]\label{th:fifth-fun-anab}
Let $K$ be a $p$-adic field and let $I_K\subseteq G_K$ (resp. $P_K\subseteq G_K$) be the inertia subgroup (resp. the wild inertia subgroup). Then $I_K$ and $P_K$ are amphoric. 
\ethm
\bp
For proofs see \cite[Prop. 1.2 and the proof of Corollary 1.3]{mochizuki-local-gro}  or \cite[Proposition 3.6]{hoshi-mono}.
\ep

\newcommand{\sF}{\mathscr{F}}
\newcommand{\fT}{\mathfrak{T}_n}
These are five fundamental theorems of classical anabelomorphy. 
\brem\label{re:finiteness}
Let $p$ be a prime, let $\bQ_p$ be an algebraic closure of $\Q_p$ and let $N\geq 1$ be a positive integer. Let $$\sF_N=\left\{K:K\subset\bQ_p \text{ and } [K:\Q_p]= N\right\}.$$ Since the degrees of $K\in\sF_N$ are fixed, the set $\sF_N$ is finite. Since anabelomorphism is an equivalence relation on $\sF_N$, one obtains a   partition $\sF_N$ into a finite, disjoint union of anabelomorphism classes (each of which is also finite). Each of these partitions is the zero-dimensional analog of the classical Teichm\"uller space (\Cref{def:zero-dim-teich}). This leads to the question: is there a ``mass formula'' analogous to that of \cite{serre-mass-formula} for each of these partitions?
\erem

\newcommand{\kr}{K_{r}}
\newcommand{\fr}{F_r}
\newcommand{\lr}{L_r}
\newcommand{\zpr}{\zeta_{p^r}}
\newcommand{\dkr}{v_p(\delta(\kr/\Q_p))}
\newcommand{\dlr}{v_p(\delta(\lr/\Q_p))}

\subsection{Monoradicality of $p$-adic fields is Amphoric}\label{se:monorad}
Let $K$ be a $p$-adic field. An extension $M/K$ is a \emph{monoradical extension} if $[M:K]=m$ and $L$ is of the  form $M=K(\sqrt[m]{x})$ for some $x\in K$ (and in this case $\sqrt[m]{x}$ is called a generator of $M/K$). The following assertion is taken from \cite[Lemma 2]{jarden79}.
\bthm\label{th:monoradical-ext} 
Suppose $K\anabmapright{\sigma} L$ is an anabelomorphism of $p$-adic fields and suppose that  $M=K(\sqrt[m]{x})$ is a monoradical extension of $K$. Suppose $H\subset G_K$ is the open subgroup corresponding to $M$. Then $H'=\sigma(H)\subset G_L$ has as its fixed field a monoradical extension $M'=L(\sqrt[m]{y})$ with $y\in L$ and $v_K(x)=v_L(y)$. 
\ethm

\nws
\section{Discriminant and Different of a $p$-adic field are  not amphoric}\label{se:disc-unamph}
For the definition of the \emph{different} and the \emph{discriminant} of a $p$-adic field see \cite[Chap III]{serre1979-local-fields}.  The following result is fundamental for quantitatively understanding anabelomorphy and especially understanding the observation that the upper numbering ramification filtration of the absolute Galois group of a $p$-adic field is a stand-in for the field structure i.e. of the subtle intertwining between the additive and multiplicative groups of the field (by \Cref{th:third-fun-anab} the additive and multiplicative groups of a $p$-adic field are amphoric). More examples of this are documented in \ssep\ref{ss:weak-anab-elliptic}. 
\bthm\label{th:discriminant-is-unamphoric} The different and the discriminant of a finite Galois extension $K/\Q_p$ are not amphoric. 
\ethm
\bp
By \cite[Chap III, Prop 6]{serre1979-local-fields} and \Cref{th:third-fun-anab}{\bf(2)}, it is sufficient to prove  that the  different of $K/\Q_p$ is not amphoric. 

Let $r\geq 1$ be an integer, $p$ an odd prime and let $\kr=\Q_p(\zpr,\sqrt[{p^r}]{p})$ so $\fr\subset\kr$ and let $\lr=\Q_p(\zpr,\sqrt[{p^r}]{1+p})$. By Lemma~\ref{le:basic-example} below,  one has an anabelomorphism $\kr\anabelmap \lr$ and hence one has $G_{\lr}\isom G_{\kr}$. But, $\kr$ and $\lr$ are not isomorphic fields so by \cite{mochizuki-local-gro} they have distinct inertia filtrations. We claim that they have  distinct differents and  discriminants. More precisely, one has the following formulae for the discriminants of $\kr/\Q_p$ (resp. $\lr/\Q_p$) \cite[Theorem 5.15 and 6.13]{viviani04}.
\beas
\dkr & = & rp^{2r-1}(p-1)+p\left(\frac{p^{2r}-1}{p+1}\right)-p\left(\frac{p^{2r-3}+1}{p+1}\right),\\
\dlr & = & p^r\left( r\cdot p^r-(r+1)\cdot p^{r-1}\right) + 2\left(\frac{p^{2r}-1}{p+1}\right).
\eeas
Note that as $(p+1,p)=1$, for $r\geq2$,   $\dkr$ is rational number whose numerator is clearly divisible by $p$, while $\dlr$ is a rational number whose numerator is not divisible by $p$.  For $r=1$, these are equal to $2p(p-1)-1$ and $p^2-2$ respectively and evidently $2p(p-1)-1\neq p^2-2$ for any odd prime $p$. Thus for all $r\geq1$, one has $\dkr\neq \dlr$. This proves the assertion.
\ep

\blem\label{le:basic-example}Let $r\geq 1$ be any integer and $p$ any odd prime.
Let  $\kr=\Q_p(\zpr,\sqrt[{p^r}]{p})$ and let $\lr=\Q_p(\zpr,\sqrt[{p^r}]{1+p})$. Then one has an anabelomorphism $$ \kr \anabelmap \lr \text{ equivalently } G_{\lr}\isom G_{\kr}.$$
\elem
\bp 
Let $\fr=\Q_p(\zpr)$. Both fields contain $\fr$ and by elementary Galois theory and Kummer theory one checks that $\fr\subset \kr$ and $\fr\subset \lr$ is the maximal abelian subfield of both $\kr,\lr$ and both $\kr,\lr$ have the same degree over $\Q_p$. Then \Cref{th:fourth-fun-anab} says that the absolute Galois groups of $\kr,\lr$ are isomorphic i.e. $\kr\anabelmap \lr$. Hence, the claim.
\ep

\section{Anabelomorphy and Galois representations}\label{se:anabel-galois-reps} \nwss
\subsection{Definitions}
Consider an auxiliary  topological field $E$ which will serve as a coefficient field for representations of $G_K$. The following list of  coefficient fields  will be more than adequate for the present discussion:  $E=\Q_\ell$ for any finite prime $\ell$ including $\ell=p$ and $E=\C$ if $\ell=\infty$ and occasionally $E$ will be a finite field (with discrete topology).

Let $K$ and $L$ be $p$-adic fields. Let $\reptrip$ be a representation of $G_K$. Let $\alpha: K\anabelmap L$ be an anabelomorphism. Then as $\alpha:G_K\isom G_L$, so  by composing with $\alpha^{-1}:G_L\to G_K$ any $G_K$-representation gives rise to a $G_L$-representation and conversely,  by composing with $\alpha:G_K\to G_L$ any $G_L$-representation gives rise to a $G_K$-representation. One sees immediately that this isomorphism induces an equivalence between categories of finite dimensional continuous $E$-representations of $G_K$ and $G_L$ respectively. This will be referred to as \textit{anabelomorphy of Galois representations.} 

\bdefn\label{def:rep-amphoric}
Let $V$ be a finite dimensional $E$-vector space (viewed as a topological vector space using the topology of $E$). Let $\rho:G_K\to \gl V$ be a continuous representation of $G_K$. We say that a quantity or an algebraic structure or a property of the triple $(G_K,\rho, V)$ is \textit{amphoric} if it is determined by the anabelomorphism class of $K$.  
\edefn

\begin{example}\label{le:irred-amphoric}
Let $K$ be a $p$-adic field and let $E$ be a coefficient field. Then from \Cref{def:rep-amphoric} one has
\benumlab
\item $\dim_E(V)$ is an amphoric quantity.
\item The category of finite dimensional $E$-representations of $G_K$ is amphoric.
\item Irreducibility of a $G_K$-representation is an amphoric property.
\eenum
\end{example}
\subsection{Unramifiedness and tame ramifiedness of a local Galois representation are amphoric}\label{se:unram-amph}
Recall that $\rho:G_K\to \gl V$ is \emph{unramified} (resp. \emph{tamely ramified}) if the image $\rho(I_K)=1$ (resp. $\rho(P_K)=1$).
\bthm\label{th:unram-amph}
Let $K$ be a $p$-adic local field. Unramifiedness (resp. tame ramifiedness) of $\reptrip$ are amphoric properties.
\ethm
\bp 
The assertion is immediate from the definition of unramifiedness (resp. tame ramifiedness) and \Cref{th:fifth-fun-anab}. 
\ep

\subsection{Ordinarity of a local Galois representation is amphoric}\label{se:ordinarity-unamph}
Let $K$ be a $p$-adic field and let $\rho:G_K\to \gl V$ be a continuous $E$-representation of $G_K$ with $E\supseteq\Q_\ell$ a finite extension of $\Q_\ell$ (and $\ell\neq p$). Then $(\rho,V)$ is said to be an \emph{ordinary representation of $G_K$} if the image $\rho(I_K)$ of the inertia subgroup of $G_K$ is unipotent. Let $\chi_p:G_K\to \Z_p^*$ be the $p$-adic cyclotomic character.
\newcommand{\gr}{\rm gr}
Recall from \cite{perrin-riou1994} that a $p$-adic representation $\reptrip $, where  $V$ is a finite dimensional $\Q_p$-vector space, is said to be an \emph{ordinary $p$-adic representation of $G_K$} if there exist $G_K$-stable filtration  $\{V_i\}$ on $V$ consisting  of $\Q_p$-subspaces of $V$ such that the action of $I_K$ on $\gr_i(V)$ is given by $\chi^i_p$ (as $G_K$-representations).

\bthm\label{th:ordinary-amphoric} 
Ordinarity of an $\ell$-adic or a $p$-adic representation $\reptrip$ is an amphoric property.
\ethm
\bp 
Let $\rho:G_K\to \gl V$ be a continuous, ordinary Galois representation on $G_K$ on a finite dimensional $E$ vector space with $E/\Q_\ell$ a finite extension. Let $L$ be a $p$-adic field with an isomorphism $\alpha:G_L\isom G_K$. By \Cref{th:fifth-fun-anab}, the inertia (resp. wild inertia) subgroups are amphoric. Then $\rho(\alpha(I_L))\subset\rho(I_K)$ so the image of $I_L$ is also unipotent. This gives the assertion for $\ell\neq p$. But the $\ell=p$ case is also similar. Recall from \Cref{th:third-fun-anab}{\bf(4)} that for any $p$-adic field $K$, the $p$-adic cyclotomic character of $G_K$ is amphoric. From \Cref{le:irred-amphoric}, the filtration $V_i$ is $G_L$-stable and from \Cref{th:third-fun-anab}, $\chi_p\circ\alpha$ is the cyclotomic character of $G_L$.  By definition, for any $v\in V_i$ and any $g\in I_K$, $$\rho(g)(v)=\chi_p^i(g) v + V_{i+1}.$$ Now given an isomorphism $\alpha:G_L\to G_K$, one has  for all $g\in G_L$
$$\rho(\alpha(g))(v)=\chi_p^i(\alpha(g)) v + V_{i+1}.$$
and thus ordinarity of $\reptrip$ is  determined solely by the isomorphism class of $G_K$.
\ep

\newcommand{\knr}{K^{\rm nr}}
\newcommand{\lnr}{L^{\rm nr}}
\subsection{Anabelomorphy of $K^t$ and $K^{nr}$}
For a $p$-adic field $K$, let $\bK$ be an algebraic closure of $K$. Let  $\knr$ (resp. $K^t$)  be the maximal unramified extension (resp. maximal tamely ramified extension) of $K$ contained in $\bK$.
\bpro\label{le:anab-knr}
Suppose  $K\anabmapright{\alpha} L$ is an anabelomorphism of $p$-adic fields. Then one has anabelomorphisms induced by $\alpha$:
$$\alpha: K^t\anabelmap L^t,$$ 
and
$$\alpha:\knr\anabelmap \lnr.$$
Moreover, these anabelomorphisms are preserved under the passage to the completions of these fields.
\epro
\bp 
By \Cref{th:fifth-fun-anab}, the inertia subgroup $I_K$ (resp.  the wild inertia subgroup $P_K$) of $G_K$ are amphoric. As $\knr$ (resp. $K^t$) is the fixed field of $I_K$ (resp. $P_K$), both the anabelomorphisms follow from $K\anabelmap L$. Since the Galois group is unaffected by passage from a rank-1 henselian valued field to its completion, the last assertion is also clear.
\ep

\newcommand{\fq}{{\mathbb{F}_q}}
\newcommand{\fqb}{{\bar{\mathbb{F}}_q}}
\subsection{Peu and Tres ramifiedness are not amphoric properties}\label{se:peu-tres-unamp}
In many theorems in the theory of Galois representations and modular forms, the notion of \textit{peu and tres ramifi\'ee extensions}, defined in \cite[Section 2.4]{Serre1987} plays an important role (for example \cite{edixhoven1992}). 

\bthm\label{th:peu-tres} 
The property of an extension $L/\Q_p$ being peu ramifi\'ee (resp. being tres ramifi\'ee) extension (resp. representation) is not amphoric.
\ethm
\bp 
By  Lemma~\ref{le:basic-example}, the fields $K=\Q_p(\zeta_p,\sqrt[p]{p})$ and $L=\Q_p(\zeta_p,\sqrt[p]{1+p})$ are strictly anabelomorphic. We claim that  $K/\Q_p$ is a tres ramifi\'ee extension while $L/\Q_p$ is a peu ramifi\'ee extension. Observe that $K/\Q_p$ and $\Q_p^{nr}/\Q_p$ are linearly disjoint over $\Q_p$ as $K/\Q_p$ is totally ramified while $\Q_p^{nr}/\Q_p$ is unramified and a similar assertion holds for $L/\Q_p$. Thus, one has extensions
$$K\Q_p^{nr}\supset \Q_p(\zeta_p)\Q_p^{nr}\supset \Q_p^{nr}$$
with the first inclusion being a totally ramified extension while the second inclusion giving tamely ramified extension.
A similar assertion holds for $L$. Then it is immediate from the definitions \cite[Section 2.4]{Serre1987} that $K\Q_p^{nr}/\Q_p^{nr}$ is tres ramifi\'ee while $L\Q_p^{nr}/\Q_p^{nr}$ is peu ramifi\'ee.\ep

\brem 
\Cref{th:peu-tres} and \cite{Serre1987} suggest that anabelomorphy affects deformation theory of Galois representations. A precise assertion is given in \Cref{th:anab-deform-thry}
\erem

\subsection{Frobenius elements are Amphoric}\label{se:frob-anab}
One has the following result of Uchida from \cite[Lemma 3]{jarden79}:
\bthm\label{th:frob-amphoric} 
Let $K \anabelmap L$ be an anabelomorphism of $p$-adic fields. If $\sigma\in G_K$ is a Frobenius element for $K$, then for any topological isomorphism $\alpha:G_K\mapright{\isom} G_L$, $\alpha(\sigma)$ is a Frobenius element for $L$.
\ethm
This has the following important corollary.
\bcor\label{th:char-poly-amphoric}
Let $K$ be a $p$-adic field and let $\reptrip$ be a finite dimensional continuous representation of $G_K$ in an $E$-vector space with $E/\Q_\ell$ a finite extension and $\ell\neq p$. Then the (local) $L$-function $L(\rho,V,T)$ is amphoric. 
\ecor
\bp
By definition, $L(\rho,V,T)=\det(1-T\cdot\rho(Frob_p)\big|_{V^{I_K}})$, hence the assertion follows from \Cref{th:frob-amphoric}.
\ep

\subsection{Amphoricity of the Iwasawa cohomology}
\bpro\label{pr:iwasawa-coh} 
Let $K$ be a $p$-adic field and let $V$ be a continuous $\Z_p$-representation of $G_K$ given by the cyclotomic character.  Then the Iwasawa cohomology $H^i_{Iw}(G_K,V)$ is amphoric for all $i\geq 0$.
\epro
\bp
Let $L\anabmapright{\alpha} K$ be an anabelomorphism of $p$-adic fields. Fix an algebraic closure $\bK$ of $K$ (resp. $\bL$ of $L$).  By \cite[D\'efinition II.1.1]{colmez1999} one has
$$ H^i_{Iw}(G_K,V) = \invlim_{n\geq 1} H^i(G_{K_n},V),$$
where for $n\geq 1$, $K_n=K(\zeta_{p^n})$ where $\zeta_{p^n}\in \bK$ is a primitive $p^n$-root of unity in $\bK$ and where the inverse image is with respect to corestriction maps.
For $n\geq 1$, write $L_n=L(\zeta_{p^n})$ where $\zeta_{p^n}\in \bL$ is a primitive $p^n$-root of unity in $\bL$, By \Cref{th:third-fun-anab}, the $p$-adic cyclotomic character is amphoric and hence,  for all $n\geq 1$, $\alpha$ induces anabelomorphisms $L_n\anabmapright{\alpha_n} K_n$ which are compatible with the cyclotomic action of $G_L$ and $G_K$ respectively. Let $V'$ be the $G_L$-representation obtained from the $G_K$-representation $V$. Hence, one obtains for  all $n\geq1$,  and for all $i\geq 0$, an isomorphism
$$ H^i(G_{K_n},V) \overset{\alpha_n}{\isom}  H^i(G_{L_n},V').$$
These isomorphisms are compatible with corestriction maps on both the sides as the corestriction map is obtained from  the transfer map which is functorial the pairs of open groups $G_{K_{n+1}}\subset G_{K_n}$ (and similar ones for $L$). Hence, passage to inverse limits gives the amphoricity assertion.
\ep

\subsection{Artin and Swan Conductor of a local Galois representation are not amphoric}\label{se:swan}
For consequences of this in the context of elliptic curves and curves in general see \ssep\ref{se:weak-anabelomorphy}. For Artin and Swan conductors see \cite{serre1979-local-fields}, \cite[Chapter 1]{katz-gauss}. The coefficient field of our $G_K$ representations will be a finite extension $E/\Q_\ell$ with $\ell\neq p$. The Artin conductor (resp. the Swan conductor) of an unramified (resp. tamely ramified) representation are zero.  The theorem is the following:
\bthm\label{th:artin-swan-unamphoric} 
Let $K$ be a $p$-adic field and let $\reptrip$ be $\Q_\ell$-adic representation of $G_K$ with finite image. Let   $Artin(\rho,V)$ (resp. $Swan(\rho,V)$) be the Artin conductor (resp. the Swan conductor) of $\reptrip$. Then
\benumlab
\item The property $\textrm{Artin}(\rho,V)=0$ and $\textrm{Swan}(\rho,V)=0$ are amphoric properties.
\item The Artin conductor  and the Swan conductor are not amphoric in general.
\eenum
\ethm
\brem\label{re:ramification-manifest-additive-str} 
The Artin and the Swan conductors depend on the ramification filtration which, by \Cref{re:standin} is a stand-in for the field structure.\erem
\bp 
Let $L\anabmapright{\alpha} K$ be an anabelomorphism. Suppose $\rho$ has finite image, then so does $\rho'=\rho\circ\alpha$. From \cite{Serre1987} one knows that $\textrm{Artin}(\rho,V)=0$ if and only if $\rho$ is unramified and $\textrm{Swan}(\rho,V)=0$ if and only if $\rho$ is tamely ramified. By \Cref{th:unram-amph}, both these properties of $\reptrip$ are amphoric. Hence, {\bf(1)} is proved.

By \cite[Chapter 19]{serre-linear-representations} or \cite[1.2]{Serre1987}, it is enough to prove that the Swan conductor is not amphoric and to prove this it is enough to give an explicit example for each prime $p$. Using \cite[Th\'eor\`eme 1.3]{henniart1988}, which provides a formula for number of Galois representations (of the sort occurring in the Local Langlands Correspondence) which have a given integer as their Swan Conductor, one sees that such representations exist. 

Explicit examples can also be constructed. Let $K_1=\Q_p(\zeta_p,\sqrt[p]{1+p})$ and $K_2=\Q_p(\zeta_p,\sqrt[p]{p})$. By \Cref{le:basic-example}, these fields are anabelomorphic and one has $\gal(K_1/\Q_p)\isom \Z/p\rtimes(\Z/p)^*\isom \gal(K_2/\Q_p)$. By the character table for this finite group (see \cite[Theorem 3.7,]{viviani04}), there is a unique irreducible character $\chi$ of a $\C$-representation $V$ of dimension $p-1$. Finiteness of these galois groups means that this representation descends to an algebraic number field (and hence provides representations with coefficient fields which are finite extensions of $\Q_{\ell}$ for all $\ell$). For $i=1,2$, let $f_i(\chi)$  denote the exponent of the Artin conductor of $\chi$. Then by \cite[Cor. 5.14 and 6.12]{viviani04} one has
\beas
f_1(\chi)&=& p\\
f_2(\chi)&=& 2p-1.
\eeas
Evidently $f_1(\chi)\neq f_2(\chi)$.
\ep

The theory of the Swan and Artin conductors   depends on the theory of break-decomposition in $\ell$-adic representations (see \cite[Chapter 1, 1.1-1.10]{katz-gauss}):
\bcor\label{co:break-decomp-unamphoric}
Let $K$ be a $p$-adic field and let $\reptrip$ be a continuous representation of $G_K$ in a $\Q_\ell$-vector space $V$ with non-trivial action of the wild inertia subgroup $P_K$.   Then the breaks in the break-decomposition  of $V$ are not amphoric in general. 
\ecor
\bp 
By \cite[Definition 1.6]{katz-gauss}, the Swan conductor is the sum, counted with multiplicity, of the breaks in the break-decomposition and by \cite[Proposition 1.1 and Definition 1.2]{katz-gauss} each break is a non-negative rational number.  By \Cref{th:artin-swan-unamphoric} (and its proof) the Swan conductor is not amphoric. This means that the breaks are not amphoric in general.
\ep

\section{Anabelomorphy and $p$-adic Hodge Theory}\label{se:p-adic-hodge}
\subsection{Crystalline-ness of a $p$-adic representation is not amphoric in general}
\Cref{th:ordinary-amphoric} should be contrasted with the following result which combines fundamental results of Mochizuki and Hoshi  \cite{mochizuki-topics1,hoshi13,hoshi18}:
\bthm\label{th:crystalline-not-amphoric}\  
\benumlab
\item Let $\alpha:K\anabelmap L$ be an anabelomorphism of $p$-adic fields.  Then the following conditions are equivalent
\benum
\item For every Hodge-Tate representation $\reptrip$, the composite $\rho\circ\alpha$ is a Hodge-Tate representation of $G_L$.
\item There exists an isomorphism of algebraic closures $\bK\isom \bar{L}$ which induces an isomorphism $K\mapright{\isom} L$.
\eenum
\item For every prime number $p$, there exist a $p$-adic field $K'$, an anabelomorphism $K'\anabmapright{\alpha}L'$ of $p$-adic fields and a crystalline representation $\reptrip$ such that $\rho\circ\alpha:G_{K'}\to \gl{V}$ is not Hodge-Tate. 
\item In particular, being  crystalline, semi-stable or de Rham  is not an amphoric property of a general $p$-adic representation.
\eenum
\ethm
\bp 
The assertion {\bf(1)} is \cite[Theorem 3.5(ii)]{mochizuki-topics1}. Let us prove ${\bf(2)}$.  Pick a $p$-adic field $K$ and a strict anabelomorphism $K\anabmapright{\alpha} L$ (by \Cref{le:basic-example}) such $K,L$ exist). For any open subgroup $H\subset G_K$, let $K'$ be the fixed field of $H$ (so $H=G_{K'}$); let $H'=\alpha(H)\subset G_L$. Then $H'$ is also open and let $L'$ be the fixed field of $H'$ so that $H'=G_{L'}$. Moreover, one has an anabelomorphism $K'\anabmapright{\alpha'} L'$ where $\alpha':H\isom H'=\alpha(H)$ is given by the restriction of $\alpha$ to $H$. This notational setup will be applied in the remainder of the proof by choosing a suitable $H$. Since the anabelomorphism $\alpha$ is strict, by \cite[Theorem 3.5(i)]{mochizuki-topics1}, there  exists a $\reptrip$ and an open subgroup $H\subset G_K$, such that the $H=G_{K'}$-representation $V'=V\big|_{H}$ is a representation given by a  Lubin-Tate character, but the $G_{L'}=H'$-representation $\rho':G_{L'}\to GL(V')$ obtained via the anabelomorphism $K'\anabmapright{\alpha'} L'$ is not Hodge-Tate.  It is standard that a representation given by a Lubin-Tate character is crystalline. Thus, the $G_{K'}$-representation $V'$ is crystalline, but the $G_{L'}$-representation $V'$ is not Hodge-Tate. This proves the assertion.
\ep

\subsection{Amphoricity of pure Hodge-Tate weight $p$-adic representations}
\newcommand{\kbh}{\C_K}
\newcommand{\lbh}{\C_L}
\newcommand{\khi}{\hat{K}_\infty}
\newcommand{\lhi}{\hat{L}_\infty}
Let $\bK\supset K$ (resp. $\bL$) be an algebraic closure of $K$ (resp. $L$), let $\C_K$ (resp. $\C_L$) be the $p$-adic completion of $\bK$ (resp. $\bL$).
We begin with a  somewhat elementary result below which is still true despite \Cref{th:crystalline-not-amphoric}.  This  is surprising because the main theorem of \cite{mochizuki-local-gro}  says that  the $p$-adic completion $\kbh$ is not amphoric in general (for example, \cite[Proposition 2.2]{mochizuki-local-gro} shows that the determination of the $G_K$-module $(\kbh,+)$ requires the ramification filtration of $G_K$). A \textit{$\kbh$-admissible representation} is a Hodge-Tate representation of $G_K$ of weight zero (see \cite[3.2]{fontaine94b}). 

\bthm\label{th:sen-thm-and-anabelomorphy}
Let $K$ be a $p$-adic field and let $\alpha:L\anabelmap K$ be an anabelomorphism. Let $\reptrip$ be a $p$-adic representation. 
\benumlab
\item Then $V$ is of Hodge-Tate weight zero (equivalently $\kbh$-admissible), if and only if, $\rho\circ\alpha$ is of Hodge-Tate weight zero (equivalently $\lbh$-admissible). 
\item In particular, $V$ is pure of Hodge-Tate weight $m$ as a $G_K$-module, if and only if, $V$ is pure of Hodge-Tate weight $m$ as a $G_L$-module.
\eenum
\ethm
\bp 
By \cite[3.2 Proposition]{fontaine94b} (and \cite[Pages 113-114]{sen80}), $V$ is $\kbh$-admissible, if and only if,  $\rho(I_K)$ is finite. By \Cref{th:fifth-fun-anab},  $\rho(I_K)$ is finite, if and only if, $\rho(\alpha(I_L))$ is finite. So the assertion {\bf(1)} is proved.

If $V$ is Hodge-Tate of weight $m$, then twisting $V$ by $\chi_p^{-m}$, one can assume that $V$ is Hodge-Tate of weight zero as a $G_K$-representation, and the assertion follows from  {\bf(1)}. This proves the assertion.
\ep

\subsection{Anabelomorphy of $(\varphi,\Gamma)$-modules and $\phisen$ is not amphoric}\label{se:phisen-unamph}
Let $K$ be a $p$-adic field, let $\bK\supset K$ be an algebraic closure of $K$. Let $\C_K$ be the $p$-adic completion of $\bK$. Let $H_K\subset G_K$ be the kernel of the composite homomorphism $G_K\mapright{\chi_K} \Z_p^*\to \Z_p^*/\textrm{Tor}(\Z_p^*)=\Z_p$, where $\chi_K$ is the  $p$-adic cyclotomic character, and  $\textrm{Tor}(\Z_p^*)\subset \Z_p$ is the torsion subgroup. Let $\Gamma=\Gamma_K=G_K/H_K$, then one has $\Gamma\isom \Z_p$. 
Let $K_\infty=\bK^{H_K}$ be the fixed field of $H_K$. Let $\khi\subset \C_K$ be the $p$-adic completion of $K_\infty$. 
In the notation of \cite{wintenberger1983}, let $X_K(K_\infty)$ (resp. $X_L(L_\infty)$ be the field of (cyclotomic) norms of $K_\infty/K$ (resp. $L_\infty/L$) and let $G_{X_K(K_\infty)}$ (resp. $G_{X_L(L_\infty)}$) be its absolute Galois group.

For a $p$-adic representation $\reptrip$ of $G_K$ let $\phisen(\rho,V)$ be the invariant defined by \cite[Theorem 4]{sen80}. These conventions will be in force in this subsection. The precise meaning of the title of the subsection is given by the following theorem.
\newcommand{\bB}{\mathbf{B}}
\newcommand{\modet}[1]{{\rm Mod}_{\bB_{#1}}(\varphi,\Gamma)}
\newcommand{\rep}[1]{{\bf Rep}_{\Q_p}(G_{#1})}
\bthm\label{th:phi-gamma}
Let $K$ be a $p$-adic field. Each anabelomorphism $L\anabmapright{\alpha} K$ of $p$-adic fields sets up
\benumlab
\item an anabelomorphism $K_\infty\anabmapright{\alpha} L_\infty$, and
\item an anabelomorphism $\hat{K}_\infty\anabmapright{\alpha} \hat{L}_\infty$ of perfectoid fields with isometric tilts $\khi^\flat\isom \lhi^\flat$, and
\item an anabelomorphism $X_K(K_\infty)\anabelmap X_L(L_\infty)$ between the cyclotomic fields of  norms of $K$ and $L$ respectively;
\item an equivalence between the category $\modet{K}$ of \'etale $(\varphi,\Gamma)$-modules over a certain field $\bB_K$  and  the category $\modet{L}$ of \'etale $(\varphi,\Gamma)$-modules over a corresponding field $\bB_L$. 
\item Under the equivalence of {\bf(4)}, the property ``$\phisen(\rho,V)$ is semisimple and has integer eigenvalues'' of a continuous $p$-adic representation $\reptrip$, is not amphoric.
\eenum
\ethm
\bp 
Let $\alpha:L\anabelmap K$ be an anabelomorphism. By \Cref{th:third-fun-anab}, the cyclotomic character is amphoric and $\chi_L=\chi_K\circ\alpha$ and one has similar quantities associated to $L$, namely $\varphi_L, \Gamma_L$, $\bL_\infty$,  $\lhi$ etc.

By the amphoricity of the cyclotomic character one has an isomorphism $H_L\overset{\alpha}{\isom} H_K$ and hence also of the quotients $G_L/H_L\isom \Z_p \isom G_K/H_K$. Hence {\bf(1)} is immediate.

By \cite[Example 2.1.1 and 2.2.2]{weinstein-aws}, the fields $\lhi$ and $\khi$ are perfectoid with absolute Galois groups $H_L$ and $H_K$ respectively. Since these groups are isomorphic, one sees that the anabelomorphism $\alpha$ induces an anabelomorphism $$\lhi\anabmapright{\alpha} \khi$$
of perfectoid fields.
Further, one also sees that the tilts $\khi^\flat\isom \lhi^\flat$ are isometric. This proves {\bf(2)}.

We claim that one has isomorphisms of topological groups:
$$G_{X_K(K_\infty)}\isom H_K\overset{\alpha}{\isom}  H_L \isom G_{X_L(L_\infty)}.$$
To see this, note that the isomorphism in the middle is given by the proof of {\bf(1)} and the outer two isomorphisms are given by \cite[Corollaire 3.2.3]{wintenberger1983}.
Thus, one has a natural anabelomorphism $X_K(K_\infty) \anabmapright{\alpha} X_L(L_\infty)$ of the fields of norms.
This proves {\bf(3)}.

Now to prove {\bf(4)}. By \cite[Part IV, Section 13.6]{conrad-brinon} (also see \cite[Th\'eor\`eme 3.4.3, 3.4.4 Remarques(c)]{fontaine2007} which uses different notation), there exists an equivalence between the category, $\rep K$, of continuous $p$-adic representations of $G_K$ and the category, $\modet{K}$, of \'etale $(\varphi,\Gamma)$-modules over a certain field $\bB_K$; a similar description holds for the category $\rep L$ of continuous $p$-adic representations of $G_L$ and the category $\modet{L}$ of \'etale $(\varphi,\Gamma)$-modules over a certain field $\bB_L$.

Now let $\reptrip$ be a continuous $p$-adic representation of $G_K$ and let $L\anabmapright{\alpha} K$ be an anabelomorphism of $p$-adic fields providing an isomorphism $G_L\mapright{\alpha} G_K$ of topological groups. Write $\rho'=\rho\circ\alpha:G_L\to \gl V$ for the $p$-adic representation of $G_L$ obtained by composition with $\alpha$. By \cite[Part IV, Page 227]{conrad-brinon}, one associates to $(\rho, V)\in \rep K$ an \'etale $(\varphi,\Gamma)$-module $M(\rho, V)\in \modet{K}$ (with certain other data which is unimportant at the moment) and this association is an equivalence of categories. The equivalence asserted by the theorem is simply the association $\modet{K}\ni M(\rho,V)\mapsto M(\rho',V)\in\modet{L}$.
This completes the proof of {\bf(4)}. 

Now to prove {\bf(5)}, consider a continuous $p$-adic representation $\reptrip$ of $G_K$. By \cite[Theorem 4]{sen80}, there exists an endomorphism $\phisen(\rho,V)\in \End((V\tensor\C_K)^{H_K})$ of the $K_\infty$-vector space $(V\tensor \C_K)^{H_K}$. By \cite[Theorem 5]{sen80}, one can always find a basis of $(V\tensor \C_K)^{H_K}$ such $\phisen$ is given by a matrix with coefficients in $K$. 

By \cite[Corollary of Theorem 6]{sen80}, the $G_K$-representation $V$ is Hodge-Tate if and only if $\phisen(\rho,V)$ is semisimple and has integer eigenvalues.  By \Cref{th:crystalline-not-amphoric}, there exists some strictly anabelomorphic pair of $p$-adic fields $L\anabmapright{\alpha} K$ and a crystalline representation $\reptrip$ such that the $G_L$-representation $(\rho',V)$ is not Hodge-Tate. Let $\phisen(\rho',V)$ be the endomorphism of the $L_\infty$-vector space $(V\tensor \C_L)^{H_L}$ arising from the $G_L$-representation $(\rho',V)$. Then $\phisen(\rho', V)$ is either not semisimple or it does not have integer eigenvalues. This completes the proof of {\bf(5)} and the theorem.
\ep

\bcor
In the notation of \Cref{th:phi-gamma} and its proof, one has a natural isomorphism (with respect to choice of anabelomorphisms $K\anabelmap L$) of topological groups 
$$ X_K(K_\infty)^* \isom  X_L(L_\infty)^*.$$
In other words, the multiplicative structures of the fields of cyclotomic norms of anabelomorphic $p$-adic fields are naturally isomorphic (with respect to choice of anabelomorphisms $K\anabelmap L$).
\ecor
\bp 
The field $K_\infty = \bigcup_{n\geq 1} K_n$ (resp. $L=\bigcup_{n\geq 1} L_n$) is an increasing union of $p$-power cyclotomic extensions $K_n$ of $K$ (resp. $L_n$ of $L$). Let  $\alpha:K\anabelmap L$ be an anabelomorphism, let $\chi_K$ (resp. $\chi_L$) be the $p$-adic cyclotomic character of $K$ (resp. $L$). Then using $\chi_L\circ\alpha=\chi_K$ one inductively defines, for each $n\geq 1$,  anabelomorphisms $K_n\anabmapright{\alpha_n} L_n$. By \cite[2.1.1]{wintenberger1983}, the multiplicative group $X_K(K_\infty)^*$ is given by
$$X_K(K_\infty)^*=\invlim_{n} K_n^*,$$
where the inverse limit is with respect to the norm homomorphisms $N_{K_{n+1}/K_n}:K_{n+1}^*\to K_n^*$. The norm homomorphism corresponds to the inclusion $G_{K_{n+1}}^{ab}\into G_{K_n}^{ab}$ of the abelianizations of $G_{K_{n+1}}$ and $G_{K_n}$ respectively and is compatible with the reciprocity homomorphism $K_{n}^*\to G_{K_n}^{ab}$ (for instance see \cite[Lemma 1.7(ii)]{hoshi-mono}) and by \cite[Proposition 3.11]{hoshi-mono} one obtains compatibility with anabelomorphisms $G_{K_n}\anabmapright{\alpha_n} G_{L_n}$. Thus, one sees that the isomorphisms $K_n^*\isom L_n^*$ are compatible with norm homomorphisms on either side and hence the assertion follows on passage to inverse limits.
\ep

\subsection{A useful lemma}
\begin{lemma}\ 
	\label{th:anab-galois-h1}
	Let $K\anabelmap L$ be two anabelomorphic $p$-adic fields. Then one has
	\benumlab
	\item for each prime number $\ell$, an isomorphism of $\Q_\ell$-vector spaces
	$$ H^1(G_K,\Q_\ell(1))\isom\Ext^1_{G_{K}}(\Q_{\ell}(0),\Q_{\ell}(1)) \isom \Ext^1_{G_{L}}(\Q_{\ell}(0),\Q_\ell(1))\isom H^1(G_L,\Q_\ell(1)),$$
	\item for $\ell=p$, an isomorphism of $\Q_p$-subspaces of the above vector spaces {\bf(1)} defined in \cite[Section 3]{bloch90}:
	\begin{enumerate}
		\item $H^1_f(G_K,\Q_p(1)) \isom H_f^1({G_L},\Q_p(1)),$
		\item $H^1_e(G_K,\Q_p(1)) \isom H_e^1({G_L},\Q_p(1)).$
	\end{enumerate}
	\eenum
\end{lemma}
\bp
Choose an anabelomorphism $\alpha: G_K\anabelmap G_L$.  Then by \cite[Proposition 4.2(iv)]{hoshi-mono} (this result is implicit in the proof of \cite[Proposition 1.1]{mochizuki-local-gro}), the Galois module  $G_K\act\hat{\Z}(1)_K$ of roots of unity in an algebraic closure of $K$ is amphoric and hence $\alpha$  carries the $\ell$-adic cyclotomic character $\chi_{L,\ell}$ to $\chi_{K,\ell}$ for each prime number $\ell$. This gives the middle isomorphism in {\bf(1)}. The outer isomorphisms are a special case of \cite[Lemme 3.3(i)]{perrin-riou1994}.

Note that {\bf(2)(b)} follows from {\bf(2)(a)} by \cite[Example 3.9]{bloch90} as $$H^1_f(G_K,\Q_p(1))\supseteq H^1_e(G_K,\Q_p(1))$$ and the two have the same dimensions as $\Q_p$-vector spaces. So it remains to prove {\bf(2)(a)}. This follows from \Cref{th:third-fun-anab}{\bf(3)} and the fact that $H^1_f(G_K,\Q_p(1))=\left(\projlim \O_K^*/\O_K^{*p^n}\right)\tensor\Q_p$ (\cite[Example 3.9]{bloch90}). 
\ep

\subsection{The $\linv$-invariant is not amphoric}\label{se:l-invariant} 
Let $K$ be a $p$-adic field and let $V$ be a 2-dimensional ordinary semi-stable representation of $G_K$ fitting in the following exact sequence of $G_K$-representations
$$0\to\Q_p(1)\to V \to \Q_p(0)\to 0.$$
This extension lives in $\Ext^1_{G_K}(\Q_p(0),\Q_p(1))\isom H^1(G_K,\Q_p(1))$.
One has the \emph{$\linv$-invariant} of $V$, denoted $\linv(V)$, (see \cite{colmez10}) which plays a central role in the theory of $p$-adic $L$-function of $V$, and which one may think of $\linv(V)$ as a quantity associated to the Hodge filtration on the $K$-vector space $D_{dR}(V)$ defined in \cite{fontaine94b}. For a more detailed discussion of $D_{dR}(V)$ for arbitrary ordinary representations see \Cref{th:ord-ddr-is-amphoric}.

One of the simplest, but important, consequences of anabelomorphy is the following:
\bthm\label{th:l-inv-unamphoric}
Let $K$ be a $p$-adic field. Let $V$ be as above. Then the $\linv$-invariant, $\linv(V)$, of $V$ is not amphoric. 
\ethm
\brem 
As pointed out in  \ssep\ref{ss:thta-mhs-linv}, there is an archimedean analog, $\linv_\infty$, of the $p$-adic $\linv$-invariant, and the archimedean version of the above result (\Cref{pr:linv-inf-unamph})  provides the simplest way of understanding \Cref{th:l-inv-unamphoric}.
\erem

\bp 
It will be sufficient to prove this under the assumption that $V$ is a non-split crystalline representation of $G_K$ of the form
$$0\to\Z_p(1)\to V \to \Z_p(0)\to 0.$$ Let $q(V_K)\in H^1_f(G_K,\Z_p(1)) \subset H^1_f(G_K,\Q_p(1))$ be the class of $V$ viewed as the given $G_K$-representation.  Since $V_K$ is non-split $q(V_K)\neq0$, one obtains a non-zero element of the $\Q_p$-vector space $H^1_f(G_K,\Q_p(1))$ and hence the $\Q_p$-linear subspace spanned by $q(V_K)$ gives a  point, which can be identified with the $\linv$-invariant $\linv(V_K)$  in  the projective space $\P(H_f^1(G_K,\Q_p(1)))$  of lines in the $\Q_p$-vector space $H_f^1(G_K,\Q_p(1))$. 

Now suppose one has an anabelomorphism $L\anabmapright{\alpha} K$. Then by \Cref{th:anab-galois-h1}, one has an induced isomorphism $H_f^1(G_L,\Q_p(1))\isom H_f^1(G_K,\Q_p(1))$ and hence an isomorphism of topological spaces $\P(H_f^1(G_L,\Q_p(1)))\overset{\alpha}{\isom} \P(H_f^1(G_K,\Q_p(1)))$. However, in general there is no natural isomorphism between $G_L\isom G_K$ and no natural isomorphism given by \Cref{th:anab-galois-h1} and hence no natural isomorphism between these projective spaces. To see this, it is sufficient to consider outer automorphisms of $G_K$. Assume that $p$ is odd, $[K:\Q_p]>1$. By \cite[Theorem 1.5]{hoshi2022}, there exists an outer automorphism $\sigma:G_K\mapright{\isom} G_K$  such that (1) for all  integers $n\geq 1$, $\sigma^n\neq 1$, and (2) the isomorphism induced by $\sigma^n$ on the $\Q_p$-vector space $H^1_f(G_K,\Q_p(1))$ satisfies $\sigma^n\neq 1$ for all $n$. Thus, in general anabelomorphisms induce non-trivial isomorphisms of this $\Q_p$-vector space. Write $V_{K,L,\alpha}$ for the (non-split) $G_L$ representation obtained from $V_K$ by composition with $\alpha$, then one obtains its image $\alpha(\linv(V_{K,L,\alpha}))\in \P(H_f^1(G_K,\Q_p(1)))$. As the anabelomorphism $G_L\overset{\alpha}{\isom} G_K$ and the field $L$ vary (keeping $K$ fixed), one obtains a set of points of $\P(H_f^1(G_K,\Q_p(1)))$ which lie in the $\textrm{image}(H_f^1(G_K,\Z_p(1))-\{0\}) \subset  \P(H_f^1(G_K,\Q_p(1)))$ and by what has been just said, in general, this set is not a one point set. 
\ep

This has the following consequence:
\bcor\label{th:hodge-filt-unamphoric}
Let $V\in\Ext^1_{G_K}(\Q_p(0),\Q_p(1))$. Then  the Hodge filtration on $D_{dR}(V)$ is not amphoric.
\ecor
\bp
From \cite{colmez10} one sees that $\linv(V)$ controls the Hodge filtration on the filtered $(\phi,N)$-module $D_{dR}(V)$. Therefore, one deduces that anabelomorphy changes the $p$-adic Hodge filtration. See Section~\ref{se:fontaine} for additional comments on this.
\ep

\newcommand{\bcris}{B_{cris}}
\newcommand{\bst}{B_{st}}
\newcommand{\bdr}{B_{dR}}
\newcommand{\dcrisv}[1]{D_{cris}(\ifstrempty{#1}{}{#1,}V)}
\newcommand{\dstv}[1]{D_{st}(\ifstrempty{#1}{}{#1,}V)}
\newcommand{\ddrv}[1]{D_{dR}(\ifstrempty{#1}{}{#1,}V)}
\newcommand{\ddr}[1]{D_{dR}({#1})}

\brem\label{re:comp-question}
In \cite{fontaine-colmez} it was shown that every weakly admissible filtered $(\phi,N)$ module is an admissible filtered $(\phi,N)$ module (this had been conjectured in \cite{fontaine94b}). The proof proceeds by changing the Hodge filtration on a filtered $(\phi,N)$-module. 

The idea of \cite{fontaine-colmez} is to replace the original Hodge filtration (which may make the module possibly inadmissible) by a new Hodge filtration  so that the new  module becomes admissible i.e. arises from a Galois representation. So in this situation the $p$-adic Hodge filtration is considered mobile while other structures remain fixed. This allows one to keep the $p$-adic field $K$ fixed. 

\emph{This should be viewed as an example of anabelomorphy but carried out on the $p$-adic Hodge structure}. 

\Cref{th:l-inv-unamphoric} says that the $\linv$-invariant of an elliptic curve over a $p$-adic field is not amphoric and the $\linv$-invariant is related to the filtration of the $(\phi,N)$-module \cite[3.1]{colmez10}. So the filtration is moving in some sense but the space on which the filtration is defined is also moving because the Hodge filtration for the $G_K$-module $V$ lives in the $K$-vector space $D_{st}(V)$, while the Hodge filtration for the $G_L$-module $V$ lives in an $L$-vector space. This leads to the question of comparing the two vector spaces\footnote{The question of comparing the $K$-vector space $D_{dR}(\rho,V)$ and $L$-vector space $D_{dR}(\rho\circ\alpha,V)$ was raised by S.~Mochizuki in his e-mail to the Author (circa March-April 2020): ``it remains a significant challenge to find containers where  the $K$-vector space $D_{dR}(\rho,V)$ and $L$-vector space $D_{dR}(\rho\circ\alpha,V)$ can be compared.''}.
\erem

\subsection{Amphoricity of $D_{dR}(V(r))$ ($r\gg0$) for an ordinary representation}\label{se:fontaine}
\Cref{th:ord-ddr-is-amphoric} below answers the question raised at the end of \Cref{re:comp-question} by showing  that  if $V$ is ordinary, then for all sufficiently large integers $r$, the $\Q_p$-vector space $D_{dR}(\rho,V(r))$ is amphoric i.e. there  is a natural $\Q_p$-vector space isomorphism between the $K$-vector space $D_{dR}(\rho,V(r))$ and $L$-vector space $D_{dR}(\rho\circ\alpha,V(r))$.

Let $K$ be a $p$-adic field and let $\alpha:L\anabelmap K$ be an anabelomorphism of $p$-adic fields. Consider $\reptrip$ of $G_K$. Suppose that $V$ is a de Rham representation of $G_K$ in the sense of \cite{fontaine94b}. By \Cref{th:crystalline-not-amphoric}, $\rho\circ\alpha: G_L\to \gl{V}$ need not be Hodge-Tate and hence need not be de Rham. Suppose $V$ is ordinary. Then by \cite[Th\'eor\`eme 1.5]{perrin-riou1994}, $V$ is semi-stable and hence also de Rham. By \Cref{th:ordinary-amphoric}, one deduces that the $G_L$-representation $\rho\circ\alpha:G_L\to\gl{V}$ is also ordinary and hence also de Rham. Write $\ddrv{\rho}$ for the $K$-vector space associated to the de Rham representation $\reptrip$ of $G_K$ and write $\ddrv{\rho\circ\alpha}$ for the $L$-vector space associated to the de Rham representation $\rho\circ\alpha:G_L\to\gl{V}$ of $G_L$.

\bthm\label{th:ord-ddr-is-amphoric}
Let $K$ be a $p$-adic field, let $\alpha:L\anabelmap K$ be an anabelomorphism of $p$-adic fields. Let $\reptrip$ be an ordinary $p$-adic representation of $G_K$ (so that, by \Cref{th:ordinary-amphoric}, $\rho\circ\alpha: G_L\to\gl{V}$ is also an ordinary $p$-adic representation of $G_L$). Then for all sufficiently large integers $r\gg0$ (depending only on  $\alpha$ and $(\rho,V)$), there exist natural isomorphisms of $\Q_p$-vector spaces
$$D_{dR}({\rho,V(r)})\overset{\scriptscriptstyle{\Cref{th:surj-bk-exp}}}{\isom} H^1(G_K,V(r))\overset{\alpha}{\isom} H^1(G_L,V(r))\overset{\scriptscriptstyle{\Cref{th:surj-bk-exp}}}{\isom} D_{dR}({\rho\circ\alpha,V(r)}).$$
Moreover, for all sufficiently large integers $r$:
$$\dim_{\Q_p}(D_{dR}({\rho,V(r)}))=\dim_{\Q_p}(D_{dR}({\rho,V}))$$  and hence is independent of $\alpha$ and also of such $r$.
\ethm

\brem 
Note that the Hodge filtration on the $K$-vector space $D_{dR}(\rho,V(r))$ is up to shifting, the filtration on the $K$-vector space $D_{dR}(\rho,V)$.   However, at the moment, we do not know how to compare the Hodge filtrations on $D_{dR}({\rho,V(r)})$ and $D_{dR}({\rho\circ\alpha,V(r)})$.
\erem

\newcommand{\ddc}{D_{cris}}

\bthm\label{th:surj-bk-exp}
Let $K$ be a $p$-adic field and let $\reptrip$ be an ordinary $p$-adic representation of $G_K$. Then for all sufficiently large integers $r\gg 0$, the Bloch-Kato exponential homomorphism 
$$\ddr{V(r)} \mapright{exp_{BK}} H^1_e(G_K,V(r))\mapright{\isom} H^1_f(G_K,V(r))\mapright{\isom} H^1(G_K,V(r)).$$
is an isomorphism of $\Q_p$-vector spaces.
\ethm
\bp
The Bloch-Kato exponential is defined in \cite[Definition 3.10]{bloch90}. Since  $H_e^1(G_K,V(r))\into H^1_f(G_K,V(r))\into H^1(G_K,V(r))$ are natural inclusions of $\Q_p$-subspaces by their definitions (\cite[3.7.2]{bloch90}), one may view the Bloch-Kato exponential as taking values in $H^1(G_K,V(r))$.

By \cite[1.15 Theorem]{nekovar93} or \cite[Corollary 3.8.4]{bloch90}, one has an exact sequence of $\Q_p$-vector spaces
$$0\to H^0(V(r))\to \ddc^{f=1}(V(r))\to \ddr{V(r)}/\ddr{V(r)}^0\mapright{exp_{BK}} H^1_e(G_K,V(r))\to0.$$
Since the twist $r$ is very large, one sees that $\ddc(V(r))$ has no subspace on which Frobenius $f$ acts by $1$ and hence $\ddc^{f=1}(V(r))=0$ (this part of the proof does not use the assumption that $V$ is ordinary and crystalline). As $V$ is ordinary and the twist $r$ is large and shifts the Hodge filtration on $\ddr{V(r)}$, one sees that the filtration $\ddr{V(r)}^0=0$. Thus, one obtains the isomorphism
$$\ddr{V(r)}\mapright{\isom} H^1_e(G_K,V(r))$$
induced by the Bloch-Kato exponential $exp_{BK}$.

Next by \cite[1.16 Corolalry]{nekovar93} one has an exact sequence
$$0 \to H^1_e(G_K,V(r))\to H^1_f(G_K,V(r)) \to \ddc(V(r))/(1-f)\ddc(V(r))\to0.$$
Again as $r\gg 0$, one sees that $\ddc(V(r))$ has no quotient on which Frobenius acts by $1$. Thus, one sees that $H^1_e(G_K,V(r))\isom H^1_f(G_K,V(r))$. 

Now  the formula \cite[3.8.5]{bloch90} says
$$\dim_{\Q_p}(H_f^1(G_K,V(r)))+\dim_{\Q_p}(H_f^1(G_K,V^*(1-r)))=H^1(G_K,V),$$
where $V^*$ is the $G_K$-representation dual to $V$. Thus to prove the theorem, it is enough to prove that 
$$H_f^1(G_K,V^*(1-r))=0.$$ Since $V(1-r)$ has sufficiently negative slopes and is an ordinary representation, this required vanishing  follows from the following lemma and this completes the proof of \Cref{th:surj-bk-exp}.
\ep

\blem Let $K$ be a $p$-adic field and let $\rho:G_K\to \gl W$ be an ordinary representation such that all the Frobenius slopes of $D_{st}(W)$ are sufficiently negative. Then $$H^1_f(G_K,W)=H^1_g(G_K,W)=0.$$
\elem
\bp
Since $H^1_f(G_K,W)\subseteq H^1_g(G_K,W)$, it is enough to show that $H^1_g(G_K,W)=0$.
This will be proved by induction on $\dim_{\Q_p}(W)$. 
If $\dim_{\Q_p}(W)=1$, and $W=\Q_p(m)$ with $m\ll 0$ and hence by the table in \cite[Example 3.9]{bloch90} one obtains the desired vanishing. But, in general, one has (for $\dim(W)=1$) that $W=\chi\tensor\Q_p(m)$ for some unramified character $\chi$ of $G_K$. We claim that if $m$ is sufficiently negative, then $H^1_g(G_K,\chi\tensor\Q_p(m))=0$.  To prove this, one uses the following formulae given by \cite[1.24 Proposition]{nekovar93} for any de Rham representation.
\begin{align}
h_f^1(W)&=h^0(W)+[K:\Q_p]\dim_K(D_{dR}(W)/F^0)\\
h_g^1(W)&=h^1_f(W)+\dim_{\Q_p}D_{cris}(W^*(1))^{f=1},
\end{align}
where $h_{*}^i(W)=\dim_{\Q_p}H^i_*(G_K,W)$ for $*\in \{g,f\}$. Thus, to prove this, we have to show that all the terms entering the formula for $h^1_g(W)$ are equal to zero. This is where the hypothesis $W=\chi\tensor\Q_p(m)$ with $\chi$ being unramified and $m\ll0$ comes into play. Clearly, $h^0(W)=0$. As $m\ll 0$, one has $D_{dR}(W)=F^0$ and hence $h^1_f(W)=0$. Again as $m\ll0$, $W^*=\chi^{-1}\tensor \Q_p(-m)$ and so $W^*(1)=\chi^{-1}\tensor \Q_p(1-m)$ with $1-m\gg 0$ and so $D_{cris}(W^*(1))^{f=1}=0$. This shows that $h^1_g(W)=0$ as claimed.
 
Now suppose $\dim_{\Q_p}(W)>1$. Then since $W$ is an ordinary representation, one has an exact sequence
$$0\to W_2\to W \to W_1\to 0$$ where $W_1=\chi\tensor \Q_p(m)$  with $\chi$ being unramified character and $m$ is also sufficiently negative and $W_2$ is also an ordinary representation with sufficiently negative slopes. By \cite[1.25]{nekovar93}, one has an exact sequence
$$H^0(G_K,W_1)\to H^1_g(G_K,W_2)\to H^1_g(G_K,W) \to H^1_g(G_K,W_1).$$
The term $H^1_g(G_K,W_2)=0$ by induction hypothesis and as shown earlier $H^1_g(G_K,W_1)=0$
and hence the middle term is zero by exactness as asserted. This completes the proof.
\ep

\bp[Proof of \Cref{th:ord-ddr-is-amphoric}]
Using the anabelomorphism $L\anabmapright{\alpha} K$ one sees that $G_L$ acts on $V$ through the isomorphism $\alpha:G_L\isom G_K$.  So $V$ is also a $G_L$-module. By the amphoricity of the cyclotomic character given by \Cref{th:third-fun-anab} one has  compatibility with Tate twists. Hence, $V(r)$ is also a $G_L$-module for any integer $r$. Then as $G_K\isom G_L$, one has an isomorphism of $\Q_p$-vector spaces (given by $\alpha$):
$$H^1(G_K,V(r))\isom H^1(G_L,V(r)).$$

By \Cref{th:surj-bk-exp}, for all $r\gg0$, one has isomorphisms of $\Q_p$-vector spaces
$$\ddr{\rho, V(r)}\isom H^1(G_K,V(r))\isom H^1(G_L,V(r))\isom \ddr{\rho\circ\alpha,V(r)}.$$ The proof of \Cref{th:surj-bk-exp} makes it clear that $r$ depends only on $(\rho,V)$.
This proves the first assertion.

So it remains to prove the last assertion. Since $V$ is ordinary,  by \Cref{th:ordinary-amphoric},  $\rho\circ\alpha$ is also ordinary and hence both $\rho$ and $\rho\circ\alpha$ are semi-stable by \cite[Th\'eor\`eme 1.5]{perrin-riou1994}. Let $K\supset K_0$ be the maximal unramified subfield of $K$. Then one has by \cite[5.1.7]{fontaine94b}, $D_{dR}(\rho,V(r))=K\tensor_{K_0}D_{st}(\rho,V(r))$. By \cite[2.2 Lemme]{perrin-riou1994}, one knows that $D_{st}(\rho,V)[r]=D_{st}(\rho,V(r))$. The twist $[r]$ on the filtered $(\phi,N)$-module $D_{st}(\rho,V)[r]$ shifts the filtration on the $K$-vector space $K\tensor_{K_0}D_{st}(\rho,V)$. Thus, the dimension of the $K_0$-vector space $D_{st}(\rho,V(r))$ is independent of $r$. Hence, the dimension of $K$-vector space $D_{dR}(\rho,V(r))$ is independent of $r$ (for all sufficiently large integers $r$). Since $[K:\Q_p]$, $[K_0:\Q_p]$ are amphoric by \Cref{th:third-fun-anab}, the last assertion is proved.
\ep

\subsection{Anabelomorphy and Deformations of Galois representations}
\newcommand{\rabk}{R^{[a,b], K'}_{V_F}}
\newcommand{\rabkp}{R^{[a,b], L'}_{V_F'}}
\newcommand{\rabok}{R^{ord}_{V_F}}
\newcommand{\rabokp}{R^{ord}_{V_F'}}
The assertion is the following:
\bthm\label{th:anab-deform-thry} 
Let $K$ be a $p$-adic field of residue characteristic $p$ for some prime number $p$. Let $a\leq b\in \Z$ be two integers. Let $\F=\F_q$ be a finite field with $q$ elements and of characteristic $p$. Let $W=W(\F)$ be the ring of Witt vectors of $\F$.  Let $V_{\F}$ be a finite dimensional representation of $G_K$ with values in $\F$ and $\End_{\F[G_K]}(V_\F)=\F$.  Let $R_{V_\F}$ be the deformation ring $R_{V_\F}$ of the $G_K$ representation $V_{\F}$. For any finite extension $K'$ of $K$, let $\rabk$ be the deformation ring of  $G_K$-representations with values in finite extensions of $W[1/p]$, which become semi-stable representations   of Hodge-Tate weights in $[a,b]$ when restricted to $G_{K'}\subseteq G_K$. Let $L\anabmapright{\alpha} K$ be an anabelomorphism of $p$-adic fields. Let $V'_F$ be the $G_L$-representation obtained from $V_F$ by composing with the anabelomorphism $\alpha$. Then
\benumlab
\item The deformation ring $R_{V_\F'}$ of $V_\F'$ also exists and the anabelomorphism $\alpha$ induces an isomorphism of the deformation rings $\alpha:R_{V_\F'}\isom R_{V_\F}$,
\item  and an isomorphism of ordinary deformation rings $\alpha:\rabokp\isom \rabok$.
\item The canonical quotient $R_{V_\F}\to\rabok$ is amphoric.
\item In general, the canonical quotient $R_{V_\F}\to {R^{[a,b], K}_{V_F}}$ is not amphoric.
\eenum 
\ethm
\bp 
One sees trivially that the anabelomorphism $G_L\overset{\alpha}{\isom} G_K$,  induces an equivalence of deformation groupoids considered in \cite{bockle2013}; moreover, $\alpha$ also induces an isomorphism $\End_{\F[G_L]}(V_\F')\isom\End_{\F[G_K]}(V_\F)=\F$, hence the deformation rings $R_{V_\F}$ and $R_{V_\F'}$ exist (\cite[Theorem 2.11]{bockle2013}) and are naturally isomorphic. Thus {\bf(1)} is immediate. 

Similarly, \Cref{th:ord-syn-thm}, and the amphoricity of the cyclotomic character \Cref{th:third-fun-anab}, shows that the anabelomorphism $\alpha$ also induces a natural equivalence between the groupoids of ordinary deformations of the $G_K$-module $V_\F$ (resp. the $G_L$-module $V_{\F}'$). Hence, it induces an isomorphism of the ordinary deformation rings $\rabok\isom \rabokp$ compatible with the isomorphism  $R_{V_\F}\isom R_{V_\F'}$. This proves {\bf(2,3)}. 

For finite extensions $K'/K$ (resp. $L'/L$), the existence of the potentially semi-stable deformation rings $\rabk, \rabkp$ is established in \cite{kisin2008}. \Cref{th:crystalline-not-amphoric}{\bf(3)} shows that among all integers $a\leq b$, all finite fields $F$ of characteristic $p$, all anabelomorphisms $L\anabmapright{\alpha} K$ of $p$-adic fields, and all $G_K$-modules $V_F$ and all potentially semi-stable lifts of $V_F$ of Hodge-Tate weights in $[a,b]$, there exists some finite field $F$ of characteristic $p$, some finite dimensional  $G_K$-module $V_F$, some anabelomorphism of $p$-adic fields $G_L\anabelmap G_K$, and some potentially semi-stable lift of $V_F$ which is not Hodge-Tate when viewed as a $G_L$-representation lifting the $G_L$-representation $V_F$. Hence, one sees that ${R^{[a,b], K}_{V_F}}$ is not amphoric in general. This  proves {\bf(4)}.
\ep

\newcommand{\sG}{\mathscr{G}}
\newcommand{\sH}{\mathscr{H}}

\section{Anabelomorphy and the local Langlands correspondence}\label{se:anab-langlands}\nwss
Let $K\anabelmap L$ be anabelomorphic $p$-adic fields.  One sees from \Cref{le:weil-deligne-groups} (below), any anabelomorphism $K \anabelmap L$ provides a natural isomorphism of Weil-Deligne groups of $K$ and $L$ respectively. The local Langlands correspondence matches certain representations of Weil-Deligne group $W_K'$ of $K$ to certain of representations of $\gln K$. But as far as one is aware, topological groups $\gln{K}$ and $\gln{L}$ are not known to be topologically homeomorphic as groups except for  $n=1$ (\Cref{th:third-fun-anab}). Thus, one is led, by the results of \ssep\ref{se:anabel-galois-reps} and \Cref{le:weil-deligne-groups} to the following questions.  Given an anabelomorphism $K\anabelmap L$ of $p$-adic fields:
\benumlab
\item  how to construct a functor (natural in anabelomorphisms) between the categories of representations of $\gln K$ and $\gln L$,  and secondly
\item how to construct a correspondence between automorphic representations of $\gln K$ and $\gln L$ (here automorphic representations  will mean admissible (equivalently smooth), complex valued representations) of $\gln{K}$ and $\gln L$.  
\eenum
One of the referees remarked that  the first question is in the cadre of Grothendieck's mysterious functor question (see \cite[Page 17]{Illusie2015GrothendieckPisa}), and one has no answer for it (except for $n=1$ where it is trivial), but the second question is answered here for principal series representations (\Cref{th:automorphic-ordinary-synchronization}) for all $n\geq1$, and in \Cref{th:syn-super-cusp-n} for all irreducible supercuspidal representations in the so-called tame case (i.e. $p$ does not divide $n$). This is the main theme of this section.  In what follows, we will use the terms `\emph{synchronize}' or `\emph{synchronization}' as synonyms for the existence, or the construction, of natural bijections between various sets (for example the set of irreducible admissible representations). 
\subsection{Anabelomorphisms of Weil and Weil-Deligne Groups}
The following results will be used in the subsequent discussions.

\blem\label{le:norm-amphoric}
Let $K$ be a $p$-adic field. Let $q_K$ be the cardinality of the residue field of $K$. Then
\benumlab
\item $q_K$ is amphoric.
\item The homomorphism $\ord_K:K^*\to \Z$ given by $x\mapsto \ord_K(x)$ is amphoric.
\item the homomorphism $\norm{-}_K:K^*\to \R^*$ defined by $\norm{x}=q_K^{-\ord_K(x)}$ is amphoric.
\eenum
\elem

\bp 
One has $q_K=p^{f_K}$, where   $f_K=[K:\Q_p]/e_K$, and by \Cref{th:third-fun-anab}{\bf(1,2)} one sees that $f_K$ is amphoric. This proves {\bf(1)}. It is clear that the third assertion follows from the second. So it is sufficient to prove the second assertion. This is proved as follows: by \Cref{th:third-fun-anab}, $K^*$ is amphoric. By \cite[Page 144]{Cassels1967}, one sees  that the composite $$K^*\mapright{rec_K} G_K^{ab}\to Gal(K^{nr}/K)=\hat{\Z}$$
is the valuation map $ord_K:K^*\to \Z$. A uniformizing element in $K^*$ maps to a Frobenius element under the reciprocity map $rec_K$ and its image in $\hat{\Z}$ is $1\in\Z$. By \cite[Proposition 3.11]{hoshi-mono}, one deduces the required amphoricity assertion of {\bf(2)}.
\ep

\bpro\label{le:weil-deligne-groups}
Let $K$ be a $p$-adic field and let $\alpha:K\anabelmap L$ be an anabelomorphism. Let $W_K$ (resp. $W_L$) be the Weil group of $K$ (resp. $L$) and let $W_K'$ (resp. $W_L'$) be the Weil-Deligne group of $K$ (resp. $L$).
Then the anabelomorphism $K \anabmapright{\alpha} L$ induces natural  topological isomorphisms of Weil groups and Weil-Deligne groups:
\benumlab 
\item $W_K\isom W_L$, and
\item $W_K' \isom W_L'$
\eenum
such that a Frobenius element of $W_K$ maps to a Frobenius element of $W_L$ (and resp. for Weil-Deligne groups).
\epro
\bp 
The anabelomorphism $\alpha : K\anabelmap L$ gives an isomorphism $\alpha: G_K\to G_L$. By \Cref{le:norm-amphoric}{\bf(1)}, the cardinality $q=q_K$ of the residue field of $K$ is amphoric. Let $\fq$ be the residue field of $K$ (and hence of $L$). The anabelomorphism $\alpha:G_K\mapright{\isom} G_L$ together with the amphoricity of the inertia subgroup $I_K\subset G_K$ gives us a commutative diagram of homomorphisms of groups in which horizontal arrows are isomorphisms:
$$
\begin{tikzcd}
G_K \arrow[r, "\alpha"] \arrow[d] & G_L\arrow[d]\\
G_K/I_K \arrow[r,"\tilde{\alpha}"] & G_L/I_L.
\end{tikzcd}
$$
Let ${\rm Frob}_K\in G_K$ be a Frobenius element for $K$.  This is a well-defined element of  $G_K/I_K$. By \cite[Proposition 3.9]{hoshi-mono}, $\alpha({\rm Frob}_K \pmod{I_K})={\rm Frob}_L\pmod{I_L}\in G_L/I_L$. 
The Weil group $W_K\subset G_K$ is the subgroup of elements $g\in G_K$ such that $g\pmod{I_K}\in {\rm Frob}_K^\Z\in G_K/I_K$. Hence, under the anabelomorphism $G_K\anabmapright{\alpha} G_L$, one has $$\tilde{\alpha}(g)\in\tilde{\alpha}({\rm Frob_K})^\Z={\rm Frob}_L^\Z.$$
Thus one sees that $\alpha(W_K)\subseteq W_L$.
Since, starting from $G_L$ and the inverse $\alpha^{-1}$ one arrives from $W_L$ into $W_K$ one sees that $\alpha$ induces an isomorphism of Weil groups $\alpha:W_K\mapright{\isom} W_L$. This proves {\bf(1)}. The assertion for Weil-Deligne groups is immediate from this, \Cref{le:norm-amphoric}{\bf(3)} and the definition of the Weil-Deligne group.
\ep

\subsection{Amphoricity of $\sS{K}$ and $\sS{K^*}$ and its consequences}\nwss
Let $K$ be a $p$-adic field. Then $(K,+)$ (resp. $(K^*,\times)$) is a locally compact topological group. Let $\sS{K}$ (resp. $\sS{K^*}$) be the space of locally constant, compactly supported, complex valued, continuous functions on $(K,+)$ (resp. $(K^*,\times)$). 

Let $\dmu{K}{x}$ (resp. $\dmus{K}{x}$) be a Haar measure on $(K,+)$ (resp. $(K^*,\times)$). One may also choose $\dmu{K}{x}$ to be normalized as in \cite[Paragraph before Theorem 2.2.2]{tate-thesis}. In \cite[Sections 2.2, 2.3]{tate-thesis}, one identifies the topological group $K^*$ as a subset of the topological group $(K,+)$. \textit{This identification uses the fact that $K$ is a $p$-adic field and $K^*=K-\{0\}$ is the subset of non-zero elements of the field $K$}. Notably \cite[Lemmas 2.2.4, 2.2.5]{tate-thesis} use the field structure of $K$ to identify the Haar measure on $(K^*,\times)$ as  
\be\label{eq:tate-haar} 
\dmus{K}{x}=\frac{\dmu{K}{x}}{\abs{x}_K}.
\ee

\bthm\label{th:meas-thry}
Let $K$ be a $p$-adic field. Let $\dmu{K}{x}$ be the Haar measure on $(K,+)$ defined by \cite{tate-thesis}, let $\dmus{K}{x}$ be the Haar measure on $(K^*,\times)$ given (using the field structure of $K$) by \eqref{eq:tate-haar}.
\benumlab
\item The spaces $\sS K$ and $\sS{K^*}$ are amphoric.
\item The pair $(\sS{K}, \dmu{K}{x})$ is  amphoric, 
\item  but  the pair $(\sS{K^*}, \dmus{K}{x})$ is not amphoric.
\item The identification  $(K,+)$ with its character group is not amphoric (in general).
\eenum
\ethm
\bp 
Let $K\anabmapright{\alpha} L$ be an anabelomorphism of $p$-adic fields.
From \Cref{th:third-fun-anab} and as noted in the proof of \Cref{th:anab-affine-spaces}, $\alpha$ functorially provides an isomorphism  $\alpha:(K,+)\mapright{\isom} (L,+)$ of topological groups. Similarly, one obtains from the anabelomorphism $\alpha$, an isomorphism of topological groups $(K^*,\times)\mapright{\isom} (L^*,\times)$. This proves {\bf(1)}. By \cite[Lemma 3.12, Summary 3.15]{hoshi-mono}, the Haar measure $\dmu{K}{x}$ is amphoric. This proves {\bf(2)}. As remarked earlier,  the measure $\dmus{K}{x}=\frac{\dmu{K}{x}}{\abs{x}}$ is defined using the field structure of $K$ via the inclusion $K^*\into K$. To prove that it is not amphoric it will suffice to prove that while $\alpha(\O_K^*)= \O_L^*$ under the isomorphism $\alpha:K^*\mapright{\isom} L^*$, $\O_K^*$ and $\O_L^*$ have distinct volumes with respect $(K^*,d\mu_K^*)$ and $(L^*, d\mu_L^*)$. This is immediate from volume computation of \cite[Lemma 2.3.3]{tate-thesis} which gives
$$\int_{\O_K^*}\dmus{K}{x}=\abs{\mathfrak{d}_K}_K^{-1/2},$$ where $\mathfrak{d}_K$ is the discriminant ideal of $K$.
As discriminants of $p$-adic fields are not amphoric in general (\Cref{th:discriminant-is-unamphoric}), one deduces {\bf(3)}. The canonical identification of $(K,+)$ with its own character group proved in \cite[Lemma 2.2.1]{tate-thesis} depends on the field structure of $K$ and hence is not amphoric in general. This proves {\bf(4)}. Thus, one has proved all the assertions.
\ep

As is standard \cite{tate-thesis}, it makes perfect sense to talk about integrals of functions $f(x)\in\sS K$
$$\int_K f(x)\dmu K x.$$

\newcommand{\sI}[1]{\mathscr{I}_{#1}}
For a $p$-adic field $K$, let $\log$ denote the $p$-adic logarithm and let $$\sI K=\frac{1}{2\cdot p}\log(\O_K^*)\subset (K,+)$$ be the \textit{log-shell} defined by Mochizuki (see \cite[Definition 1.1]{hoshi-mono}).

\bcor\label{cor:l-func-int}
Let $K \anabmapright{\alpha} L$ be an anabelomorphism of $p$-adic fields. Then
\benumlab
\item For all $f\in\sS{K}$ one has
$$\int_K f(x)\dmu{K}{x} =\int_L f(\alpha(x))\dmu{L}{\alpha(x)}.$$
\item In general $\alpha(\O_K)\neq \O_L$ under the isomorphism $\alpha:(K,+)\mapright{\isom} (L,+)$.
\item The log-shell $\sI K$ is an amphoric, topological $\Z_p$-submodule of $(K,+)$.
\eenum
\ecor
\bp 
The first assertion is clear. By the choice of normalization of the Haar measure $\dmu{K}{x}$ in \cite[Paragraph before Theorem 2.2.2]{tate-thesis} one has
$$\int_{\O_K}\dmu{K}{x}=\abs{\mathfrak{d}_K}_K^{-1/2},$$
where $\mathfrak{d}_K$ is the discriminant of $K$ and $\abs{\mathfrak{d}_K}_K$ is its absolute value, and
 $$\int_{\O_L}\dmu{L}{x}=\abs{\mathfrak{d}_L}_L^{-1/2}.$$
If $\alpha(\O_K)=\O_L$ then by {\bf(1)}, the two would have the same volume. But as discriminants are not amphoric (in general) by \Cref{th:discriminant-is-unamphoric}, the assertion {\bf(2)} follows.  The assertion {\bf(3)} is due to Mochizuki, a proof is given in \cite[Proposition 3.11(iv)]{hoshi-mono}.
\ep

\brem \Cref{cor:l-func-int}{\bf(2)} provides a natural numerical proof of the fact that $(\O_K,+)\subset (K,+)$ is not amphoric; for a less direct argument see \cite[Remark 4.3.1(ii)]{hoshi-mono}.
\erem

\subsection{Anabelomorphic Synchronization of Principal Series Representations}
For a $p$-adic field $K$, a \textit{quasicharacter} of $\glo K$ is a continuous homomorphism $K^*\to \C^*$. An admissible irreducible representation of $\glo K$ is the same as a quasicharacter of $\glo K$.
\bthm\label{th:anab-glo-reps} 
Let $\alpha: L\anabelmap K$ be an anabelomorphism of $p$-adic fields. Then $\chi\mapsto \chi\circ \alpha$ sets up a bijection between irreducible admissible representations of $\glo K$ and $\glo L$ respectively, under which $L$-functions are amphoric, but conductors and $\varepsilon$-factors are not amphoric in general.
\ethm
\bp 
The local Langlands correspondence sets up a bijection between admissible representations of $\glo K$ and one dimensional representations of the Weil-Deligne group $W_K$ (with $N=0$) which matches $L$-functions, conductors and $\varepsilon$-factors. 

Any anabelomorphism $\alpha: L\anabelmap K$ induces an isomorphism $\alpha:L^*=\glo L\to \glo K =K^*$ (\Cref{th:third-fun-anab}) and by \Cref{le:weil-deligne-groups} one also has an induced isomorphism $W_L'\isom W_K'$ of Weil-Deligne groups. Thus, one obtains a bijection between irreducible admissible $\glo K$ representations and irreducible admissible representations of $\glo L$ which is compatible with the local Langlands correspondence.

By \Cref{th:char-poly-amphoric}, $L$-functions of  Galois representations are amphoric. Hence, $L$-functions of irreducible admissible representations of $\glo K$ are amphoric.  So it remains to prove the assertion about conductors and $\varepsilon$-factors. Let $\chi:K^*\to \C^*$ be a quasicharacter and let $\varpi\in\O_K$ be a uniformizer. The conductor of $\chi$  is the smallest integer $n\geq 0$ such that $\chi(1+\varpi^n\O_{K})=1$ but $\chi(1+\varpi^{n-1}\O_{K})\neq 1$. Now suppose $\alpha:K\anabelmap L$ is an anabelomorphism of $p$-adic fields. Then by \Cref{th:third-fun-anab}, one has an induced isomorphism $\alpha: \O_K^*\isom \O_L^*$. Thus, groups $\O_K^*,\O_L^*$ have isomorphic character groups. However, $\alpha:\O_K^*\isom \O_L^*$ does not preserve the ramification filtration on $\O_K^*,\O_L^*$ in general (this is a consequence of the proof of the main theorem of \cite{mochizuki-local-gro}). A direct proof of this fact can be found in the proofs of the explicit examples of strictly anabelomorphic $p$-adic fields given in \cite[\ssep 2 Examples, Theorem]{yamagata76} also show, in general this induced isomorphism does  not  preserve the ramification filtrations on these groups. This means that, in general, the conductor of a character of $\O_L^*$ need not be the same as that of the character of $\O_K^*$ obtained by composing with the isomorphism $\O_K^*\isom \O_L^*$.  Thus, conductors of quasicharacters are not amphoric in general.

Now  let us establish the assertion for $\varepsilon$-factors. Any anabelomorphism $\alpha:K\anabelmap L$ gives an isomorphism of the additive groups $\alpha:(K,+)\isom (L,+)$ (\Cref{th:third-fun-anab}). However, this isomorphism does not preserve the topological subgroups $\O_K$ and $\O_L$ (\Cref{cor:l-func-int}{\bf(2)}), and hence does not preserve the ring structure of $K$, $L$ in general (\cite[Remark 4.3.1(iii)]{hoshi-mono}) and hence it does not preserve the filtration by the powers of the respective maximal ideals (in general).  By \cite[3.4.3.4]{deligne1973b}, one sees that the $\varepsilon$-factor of a quasicharacter $\chi$ of $K^*$  depends on the conductor of $\chi$, and the conductor of a chosen additive character $\psi:(K,+)\to\C$.  Since  the natural filtration on $\O_K^*$ (resp. on $\O_K$) is not amphoric, and the conductor is not amphoric, one sees that $\varepsilon$-factors of quasicharacters are not amphoric in general.
\ep

\brem
Since $L$-functions are defined in \cite{tate-thesis} using the pair $(\sS{K^*}, \dmus{K}{x})$,  \Cref{th:meas-thry}{\bf(3)} indicates that amphoricity of $L$-functions proved above is quite subtle and one does not have a direct way of establishing it.
\erem

The following theorem is the local automorphic analog of \Cref{th:ord-syn-thm}.
\bthm[Automorphic Ordinary Synchronization Theorem]\label{th:automorphic-ordinary-synchronization} 
Let $\alpha:K\anabelmap L$ be an anabelomorphism of $p$-adic fields. Then there is a natural bijection between principal series representations of $\gln{K}$ and principal series representations of $\gln{L}$ which is given by $$\pi(\chi_1,\ldots,\chi_n)\mapsto \pi(\chi_1\circ\alpha,\ldots,\chi_n\circ\alpha).$$ This correspondence takes irreducible principal series representations of $\gln K$ to irreducible principal series representations of $\gln L$.
\ethm
\bp 
The datum required to give a principal series representations of $\gln{K}$ consists of an $n$-tuple of quasicharacters $(\chi_1,\ldots,\chi_n)$ of $K^*$ with values in $\C^*$.  The associated principal series representation is denoted by $\pi(\chi_1,\ldots,\chi_n)$ and every principal series representation is of this type. 

Now let $\alpha:L\anabelmap K$ be an anabelomorphism, so one has the induced isomorphism $\alpha:L^*\mapright{\isom} K^*$. The correspondence $(\chi_1,\ldots,\chi_n)\mapsto (\chi_1\circ\alpha,\ldots,\chi_n\circ\alpha)$ sets up a bijection between $n$-tuples of quasicharacters of $K^*\to \C^*$ and $L^*\to\C^*$. Since every principal series representation $\pi$ of $\gln{K}$ is of the form $\pi=\pi(\chi_1,\ldots,\chi_n)$ (similarly for $\gln{L}$), the first part of the assertion is immediate. 

Now it remains to prove that, under this correspondence, an  irreducible principal series representation is mapped to an irreducible principal series representation. For this it is sufficient to note that if $\chi_i\cdot \chi_j=\norm{-}_K^{\pm1}$, then by \Cref{le:norm-amphoric}, so is $(\chi_i\circ\alpha)\cdot(\chi_j\circ\alpha)=\norm{-}_K^{\pm1}\circ\alpha=\norm{-}_L^{\pm1}$. So under this correspondence an irreducible principal series representation $\pi$ is mapped to an irreducible principal series representation.
\ep

\bthm 
Let $K\anabmapright{\alpha} L$ be an anabelomorphism of $p$-adic fields. Let $$\sH_{K,n}=\sH(\gln K, \gln{\O_K})$$ be the Hecke algebra of $\gln K$ with respect to the maximal compact subgroup $\gln{\O_K}$, similarly let $\sH_{L,n}$ be the standard Hecke algebra  of $\gln L$ with respect to $\gln{\O_L}$. Then for all $n\geq 1$, $\alpha$ induces a natural isomorphism $$\alpha:\sH_{K,n} \mapright{\isom} \sH_{L,n}.$$
\ethm
\bp 
Let $S_n$ be the symmetric group on $n$ letters. From \cite{satake1963} one has the isomorphism:
$$\sH_{K,n}=\sH(\gln K, \gln{\O_K}) \isom \C[\underbrace{(K^*/\O_K^*), \ldots, (K^*/\O_K^*)}_{n}]^{S_n}.$$
Hence, one sees from \Cref{th:third-fun-anab} that one has a natural isomorphism 
$$\sH_{K,n}= \C[\underbrace{(K^*/\O_K^*), \ldots, (K^*/\O_K^*)}_{n\text{ times}}]^{S_n} \isom \C[\underbrace{(L^*/\O_L^*), \ldots, (L^*/\O_L^*)}_{n\text{ times}}]^{S_n}=\sH_{L,n}.$$
\ep

\brem 
One expects that any anabelomorphism $K\anabelmap L$ of $p$-adic fields sets up a natural bijection between the sets of compact open subgroups of $\gln K$ and $\gln L$. Hence one should expect general versions of the above result.
\erem

\subsection{Anabelomorphy and supercuspidal representations of $\gln K$ for $(p,n)=1$}
\bthm\label{th:syn-super-cusp-n} 
Let $L\anabelmap K$ be an anabelomorphism of $p$-adic fields. Let $n\geq 1$ be an integer and assume that the residue characteristic $p$ of $K$ satisfies $(p,n)=1$ i.e. $p$ is coprime to $n$. Then any anabelomorphism $L\anabmapright{\alpha} K$ induces a natural bijection between isomorphism classes of irreducible supercuspidal representations of $\gln K$ and $\gln L$ respectively.
\ethm
\bp
Let $K,L$ be anabelomorphic $p$-adic fields. Let $U_K^1\subset \O_K^*$ (resp. $U_L^1\subset \O_L^*$) be the subgroup of $1$-units of $K$ (resp. $L$). By \Cref{th:third-fun-anab}, the group $K^*$ is amphoric and by \cite[Summary 3.15]{hoshi-mono}, the group $U_K^1$ is amphoric. 

Let $K_1/K$ be a finite extension of $p$-adic fields. 
Recall, from \cite{howe1977}, that a quasicharacter $\chi:K_1^*\to \C^*$ is said to be \textit{admissible} if (1) $\chi$ does not factor through the norm homomorphism $N_{K_1/F}:K_1^*\to F^*$ for any subfield $K_1\supset F\supseteq K$ and a quasicharacter $\varphi:F^*\to\C^*$, and (2) if $\chi|_{U_{K_1}^1}$ arises from the norm $N_{K_1/F}$ for some subfield $K_1\supseteq F\supseteq K$, then $F/K$ is unramified. 

Now we claim the following: suppose $L\anabmapright{\alpha}K$ is an anabelomorphism of $p$-adic fields and suppose $(K_1/K,\chi)$ is a quasicharacter of $K_1$ with $[K_1:K]=n$. Then there exists a quasicharacter $(L_1/L,\psi)$
with an anabelomorphism $L_1\anabelmap K_1$ induced by $L\anabmapright{\alpha}K$ (hence $[L_1:L]=n$).  To construct $L_1/L$ one proceeds as follows. Let $H\subset G_K$ be the open subgroup corresponding to $K_1/K$. Then let $H'=\alpha^{-1}(H)\subset G_L$ be the open subgroup of index $n$. Let $L_1=\bL^{H'}$ be the fixed field of $H'$. Clearly $\alpha_1:H'\mapright{\alpha} H$ is an isomorphism induced by $\alpha$ and hence $L_1\anabelmap K_1$ and this is compatible with $\alpha$.  If one has a  subfield $K_1\supseteq F\supseteq K$ then via the isomorphism $\alpha_1:L_1^*\isom K_1^*$ induced by the anabelomorphism $L_1\anabmapright{\alpha_1} K_1$, one obtains an anabelomorphic subfield $L_1\supset F'\supseteq L$. As noted earlier,  one has amphoricity of the unit groups and $1$-units and hence one has $L_1^*\isom K_1^*$, $(F')^*\isom F^*$ and $U_{F'}^1\isom U_{F}^1$. Furthermore, if $F/K$ is unramified, then $F'/L$ is unramified. This is true because the absolute ramification degrees $e_{K}, e_{K_1}, e_F$ (of $K_1,K,F$) are all amphoric (by \Cref{th:third-fun-anab}). But one has $e_{F/K}\cdot e_K=e_{F}$ and hence $e_{F/K}$ is also amphoric. This implies that if $F/K$ is unramified, then so is $F'/L$.

The construction of the quasicharacter $\chi':L_1^*\to\C^*$ from the datum $(K_1,\chi)$ is now clear.
Since the quasicharacter $(L_1,\chi')$ is constructed using the anabelomorphism $\alpha_1:L_1\anabelmap K_1$ and the quasicharacter $\chi:K_1^*\to \C^*$.  The  argument of the preceding paragraph shows that this construction takes an admissible quasicharacter datum $(K_1,\chi)$ to an admissible quasicharacter datum $(L_1,\chi')$. 

One knows that, for $(n,p)=1$,  each galois conjugacy class of pairs $(K_1,\chi)$ (with finite extensions $K_1/K$ satisfying $[K_1:K]=n$), \cite[Theorem 2]{howe1977}  associates an irreducible supercuspidal representation $\pi(K_1,\chi)$  of $\gln K$; and by \cite[Corollary 3.4.9]{moy1986} one knows, for $(n,p)=1$, that  all irreducible supercuspidal representations of $\gln K$ arise this way. 

Thus the correspondence asserted by the theorem is the correspondence $\pi(K_1,\chi)\mapsto \pi(L_1,\chi')$ obtained by associating $(K_1,\chi)\mapsto (L_1,\chi')$ described above. This completes the proof of the theorem.
\ep

\brem
One expects that the above result is also true for $p|n$, but from the discussion of the $p=2$ and $\GL_2$ case in \cite{bushnell-book}, it is likely that proofs will be complicated.
\erem

\subsection{Anabelomorphic synchronization of Weil representations of $\glt K$}
\bthm\label{th:weil-rep-gl2}
Let $L\anabmapright{\alpha} K$ be an anabelomorphism of $p$-adic fields. Then $\alpha$ induces a natural bijective correspondence between Weil representations of $\glt K$ and $\glt L$ respectively.
\ethm
\bp 
From \cite[Theorem 4.8.6]{bump-auto-book} one knows that every 
quadratic extension $K_1\supseteq K$ and a quasicharacter $\chi:K_1^*\to \C^*$ which does not factor through the norm homomorphism $N_{K_1/K}:K_1^*\to K^*$ (i.e. $\chi$ is a character such that if $\tau\in\gal(K_1/K)$ is the unique non-trivial element then $\chi^\tau\neq \chi$) gives rise to an irreducible, supercuspidal representation called the \textit{Weil representation} ${\rm Weil}(K_1/K,\chi)$ of $\glt K$.

Now suppose $\alpha: L\anabelmap K$. Then by the proof of \Cref{th:syn-super-cusp-n}, there exists a unique quadratic field $L_1/L$ and an anabelomorphism $L_1\anabmapright{\alpha_1} K_1$ which is induced by $\alpha$.
By \Cref{th:third-fun-anab}, the anabelomorphism $L_1\anabmapright{\alpha_1} K_1$ provides an isomorphism $\alpha:L_1^*\mapright{\isom} K_1^*$. 

Hence, by composing with $\alpha:L_1^*\mapright{\isom} K_1^*$, a quasicharacter  $\chi:K_1^*\to\C^*$ provides a quasicharacter $L_1^*\to\C^*$.  If $\tau':\gal(L_1/L)$ is the unique non-trivial element then evidently $(\chi\circ\alpha)^{\tau'}\neq \chi\circ\alpha$. Hence, one obtains a Weil representation ${\rm Weil}(L_1/L,\chi')$ where $\chi'=\chi\circ\alpha$. Thus, under anabelomorphy $L\anabelmap K$, one has set up a correspondence $${\rm Weil}(K_1/K,\chi)\mapsto {\rm Weil}(L_1/L,\chi').$$ This procedure is symmetrical in $L$ and $K$, so this establishes the asserted bijection between Weil representations. 
\ep

\bpro\label{pr:central-simple}
Let $K$ be a $p$-adic field and let $K\anabelmap L$ be a choice of an anabelomorphism of $p$-adic fields.
Then one has a natural, dimension preserving bijection between isomorphism classes of finite dimensional central division algebras over $K$ and $L$ respectively.
\epro
\bp
By Class Field Theory \cite[Theorem 1 and Theorem 2]{serre1967}, one has a natural isomorphism of Brauer groups $$Br(K)\isom H^2(\hat{\Z},\Z) \isom Br(L).$$ Hence, the Brauer group $Br(K)$ is amphoric. The asserted correspondence is given by defining $D_L$ to be the unique central division algebra whose isomorphism class coincides with that of $D_K$ in $H^2(\hat{\Z},\Z)$ under the above isomorphism. If $D_K$ has dimension $n^2$, then $[D_K]\in Br(K)$ is an element of order $n$ (\cite[Chapter 13, \ssep 3, Corollary 3]{serre1979-local-fields}) and as $Br(K)=\Q/\Z=Br(L)$, $[D_L]\in Br(L)$ also has order $n$ and dimension $n^2$. This proves the assertion.
\ep

\brem 
Using the above proposition one may hope to synchronize the Jacquet-Langlands correspondence. 
\erem

\subsection{Anabelomorphic Synchronization Theorem for $\GL_2$}
\bthm[Automorphic Synchronization Theorem]\label{th:local-langlands-amphoric}
Let $p$ be an odd prime and let $L\anabelmap K$ be an anabelomorphism of $p$-adic fields. Then this anabelomorphism induces a natural bijection between irreducible admissible representations of $\glt K$ and $\glt L$. This correspondence takes (twists of) irreducible principal series to irreducible principal series, Steinberg to Steinberg and supercuspidal to supercuspidal representations. 
\ethm
\bp 
From \Cref{th:automorphic-ordinary-synchronization} each anabelomorphism $L\anabelmap K$ established a correspondence between  principal series representations of $\glt{K}$ and $\glt{L}$, under which the Steinberg representation of $\glt K$ corresponding to the irreducible sub (resp. quotient) of $\pi(1,\norm{-})$ (resp. $\pi(1,\norm{-}^{-1})$) is mapped to the corresponding object of $\glt L$. By \Cref{th:syn-super-cusp-n}, one sees that this correspondence maps an irreducible supercuspidal representation of $\glt K$ to an irreducible supercuspidal representation of $\glt L$.

Moreover, up to twisting by one dimensional characters, every irreducible admissible representation of $\glt K$ is one of the three types: irreducible principal series representation, a Steinberg representation or a supercuspidal representation. Further,  any twist of an irreducible admissible representation of $\glt K$ is mapped to the corresponding twist of the appropriate irreducible admissible representation. Hence, the assertion is proved.
\ep

For $\GL_2$ and $p\neq 2$ one obtains a fairly complete result:
\bthm[Compatibility of the local Langlands Correspondence]\label{th:anab-local-l-gl2-odd}\  
Let $p$ be an odd prime  and let $L\anabelmap K$ be anabelomorphic $p$-adic fields. Then the local Langlands correspondence for $\glt{K}$ is compatible with the automorphic synchronization provided by \Cref{th:local-langlands-amphoric}.  $L$-functions are amphoric but the conductors and epsilon factors of Weil-Deligne representations and irreducible, admissible representations are not amphoric in general.
\ethm 

\bp 
From \cite[33.1]{bushnell-book} one knows that the local Langlands correspondence for $\glt K$  is a bijection between complex, semisimple, two-dimensional representations of Weil-Deligne group $W_K'$  and irreducible, admissible representations of $\glt K$. This correspondence preserves $\varepsilon$-factors and $L$-functions and maps an irreducible principal series $\pi(\chi_1,\chi_2)$ to $\chi_1\oplus\chi_2$ ($\chi_i$ are considered as characters of $W_K'$ via the Artin map), the Steinberg representation maps to the special representation $sp(2)$ of $W_K'$. A supercuspidal representation $\pi(K_1/K,\chi)$ (\Cref{th:syn-super-cusp-n}) is mapped to the irreducible $W_K'$ representation which is obtained by induction of $\chi$ from $W_{K_1}$ to $W_K$.

Suppose  $\alpha:K\anabelmap L$ is an anabelomorphism. Then \Cref{le:weil-deligne-groups} gives an isomorphism $\alpha:W_K'\isom W_L'$ of Weil-Deligne groups. Given  Weil-Deligne representation $\rho:W_L'\to \gl V$, one can associate to it,  the Weil-Deligne representation $\rho\circ\alpha: W_K' \to \gl V$. This evidently takes semisimple representations to semisimple representations and by  construction, it is compatible with the local Langlands correspondence on both the sides via \Cref{th:local-langlands-amphoric}. Note that the local Langlands correspondence for $\glt K$ matches $L$-functions of representations of $\glt K$ with the $L$-functions associated to $W_K'$ representations. By \Cref{th:char-poly-amphoric}, $L$-functions on the Galois side  (i.e. of $W_K'$ representations) are amphoric. Thus, one deduces that under the correspondence established in \Cref{th:local-langlands-amphoric},  $L$-functions are amphoric.

The simplest way of establishing the assertion about conductors is to work with principal series representations. Suppose $\pi=\pi(\chi_1,\chi_2)$ is a principle series representation of $\glt K$. Then the associated Weil-Deligne representation is a direct sum of two characters of $G_K$ and hence provides two characters of $K^*$. Similarly, by \cite[Th\'eor\`eme 4.1]{deligne1973b}, the $\varepsilon$-factor of $\pi(\chi_1,\chi_2)$ is the product of $\varepsilon$-factors of $\chi_1,\chi_2$.  Thus, to prove the assertion that conductors and $\varepsilon$ of principal series representations are not amphoric in general, it is enough to consider the case of a single character $\chi:K^*\to \C^*$. This is proved in \Cref{th:anab-glo-reps}.  This completes the proof.   
\ep

\subsection{Anabelomorphy and the $p$-adic Langlands correspondence}\label{ss:anab-p-adic-langlands}
In the preceding subsections  we  have not discussed the case of $p$-adic Langlands correspondence. This subsection   outlines our expectations in the $p$-adic case based on \citep*{emerton2025}. 

Let $E$ be a $p$-adic field, $\O_E$ be its ring of integers. The field $E$ (resp. the ring $\O_E$) serve as the coefficient field (resp. coefficient ring) for representations considered here. 
\bcon\label{con:p-adic-langlands1}
Any anabelomorphism $K\anabmapright{\alpha}L$ of $p$-adic fields induces a natural equivalence between the stable $\infty$-categories of complexes of smooth representations of $\gln{K}$ on $p$-torsion $\O_E$-modules and the stable $\infty$-categories of complexes of smooth representations of $\gln{L}$ on $p$-torsion $\O_E$-modules respectively.
In other words, the $\infty$-category of complexes of smooth representations of $\gln K$ on $p$-torsion $\O_E$-modules   is amphoric.
\econ
We remark that \Cref{con:p-adic-langlands1} is true for $n=1$. This is because $K\anabmapright{\alpha}L$ implies, by \Cref{th:third-fun-anab}, that one has an isomorphism of $\glo K=K^*\isom L^*=\glo L$ and hence the corresponding $\infty$-categories of complexes of smooth representations of $\glo K$ (resp. $\glo L$) are naturally equivalent.

\bcon\label{con:p-adic-langlands2}
Let $K$ be a $p$-adic field. For any integer $n\geq 1$, let $\mathcal{X}_{n,K}$ be the formal algebraic stack in \citep*{emerton2025} and occurring in the statement of the categorical Langlands conjecture \citep*[Conjecture 6,1.14]{emerton2025}. Let $K\anabmapright{\alpha}L$ be an anabelomorphism of $p$-adic fields. Then for each $n\geq1$,  the   anabelomorphism $\alpha$ induces an isomorphism of $\O_E$-formal algebraic stacks $\mathcal{X}_{n,K}\isom \mathcal{X}_{n,L}$. In other words, for each $n\geq 1$, the stack $\mathcal{X}_{n,K}$ is amphoric.
\econ
We remark that \Cref{con:p-adic-langlands2} is true for $n=1$ i.e. for $\glo K$ (resp. $\glo L$). This is proved
using the explicit description of the $\O_E$-formal algebraic stacks $\mathcal{X}_{1,K}$ (resp. $\mathcal{X}_{1,L}$) given in \citep*[7.1]{emerton2025}. The asserted isomorphism $$\mathcal{X}_{1,K}\overset{\alpha}{\isom}\mathcal{X}_{1,L}.$$ is obtained by invoking \citep*[7.1.1]{emerton2025} and \Cref{th:third-fun-anab}.

\section{Constructions of varieties via anabelomorphy}\nwss
\subsection{Anabelomorphy and  affine spaces, projective spaces and toric varieties}\label{se:anab-proj}
\newcommand{\ocs}{\O^\circledast}
\newcommand{\ot}{\O^\triangleright}
\newcommand{\ocsk}{\ocs_K}
\newcommand{\otk}{\ot_K}

As noted in \Cref{re:anab-Pn} if $K\anabelmap L$ are anabelomorphic $p$-adic fields, then $\P^n/K$ and $\P^n/L$ are anabelomorphic varieties.    \Cref{th:anab-toric} and \Cref{th:anab-proj-spaces}  show that some topological properties of some anabelomorphic varieties are preserved under anabelomorphisms.

\bthm\label{th:anab-affine-spaces}
Let $\A^n$ (resp. $\Gm$) be the affine space (resp. the multiplicative group) considered as algebraic varieties over a field of choice. Let $K\anabelmap L$ be an anabelomorphism of $p$-adic fields.L et $a,b\geq 0$ be integers, and let $X_K^{a,b}=\A^a\times \Gm^b$ (resp. $X_L^{a,b}=\A^a\times \Gm^b$) considered as an algebraic variety over $K$ (resp.  $L$) with $X^{a,0}=\A^a$ and $X^{0,b}=\Gm^b$. 
Then one has a homeomorphism of topological spaces
$$X^{a,b}(K)=\A^a(K)\times \Gm^b(K) \isom \A^a(L)\times \Gm^b(L)=X^{a,b}(L).$$
In particular for $b=0$, $X^{a,0}/K=\A^a/K$ (resp. $X^{a,0}/L=\A^a/L$) one has a homeomorphism
$$\alpha:\A^a(K)\mapright{\isom}\A^a(L),$$
and for $a=0$ one has a homeomorphism
$$\alpha:\Gm^b(K)=(K^*)^b\mapright{\isom}(L^*)^b=\Gm^b(L).$$
\ethm
\bp 
By \Cref{th:third-fun-anab}{\bf(3)}, one has an isomorphism of topological groups 
$$\Gm(K)=K^*\isom L^*=\Gm(L)$$ and $$(K,+)\isom (L,+).$$ Hence, for any $a\geq 0,b\geq0$,  one has homeomorphisms
$$\A^a(K)=K^a\isom L^a=\A^a(L),$$
and similarly $$\Gm^b(K)=(K^*)^b\isom (L^*)^b=\Gm^b(L).$$
The assertion for $X^{a,b}$ is now clear.
\ep

\subsection{Anabelomorphy, projective spaces and toric varieties}
From \Cref{th:anab-affine-spaces} one obtains the following theorem for smooth, complete toric varieties (\Cref{th:anab-toric}) and projective spaces (\Cref{th:anab-proj-spaces}) over $p$-adic fields.
Both these results were motivated by \cite[Theorem 8.5(iii,iv)]{scholze12-perfectoid-ihes} and are the $p$-adic analog of that result (recalled here as \Cref{th:anab-proj-spaces-perf}). The surprising part of these results is that the $p$-adic fields involved need not be isomorphic.

\bthm\label{th:anab-toric}
Let $X/K=X_{\Sigma}/K$ be a smooth, complete toric variety over a $p$-adic field $K$ obtained from a fan $\Sigma$. Suppose $K\anabmapright{\alpha} L$ is an anabelomorphism. Let $Y$ be the smooth, complete toric variety over $L$ constructed using the fan $\Sigma$. Then one has an anabelomorphism 
$$X/K \anabmapright{\alpha} Y/L$$
and a homeomorphism of  topological spaces
$$X(K)\isom Y(L).$$
Moreover, if $K\anabmapright{\alpha} L$ is a strict anabelomorphism of $p$-adic fields, then $X/K \anabmapright{\alpha} Y/L$ is a strict anabelomorphism of smooth, projective toric varieties.
\ethm

\bp 
By \cite{fulton-toric-varieties}, \cite{danilov1978}, one knows that properties of the fan $\Sigma$ determine the geometric properties of $X_{\Sigma}$ such as smoothness, completeness etc. Hence, $Y/L$ exists and also has these properties. Let $\bK$ (resp. $\bL$) be an algebraic closure of $K$ (resp. $L$). Then by \cite[Theorem 9.1]{danilov1978} $X\times_K\bK$ (resp. $Y\times_L\bL$) is simply connected. This implies that the geometric \'etale fundamental group of $X/K$ is trivial. Hence, $\pi_1(X)\isom G_K$. Similarly $\pi_1(Y)\isom G_L$. Since one has an anabelomorphism $G_K\isom G_L$, one sees that $X/K$ and $Y/L$ are anabelomorphic as asserted.

Now by construction, $X$ is covered by affine opens of the form $X_\sigma$ for cones $\sigma\in \Sigma$. By the smoothness assumption and \cite[Proposition, Page 21]{fulton-toric-varieties}, one sees that
$$X_\sigma\isom  \A^a\times \Gm^b$$
for suitable integers $a,b$ depending on the cone $\sigma$. Thus, one obtains $X_\sigma(K)=\A^a(K)\times \Gm^b(K)$ for each cone $\sigma\in \Sigma$. Hence, one sees, by \Cref{th:anab-affine-spaces}, that one has a homeomorphism $X_\sigma(K)\isom Y_\sigma(L)$ for each $\sigma\in \Sigma$. Now $X$ is obtained from $X_\sigma$ as follows. For cones $\sigma,\tau\in \Sigma$, one has  $X_\sigma\cap X_\tau=X_{\sigma\cap\tau}$ and this is affine and open in both $X_\sigma$ and $X_\tau$, and  $X_\sigma,X_\tau$ are glued along $X_{\sigma\cap\tau}$ (see \cite[Chapter 1, 1.4]{fulton-toric-varieties} for details). Thus, $X(K)$ is obtained by gluing the topological spaces $\{X_\sigma(K):\sigma\in\Sigma\}$ as described (and a similar assertion holds for $Y(L)$) one obtains the asserted homeomorphism $X(K)\isom Y(L)$. 

Now it remains to prove the last assertion. This is proved by the method of proof of \cite{joshi-gconj}. If the anabelomorphism $\alpha$ induces an isomorphism $\beta:X/K\isom Y/L$ of $\Z$-schemes, then $\beta$ induces an isomorphism  $H^0(X,\O_{X})=K\isom L=H^0(Y,\O_{Y})$ of rings. Hence, $K\isom L$ as $p$-adic fields. This contradicts the assumption that $K\anabelmap L$ is strict. This completes the proof.
\ep

\bcor\label{th:anab-proj-spaces}
Let $\alpha:K\anabelmap L$ be an anabelomorphism of $p$-adic fields. Let $n\geq 1$ be an integer. Then a choice of an anabelomorphism $\alpha:K\anabelmap L$ induces  an anabelomorphism $\P^n/K\anabelmap \P^n/L$ of projective varieties and a homeomorphism of associated topological spaces:
$$\alpha:\P^n(K)\mapright{\isom}\P^n(L).$$
\ecor
\bp 
Projectives spaces are toric varieties \cite[Chapter 1]{fulton-toric-varieties}, and so the result follows from \Cref{th:anab-toric}. A direct proof using the familiar construction of projective spaces via gluing affine spaces can also be given using \Cref{th:anab-affine-spaces}.
\ep

\subsection{Anabelomorphy and abelian varieties with split multiplicative reduction}\label{se:anab-tate-abel}
A well-known theorem,  due to John Tate  for dimension one (\cite[Chapter V, Theorem 3.1]{silverman-advanced}) and due to David Mumford  in all dimensions \cite[Corollary 4.9]{mumford72-abelian},  establishes a natural uniformization theorem for abelian varieties  with split multiplicative reduction over valued fields. The treatment of this uniformization theorem for the case of rank one valued fields in \cite[Chapter 6]{fresnel-vanderPut-book2004} is adequate for our next result. 

Let $K$ be a $p$-adic field. Let $\mathbb{G}_{m,K}^{an}$ denote the multiplicative group over $K$ considered as a rigid analytic space over $K$. Write $T_K=(\mathbb{G}_{m,K}^{an})^g$ for the rigid analytic torus of dimension $g$ over $K$. A lattice $\Lambda\subset T_K(K) = (K^*)^g=T_K(K)$ is torsion-free subgroup, isomorphic to $\Z^g$ and  such that the homomorphism $\Lambda\to \R^n$ given by $(x_1,\ldots,x_g)\mapsto(-\log\abs{x_1}_K, \ldots, -\log\abs{x_g}_K)$ is injective and its image is a lattice in $\R^n$ in the usual sense (see \cite[6.4]{fresnel-vanderPut-book2004} for more details).

\bthm\label{th:anabelomorphy-abelian-var} 
Let $K$ be a $p$-adic field and let $A/K$ be a $g$-dimensional $K$-analytic torus given as the (rigid analytic) quotient  $$A_K=T_K/\Lambda_{A_K}$$ by a lattice $\Lambda_{A_K}\subset (K^*)^g=K^*\times \cdots \times K^*$.  For each anabelomorphism $\alpha:G_K\isom G_L$ of $p$-adic fields $K, L$ one has:
\benumlab
\item  a (rigid analytic) torus  $A'_\alpha/L$ given  as the quotient $A_{L,\alpha}'=T_L/\Lambda_{L,\alpha}$, where the lattice $\Lambda_{L,\alpha}=\alpha(\Lambda_{A_K})\subset (L^*)^g =T_L(L)$ is the image of $\Lambda_{A_K}$  under the functorial isomorphism $\alpha:(K^*)^g \mapright{\isom} (L^*)^g$ given by the amphoricity of $K^*$, and one has a 
homeomorphism of topological groups
$$f_\alpha:A_K(K) = (K^*)^g/\Lambda_{A_K}\isom (L^*)^g/\Lambda_{L,\alpha}=A_{L,\alpha}'(L);$$  
\item moreover, the construction of the torus $A_\alpha'/L$,  the lattice parameter $\Lambda_{L,\alpha}$ and the homeomorphism $f_\alpha$, are all functorial in all the variables $L,\alpha$ and independent of the choice of the lattice $\Lambda_K$ giving rise to $A_K$.
\item If $A_K/K$ is an abelian variety over $K$, then  $A'_{L,\alpha}$ an abelian variety over $L$.
\eenum
\ethm
\bp 
The assertions {\bf(1), (2)} are clear. Let us prove the remaining assertion. This is done using \cite[Theorem 6.6.1]{fresnel-vanderPut-book2004}. Let $T_K$ be the rigid analytic torus with $T_K(K)=(K^*)^g$. Giving $T_K$ is equivalent to giving its character group i.e. giving a free $\Z$-module $X(T_K)$ of rank $g$, equipped with a continuous action of $G_K$. The anabelomorphism $\alpha^{-1}:G_L\isom G_K$ allows us to view this as a free $\Z$-module of rank $g$ equipped with a continuous action of $G_L$. This module is the character group of a torus $T_{L,\alpha}$ with $T_{L,\alpha}(L)=(L^*)^g$. Since $A_K$ is an abelian variety, by \cite[Theorem 6.6.1]{fresnel-vanderPut-book2004}, there exists a homomorphism 
$$\sigma: \Lambda_{A_K}\to X(T_K)$$
such that (a) $\sigma(\lambda)(\lambda')=\sigma(\lambda')(\lambda)$ for all $\lambda,\lambda'\in \Lambda_{A_K}$ and (b) the bilinear form $\langle \lambda,\lambda'\rangle=-\log\abs{\sigma(\lambda')(\lambda)}$ is positive definite.
These properties are unaltered by the isomorphism $\alpha:\Lambda_{A_K} \isom \Lambda_{L,\alpha}$ and the anabelomorphism $\alpha^{-1}:G_L\to G_K$ which provides the $G_L$-module structure on the $G_K$-module $X(T_K)$. Thus, by \cite[Theorem 6.6.1]{fresnel-vanderPut-book2004}, $A'_{L,\alpha}$ is an abelian variety.
\ep

The following corollary is immediate:
\bcor\label{co:anab- ab.-var}
In the notation and the hypothesis of \Cref{th:anabelomorphy-abelian-var}, one has an isomorphism of \'etale fundamental groups:
$$\pi_1(A)\isom \pi_1(A'_{L,\alpha}),$$
which is functorial in $\alpha$. 
In other words, $A/K$ and $A'_{L,\alpha}$ are anabelomorphic abelian varieties.
\ecor
\bp 
\newcommand{\bkh}{\widehat{\bK}}
\newcommand{\blh}{\widehat{\bL}}
\newcommand{\bA}{\bar{A}}
\newcommand{\aan}{{A}^{an}}
\newcommand{\Zh}{\hat{\Z}}
Let $g=\dim(A)$ and suppose that $\ell\neq p$ is a prime number.  
The asserted isomorphism will be first established for tempered fundamental groups. Since  \'etale fundamental groups are profinite completions of the respective tempered fundamental groups \cite[Proposition 4.4.1]{andre03}, one obtains the stated isomorphism.

Let $\bK\supset K$ (resp. $\bL\supset L$) be algebraic closures of $K$ and $L$ respectively. Let $\C_K$ (resp. $\C_L$) be the completion of $\bK$ (resp. $\bL$). For computing tempered fundamental groups, we will use geometric  base-points with values in $\C_K$ (resp. $\C_L$) and let $\mathcal{M}(K)$ (resp. $\mathcal{M}(\C_K)$) be the Berkovich spectrum of $K$ (resp. $\C_K$). Let $\bA^{an}=A^{an}_K\times_{\mathcal{M}(K)} \mathcal{M}(\C_K)$ be the base extension of the analytic space $\aan_K$ associated to $A/K$. For notational simplicity, write $B=A'_{L,\alpha}$ and $\bar{B}^{an}=B^{an}\times_{\mathcal{M}(L)}\mathcal{M}(\C_L)$.

Since one has a rigid analytic quotient isomorphism $\C_K^{*g}/\Lambda_{K}\to \bA$,  by \cite[2.4.1]{lepage-thesis} one has the following description of the geometric tempered fundamental group of $A$ i.e. of the tempered fundamental group of $\aan_K$:
$$\pi_1^{temp}(\bA)=\Lambda_K\times \hat{\Lambda}_{K}(1)^g \isom \Z^g\times \Zh(1)^g_K$$
where $G_K\act \Zh(1)_K$ is the Galois module of the roots of unity contained in $K$ i.e. the free  $\Zh$-module of rank $g$ with an action by the cyclotomic character of $G_K$ and $\hat{\Lambda}_{K}$ is the profinite completion of $\Lambda_{K}$. Moreover, the tempered fundamental group of $A/K$ fits in the exact sequence
$$0\to\Z^g\times\Zh(1)^g=\pi_1^{temp}(\bA)\to \pi_1^{temp}(A)\to G_K \to 1.$$
This exact sequence splits (using the $K$-rational point given by the identity element of $A(K)$). 

The abelian variety $B$ also provides a similar sequence over $L$. Now the asserted isomorphism $\pi_1^{temp}(A/K)\isom \pi_1^{temp}(B/L)$ follows from the following proven facts regarding these objects (1) the construction of $\Lambda_{L,\alpha}$ from $\Lambda_K$ (\Cref{th:anabelomorphy-abelian-var}) (2) the construction of $B$ from $L^{*g}/\Lambda_{L,\alpha}$ (\Cref{th:anabelomorphy-abelian-var})  (3) the isomorphism $A(K)\isom B(L)$ given by \Cref{th:anabelomorphy-abelian-var} which being a homeomorphism, maps the identity element of $A(K)$ to that of $B(L)$, (4) the splitting of the exact sequence for $\pi_1^{temp}(B/L)$, and (5) the amphoricity of the cyclotomic character of $G_K$ (\Cref{th:third-fun-anab}{\bf(4)}).

\ep

For $g=1$,   a lattice $\Lambda\subset \mathbb{G}_{m,K}^{an}$ is given by $\Lambda=(q_K^\Z)\subset K^*$ with $q_K\in K^*$ which one takes to satisfy $\abs{q_K}_K<1$ and \Cref{th:anabelomorphy-abelian-var} gives the Tate elliptic curve (\cite[Chapter V, Theorem 3.1]{silverman-advanced}):

\bcor\label{th:anabelomorphy-tate-curves}  
Let $K$ be a $p$-adic field and let $E/K$ be a Tate elliptic curve over $K$ with Tate parameter $q_K\in K^*$.
Let $\alpha:G_K\isom G_L$ be an anabelomorphism of $p$-adic fields and let $\alpha: K^*\to L^*$ be the functorial isomorphism given the anabelomorphism $\alpha$. Then there exists a Tate elliptic curve $E'_\alpha/L$ with Tate parameter $q_{L,\alpha}=\alpha(q_K)$ and a
homeomorphism of topological groups
$$f_\alpha:E(K)\isom E_\alpha'(L).$$  The construction of  $E_\alpha'/L$,  the Tate parameter $q_{L,\alpha}$ and the homeomorphism $f_\alpha$, are all functorial in $L,\alpha$ and independent of the choice of  $q_K$.  Explicitly, 
the elliptic curve $E_\alpha'/L$ is given by Tate's equation
$$y^2+xy=x^3+a_4(q_{L,\alpha})x+a_6(q_{L,\alpha}).$$
\ecor
\bp 
All the assertions are immediate from \Cref{th:anabelomorphy-abelian-var}. That $\abs{q_L}_L<1$ follows from \Cref{le:norm-amphoric}{\bf(3)} and the explicit formula for the equation of the Tate elliptic curve which is given by \cite[Chap. V, Theorem 3.1]{silverman-advanced}.
\ep

\subsection{Anabelomorphy of finite, flat group schemes of order $p$ over $p$-adic fields}
As a warm-up to the main result of the next section, let us establish the following:
\bthm\label{th:grp-ord-p}
Let $K,L$ be anabelomorphic $p$-adic fields. Then each anabelomorphism $\alpha:G_L\isom G_K$ provides a natural bijection between isomorphism classes of finite flat group schemes of order $p$ over $\O_K$ and $\O_L$ respectively.
\ethm
\newcommand{\okt}{\O_K^\triangleright}
\newcommand{\olt}{\O_L^\triangleright}
\bp 
This will be proved using anabelomorphy and the classification theorem of finite flat group scheme of order $p$ proved in \cite{oort1970}.  Let $\okt$ (resp. $\olt$) be the multiplicative monoid of the ring $\O_K$ (resp. $\O_L$). Then by  \cite[Proposition 3.11, Summary 3.15]{hoshi-mono}, the monoid $\okt$ is amphoric and the anabelomorphism $\alpha$ induces an isomorphism of topological monoids $\alpha:\olt\mapright{\isom}\okt$ which takes $p\in\olt$ to $p\in\okt$ and moreover, this is compatible with the inclusion of the respective unit subgroups  $\O_L^*, \O_K^*$ and the isomorphism $\alpha: \O_L^*\isom \O_K^*$ provided by the amphoricity of $\O_K^*$ (\Cref{th:third-fun-anab}). 

By \cite[Theorem 2 and Remark (5), Pages 15--16]{oort1970}, for each pair of elements $a,b\in\olt$ satisfying $a\cdot b=p$, one has a finite, flat group scheme $\sG_a^b$  of order $p$ over $\O_L$. Moreover, if $\sG_{c}^d$ (with $c\cdot d=p$) is another finite group scheme of order $p$ over $\O_L$, then $\sG_a^b$, $\sG_{c}^d$ are $\O_L$-isomorphic if and only if there exists a unit $u\in\O_L^*$ such that $c=u^{p-1}\cdot a$ and $d=u^{1-p}\cdot b$. 

Now these data and the relationship between them is preserved by the isomorphism of multiplicative monoids $\olt\mapright{\alpha}\okt$ (compatibly with the isomorphism $\O_L^*\mapright{\alpha}\O_K^*$). Thus, if one writes $a'=\alpha(a), b'=\alpha(b)$, then one has $a'\cdot b'=\alpha(a)\cdot \alpha(b)=\alpha(a\cdot b)=\alpha(p)=p$ and hence one has a finite flat group scheme $\sH_{a'}^{b'}$ over $\O_K$ of order $p$. Clearly, under this mapping $\sG_a^b\mapsto \sH_{a'}^{b'}$,  any $\O_L$-group scheme isomorphic to $\sG_a^b$ is mapped to an $\O_K$-group scheme  of order $p$ isomorphic to $\sH_{a'}^{b'}$. This proves the theorem.
\ep

\subsection{Anabelomorphy of $\F_q$-vector space schemes over $p$-adic fields}\label{se-group-schemes}
There is a variant of \Cref{th:grp-ord-p} based on \cite{raynaud74}.  By \cite[D\'efinition 1.2.1]{raynaud74}, an  $\F_q$-vector space scheme over a base scheme $S$ is a contravariant, representable functor from the category of $S$-schemes to the category of $\F_q$-vector spaces. All $\F_q$-vector space schemes discussed here are assumed to be finite, flat and of finite presentation over the relevant base scheme $S$. An $\F_q$-vector space scheme over $S$ is thus a finite, flat, commutative  group scheme (of finite presentation) which is annhilated by multiplication by $p$. 

For a $p$-adic field $K$, write  $\Q_p\subset K_0\subset K$ for its maximal unramified subfield and write $r=[K_0:\Q_p]$ (choice of letter $r$ instead of the conventional $f$ for this number is for compatibility with \cite{raynaud74}). Then the residue field $\F_q$ of $K_0$ has cardinality $q=p^r$. For the remainder of the section, let $q=p^r$ for this choice of $r$. In this section, one is interested in $\F_q$-vector space schemes over the base scheme $\Spec(\O_K)$ where $K$ is a $p$-adic field and $\O_K$ is the ring of integers of $K$. 

\bthm\label{th:group-schemes} Let $K,L$ be anabelomorphic $p$-adic fields. Let $K\supset K_0\supset\Q_p$ be the maximal unramified subfield of $K$.  Let $q=p^r$ where $r=[K_0:\Q_p]$. Then any anabelomorphism $\sigma:G_L\mapright{\isom} G_K$ provides a natural bijection between isomorphism classes of $\F_q$-vector space schemes of rank one over $\O_L$ and $\O_K$  respectively.
\ethm
\brem 
Note that \Cref{th:group-schemes} does not imply \Cref{th:grp-ord-p} because for a general $p$-adic field one has $r>1$.
\erem
\bp  
Let $L_0\subset L$ (resp. $K_0\subset K$) be the maximal unramified subfield of $L$ (resp. $K$). By \Cref{th:third-fun-anab}, $[K_0:\Q_p]$ is an amphoric quantity.  Then $[L_0:\Q_p]=r=[K_0:\Q_p]$. Hence, $q=p^r$ is an amphoric quantity.
By \cite[Proposition 3.11(iii)]{hoshi-mono}, any anabelomorphism $K\anabelmap L$, provides a natural isomorphism between the multiplicative monoids of non-zero elements of the residue fields of $K$ and $L$ respectively (these monoids are groups and the proof of the cited proposition shows that this isomorphism of monoids (each identified with $\F_q^*$) is a natural isomorphism of groups). 

Let $D'=\Z[\zeta_{q-1}]$ where $\zeta_{q-1}$ is a primitive $(q-1)^{th}$-root of unity in some algebraic closure of $\Q$. Choose a generator for the cyclic group $\F_q^*$. One has a ring homomorphism $D'\to \O_{K}$ which maps $\zeta_{q-1}$ to the Teichm\"uller lift, in $\O_{K_0}$, of the chosen generator of the cyclic group $\F_q^*$.  Since $q-1$ is invertible in $\O_{K_0}\subset \O_K$, this homomorphism factors through $\Z[\zeta_{q-1},\frac{1}{q-1}]$ and singles out a unique prime ideal $\wp$ of $D'$ lying over $p$, and further factors through the ring  $D\subset \Q(\zeta_{q-1})$ defined in \cite[Section 1.1]{raynaud74}, in which $q-1$ is a unit and consists of $x\in\Q(\zeta_{q-1})$ which are $\wp'$-integral for all $\wp'|p$ except possibly at $\wp$. Thus, one has $\Z[\zeta_{q-1}]\into \Z[\zeta_{q-1},\frac{1}{q-1}]\into D$. Similar assertion holds for $L_0$ (the chosen generator of $\F_q^*$ is mapped to its image under the isomorphism of multiplicative monoids of the residue fields provided in the previous paragraph). Hence, one has a ring homomorphism $D\to \O_{L_0}$. Note that the construction of $D',D$ is independent of the fields $K,L$. Thus one can view $\O_K,\O_L$ as $D$-algebras. 

Note that since $\O_K,\O_L$ are Noetherian, complete local rings of characteristic zero,  Raynaud's condition $(**)$ \cite[Page 246]{raynaud74} holds by \cite[Proposition 1.2.2]{raynaud74}  for the sort of group schemes being considered in this theorem.

\newcommand{\bfa}{\boldsymbol{\alpha}}
\newcommand{\bfb}{\boldsymbol{\beta}}
It will be convenient to reformulate Raynaud's result in the style of \cite{oort1970} discussed above. For this purpose note that the elements $w,u\in D$ satisfying $w=p\cdot u\in D$, with $u\in D^*$ being a unit, and defined by \cite[Equation (17), Proposition 1.3.1]{raynaud74} are independent of $K,L$. We claim that  given a collection \be\bfa= \left\{(a_i,b_i): a_i\cdot b_i=p \text{ and } a_i,b_i\in \O_L  \right\}_{0\leq i\leq r-1},\ee (with the convention that index $i$ is read as $i\bmod(r)$), there exists an $\F_q$-vector space scheme $\sG(\bfa)$ over $\O_L$. Indeed, given such a system of elements $\bfa$, writing $\gamma_i=a_i$ and $\delta_i=b_i\cdot u$, one obtains a Raynaud system 
\be\left\{(\gamma_i,\delta_i): \gamma_i\cdot \delta_i=p\cdot u=w \text{ with } \gamma_i,\delta_i\in \O_L  \right\}_{1\leq i\leq r},\ee
to which one may apply \cite[Corollaire 1.5.1]{raynaud74}, to obtain the claimed group scheme $\sG(\bfa)$ over $\O_L$.

Now suppose 
\be \bfb= \left\{(a_i',b_i'): a_i'\cdot b_i'=p \text{ and } a_i',b_i'\in \O_L  \right\}_{0\leq i\leq r-1},
\ee is another system giving the $\F_q$-vector space scheme $\sG(\bfb)$. Then  from the relations given \cite[Corollaire 1.51]{raynaud74} one sees that, $\sG(\bfa)\isom \sG(\bfb)$ as  $\F_q$-vector space schemes over $\O_L$ if and only if there exists a system of units $\{ u_i\in \O_L^*\}_{0\leq i\leq r-1}$ such that
\begin{align*}
a_i'&=u_{i+1}\cdot a_i\cdot u_i^{-p}\\
b_i'&=u_i^p\cdot b_i \cdot u_{i+1}^{-1}.
\end{align*}

Now suppose $\sigma:G_L\mapright{\isom} G_K$ is any anabelomorphism of $p$-adic fields. Let $\sigma:\olt\mapright{\isom}\okt$ be the isomorphism provided by amphoricity of $\olt$. As remarked in the context of the proof of \Cref{th:grp-ord-p}, this is compatible with the isomorphism $\sigma:\O_L^*\mapright{\sigma}\O_K^*$ provided by the amphoricity of $\O_L^*$. The system $\bfa$ is a system of elements of $\olt$ and applying $\sigma$ gives a system $\bfa'=\sigma(\bfa)$ in $\okt$ and hence, by \cite[Corollaire 1.5.1]{raynaud74}, an $\O_K$-group scheme $\sH(\bfa')$ which is an $\F_q$-vector space scheme and any $\F_q$-vector space scheme $\sG(\bfb)$ which is isomorphic to $\sG(\bfa)$ is mapped isomorphically to an $\F_q$-vector space scheme $\sH(\bfb')$ isomorphic to $\sH(\bfa')$. This proves the theorem.
\ep

\section{Anabelomorphic Connectivity Theorem for Number Fields}\label{se:anab-connect-thm}\nwss
The notion of anabelomorphy suggests the possibility of anabelomorphically modifying a number field at a finite number of places to create another number field which is anabelomorphically glued to the original one at a finite number of places. Anabelomorphic connectivity theorems  provide a way of passing geometric information between two such connected fields. This is the main theme of this section.

\subsection{Definition and examples}
\bdefn\label{def:ef:anab-connected} 
We say that two number fields $K,K'$ are \emph{anabelomorphically connected along non-archimedean places $\vseq v n$ of $K$ and $\vseq w n$ of $K'$} if,  for each $i=1,\ldots,n$, there exists an anabelomorphism $K_{v_i}\anabelmap K'_{w_i}$. We will simply denote this as $$\anab{K}{K'}{v_1,\ldots,v_n}{w_1,\ldots,w_n}.$$
If any of the anabelomorphisms $K_{v_i}\anabelmap K'_{w_i}$ is a strict anabelomorphism, then one says that $K,K'$ are \emph{strictly anabelomorphically connected along $\vseq v n$ and $\vseq w n$.}
\edefn

\begin{example}\label{ex:basic-ex2}
Here is a basic collection of examples for \Cref{def:ef:anab-connected}. Let $p$ be an odd prime, let $r\geq 1$ be an integer. Let $K_r=\Q(\zeta_{p^r},\sqrt[{p^r}]{p})$, $K_r'=\Q(\zeta_{p^r},\sqrt[{p^r}]{1+p})$. These are totally ramified at $p$ (see \cite[Theorem 5.5]{viviani04}). Let $\wp$ (resp. $\wp'$) be the unique prime of $K_r$ prime lying over $p$ (resp. the unique prime of $K'_r$ lying over $p$). The completions of $K_r$ (resp. $K_r'$) with respect to these unique primes are $K_{r,\wp}=\Q_p(\zeta_{p^r},\sqrt[{p^r}]{p})$ and $K'_{r,\wp'}=\Q_p(\zeta_{p^r},\sqrt[{p^r}]{p})$ respectively. By Lemma~\ref{le:basic-example}, one has a (strict) anabelomorphism $$K_{r,\wp}\anabelmap K'_{r,\wp'}.$$  
In particular, for any $r\geq1$, the number fields  $K_r=\Q(\zeta_{p^r},\sqrt[{p^r}]{p}), K'_r=\Q(\zeta_{p^r},\sqrt[{p^r}]{1+p})$ (and the unique primes $\wp_r,\wp_r'$ lying over $p$ in $K_r,K'_r$) are (strictly) anabelomorphically connected along $\wp_r$ and $\wp_r'$:
	$$\anab{K_r}{K'_r}{\wp_r}{\wp'_r}.$$ 
\end{example}

\subsection{Existence of anabelomorphically connected number fields}
The next step is to establish (in \Cref{th:anab-connectivity-thm2}) the existence of strictly anabelomorphically connected number fields. This provides a systematic  way of producing examples of anabelomorphically connected number fields starting with a given number field. 

In what follows, we will say that a number field $M$ is \emph{dense} in a $p$-adic field $L$ if there exists a place $v$ of $M$ such that the completion $M_v$ of $M$ at $v$ is isomorphic to $L$.

\newcommand{\nonarch}{non-archimedean}
We will use the following terminology: a \emph{non-archimedean local field} is a finite extension of $\Q_p$ for some (unspecified) prime  $p$. 
\begin{defn}
	We say that a non-empty finite set  of \nonarch\ local fields $\{L_1,\cdots,L_n\}$ (some of which may be pairwise isomorphic and some may have distinct residue characteristics) is a \emph{cohesive set of \nonarch\ local fields} if there exists a number field $M$ and an inclusion $M\into L_i$ which is dense for all $i$,  such that the induced valuations on $M$ are pairwise inequivalent.	
\end{defn}

\blem[Potential Cohesivity Lemma]\label{le:coh}
For every non-empty finite set $\{L_1,\cdots,L_n\}$ of \nonarch\ local fields, some of which may be pairwise isomorphic and some may have pairwise distinct residue characteristics, there exist finite extensions $L_i'/L_i$ such that $\{L_1',\cdots,L_n'\}$ is a cohesive system of \nonarch\ local fields.
\elem

\bp
By Krasner's Lemma (\cite[Chapter 3, Section 3]{koblitz-p-adic-book}) every non-archimedean field contains a dense number field and so the result is true for $n=1$ on taking $L_1'=L_1$. The general case will be proved by induction on $n$. Suppose that the result has been established for the case of $n-1$ fields with $n>1$. So for every set $L_1,\ldots,L_{n-1}$ of \nonarch\ fields there exists finite extensions $L_1',\ldots,L'_{n-1}$ of \nonarch\ fields and a number field $M\subset L_i'$ which is  dense for $i=1,\ldots,n-1$ and the valuations induced on $M$ are all inequivalent. By the primitive element theorem, one can choose $\alpha\in M$ such that $\Q(\alpha)=M$.

Now suppose that $p$ is the residue characteristic of $L_n$. By Krasner's Lemma one can choose $\beta\in L_n$ to be algebraic and such that $L_n=\Q_p(\beta)$. Then $\Q(\beta)\subset L_n$ is a dense inclusion of a number field in $L_n$. Now consider the finite extensions $L_n'=L_n(\alpha)$ and $L_i''=L_i'(\beta)$ (if the minimal polynomial of $\alpha$ over $\Q$ is not irreducible over $L_n$, then pick a direct factor of $L_n\tensor \Q(\alpha)$, as this is a product of fields each of which is a finite extension of $L_n$ equipped with an embedding of $\Q(\alpha)$, and similarly for $\beta$, for $i=1,\ldots,n-1$). Then $\Q(\alpha,\beta)\subset L_i''$ for $i=1,\ldots,n-1$ and $\Q(\alpha,\beta)\subset L_n'$. Write $L_n''=L_n'$ (for symmetry of notation). Then one sees that there exists a common number field $M$ contained in all the $L_i''$. If $M$ is not dense in each of $L_i''$ one can extend $M$ further to achieve density. Similarly, if the induced valuations on $M$ are not all inequivalent, one can extend $M$ further to achieve this as well. Let us explain how these extensions in the last two steps are carried out. 

To avoid notational chaos, we will prove both these assertions for $n=2$. So the situation is that one has  two \nonarch\ fields $L_1,L_2$ and a common number field $M$ contained in both of them. There are two possibilities: either residue characteristics of $L_1, L_2$ are equal or they are not equal. First assume that the residue characteristics are equal (say equal to $p$). Then $L_1,L_2$ are both finite extensions of $\Q_p$ and so there exists a finite extension $L$ containing both of them as subfields. Pick such an $L$. Then there is a number field $M'$ dense in $L$. Now choose a number field $F$, with $[F:\Q]>1$, which is totally split at $p$ and such that $M',F$ are linearly disjoint over $\Q$. Then let $M''=MF\into L$ and since $F$ is completely split there exist two primes $v_1\neq v_2$ of $M''$ lying over $p$ such that $M''_{v_1}=L$ and $M''_{v_2}=L$. Thus the system $L_1=L, L_2=L$ is now cohesive as $M''\into L_1=L$ and $M''\into L_2=L$ are dense inclusions corresponding to distinct primes of $M''$.

Now assume $L_1,L_2$ have distinct residue characteristics and $M$ is a number field contained in both of them.  If $v_1$ (resp. $v_2$) is the prime of $M$ corresponding to the inclusion $M\into L_1$ (resp. $M\into L_2$), then $M_{v_1}\into L_1$ and $M_{v_2}\into L_2$ are finite extensions of \nonarch\ fields. One proceeds by descending induction on the degrees $[L_1:M_{v_1}], [L_2:M_{v_2}]$. By the primitive element theorem there exists an $x_1\in L_1$ (resp. $x_2\in L_2$) such that $L_1=M_{v_1}(x_1)$  (resp. $L_2=M_{v_2}(x_2)$). Choose an irreducible polynomial $f\in M[X]$ which is sufficiently close to the minimal polynomials of $x_1$ (resp. $x_2$) in $L_1[X]$ and $L_2[X]$ respectively. Then $f$ has a root in both $L_1,L_2$ (by Krasner's Lemma). The field $M'=M[X]/(f)$ embeds in both $L_1,L_2$ and if $v_1'$ (resp. $v_2'$) is the prime lying over $v_1$ (resp. $v_2$) corresponding to the  inclusion $M'\into L_1$ and $M'\into L_2$ are dense inclusions of $M'$ in $M'_{v_1}\subset L_1$ (resp. $M$ in $M'_{v_2}\into L_2$)  and  $[L_1:M'_{v_1'}]<[L_1:M_{v_1}]$ and similarly for $L_2$. Thus by enlarging $M$ in this fashion one is eventually led to a cohesive system as claimed.
\ep

Now we can state and prove the general anabelomorphic connectivity theorem for number fields.

\bthm[Anabelomorphic Connectivity Theorem]\label{th:anab-connectivity-thm2}
Let $K$ be a number field. Let $\vseq{v}{n}$ be a finite set of non-archimedean places of $K$. Let  $\alpha_i:K_{v_i}\anabelmap L_i$ be arbitrary anabelomorphisms with \nonarch\ local fields $L_1,\ldots,L_n$. 
Then there exist
\benumlab
\item   finite extensions $L_i'/L_i$ (for all i) and a dense embedding of a number field $M'\subset L_i'$ and places $\vseq{w}{n}$ of $M'$ induced by the embeddings $M'\into L_i'$ (i.e. the collection $\{L_i'\}$ of \nonarch\ fields is cohesive via $M'$) and 
\item  a finite extension $K'/K$ and, for all $i$, places $\vseq{u}{n}$ of $K'$ lying over the places $v_i$ of $K$ together with anabelomorphisms $K'_{u_i}\anabelmap L_i'$.
\item In particular,   $\anab{K'}{M'}{\vseq{u}{n}}{\vseq{w}{n}}$ and $u_i|v_i$ for all $i=1,\ldots,n$.
\eenum
\ethm
\bp 
Fix an algebraic closure $\bK$ of $K$ and an algebraic closure $\bK_{v_i}$ of $K_{v_i}$ for each $i$  so that $G_K=\gal(\bK/K)$ and $G_{K_{v_i}}=\gal(\bK_{v_i}/K_{v_i})$ and fix embeddings $K\into K_{v_i}$.
By the Cohesivity Lemma (Lemma~\ref{le:coh}) one can replace $L_1,\ldots,L_n$ by a cohesive collection $L_1',\ldots,L_n'$ with $L_i'/L_i$ finite extensions and a number field $M'\subset L_i'$ dense in each $L_i'$ such that the induced valuations on $M'$ are all inequivalent. The finite extensions $L_i'/L_i$ provide open subgroups $H_i'\subset  G_{L_i}$ of $G_{L_i}$. Since one has anabelomorphisms $\alpha_i:K_{v_i}\anabelmap L_i$, let $H_i=\alpha^{-1}(H_i')$ be the inverse image of $H_i'$ in $G_{v_i}$. Since $\alpha_i$ is continuous, $H_i'$ is an open subgroup of $G_{v_i}$ (for each $i$). Let $\bK_{v_i}\supset F_i\supseteq K_{v_i}$ be the finite extension of $K_{v_i}$ corresponding to $H_i'$. 

Then we claim that there exists a finite extension $K'/K$ and primes $u_i|v_i$ of $K'$ such that $K'_{u_i}=F_i$. This is seen as follows. By the primitive element theorem \cite[Chap. V, \ssep 4, Theorem 4.6]{lang-algebra}, one can assume $F_i=K_{v_i}(\gamma_i)$. Let $f_i(X)\in K_{v_i}[X]$ be the monic minimal polynomial of $\gamma_i$. Then by the weak approximation theorem \cite[Chapter 1, 1.2.2, Theorem 1.4]{platonov-book}, there exists a polynomial $f(X)\in K[X]$ which is arbitrarily close to $f_i(X)$ at $v_i$. Then there exists a root $\alpha\in \bK$ of $f(X)=0$ such that the equality $K_{v_i}(\alpha)= K_{v_i}(\gamma_i)$ (in $\bK_{v_i}$) holds by Krasner's Lemma \cite[Chap II, \ssep 2, Proposition 4]{Lang1970}. Let $K'=K(\alpha)$. Thus, by construction one has $K_{v_i}(\alpha)= K_{v_i}(\gamma_i)=  F_i$ and $K'\into K_{v_i}(\gamma_i)=F_i$  is dense for each $i$. Hence, there exists primes $u_i|v_i$ of $K'$ such that $K'_{u_i}=K_{v_i}(\gamma_i)=F_i$ and $G_{K'_{u_i}}=H_i'$ as claimed,
and  $$G_{K'_{u_i}}\isom H_i\isom H'_i\isom G_{L_i'}\isom G_{M_{w_i}}.$$ Hence, one has established that $$\anab{K'}{M'}{\vseq{u}{n}}{\vseq{v}{n}}.$$ This completes the proof.
\ep

\newcommand{\gm}{\mathfrak{m}}
\newcommand{\fp}{\mathfrak{p}}
\newcommand{\wq}{\mathfrak{q}}

\bthm\label{th:anab-con-inf2} 
Let $\anab{K}{M}{\vseq{u}{n}}{\vseq{v}{n}}$ be anabelomorphically connected number fields. Then there exists  anabelomorphically connected number fields $$(K,\{\vseq{u}{n}\})\anabelmap\anab{M}{M'}{\vseq{v}{n}}{\vseq{w}{n}}$$
such that $\deg(M')>\deg(M)$ and, if $\anab{K}{M}{\vseq{u}{n}}{\vseq{v}{n}}$ is a strict anabelomorphic connectivity, then so is $\anab{K}{M'}{\vseq{u}{n}}{\vseq{w}{n}}$. In particular, the class of number fields which are (strictly) anabelomorphically connected with $(K,\{\vseq{u}{n}\})$ is infinite and degree is unbounded in this class.
\ethm
\bp

For $1\leq i\leq n$, let $p_i$ be the common residue characteristic of $K_{u_i}$ and $M_{v_i}$ (the primes $p_1,\ldots,p_n$ may not be all pairwise distinct). Let $F$ be a quadratic field such that $F/\Q$ is totally split at all the pairwise distinct primes from among $p_1,\ldots,p_n$ and also totally split at all primes which are ramified in $M/\Q$. Then $F\cap M$ has no ramified primes and hence $F\cap M=\Q$. Let $M'=MF$ be the compositum. Then by construction $M'/M$ is totally split at $v_1,\ldots,v_n$. For $1\leq i\leq n$, let $w_i$ be a prime of $M'$ lying over $v_i$. Then for each $i$, one has an isomorphism $M_{v_i}\isom M'_{w_i}$ of $p_i$-adic fields. Hence, one has anabelomorphisms $K_{u_i}\anabelmap M_{v_i}\anabelmap M'_{w_i}$ and if $K_{u_i}\anabelmap M_{v_i}$ is a strict anabelomorphism for some $i$, then so is $K_{u_i}\anabelmap M'_{w_i}$. Thus, one has anabelomorphically connected number fields $\anab{M}{M'}{\vseq{v}{n}}{\vseq{w}{n}}$ and hence 
$(K,\vseq{u}{n})\anabelmap\anab{M}{M'}{\vseq{u}{n}}{\vseq{w}{n}}$. If $\anab{K}{M}{\vseq{u}{n}}{\vseq{v}{n}}$ is a strict anabelomorphic connectivity then so is $\anab{K}{M'}{\vseq{u}{n}}{\vseq{w}{n}}$. Moreover, by construction $\deg(M')>\deg(M)$. This proves the theorem.
\ep

\subsection{The Ordinary Synchronization Theorem}\label{se:ordinary-syn}
The next result is consequence of \Cref{th:anab-connectivity-thm2} and \Cref{th:ordinary-amphoric}. Special cases of this result, in which anabelomrophic connectivity arises from isomorphisms of $p$-adic fields, can be found scattered in the literature on automorphy of Galois representations. The theorem takes its name from ``Synchronization of Geometric Cyclotomes'' discovered by Mochizuki for e.g. see \cite[Section 6]{hoshi-mono}.     The theorem is the following:
\bthm[The Ordinary Synchronization Theorem]\label{th:ord-syn-thm}
Let $$\anab{K}{K'}{\vseq{v}{n}}{\vseq{w}{n}}$$ be a pair of anabelomorphically connected number fields. Let $p_i$ be the common residue characteristic of $v_i,w_i$ for $1\leq i\leq n$. Then one has for all primes $\ell$ (including $p$) and for all $i$:
\benumlab
\item an equivalence of categories of ordinary $\ell$-adic ($\ell=p_i$ included) 
$G_{K'_{w_i}}$- and $G_{K_{v_i}}$-representations respectively;
\item for two-dimensional ordinary representations, one has 
isomorphisms of $\Q_\ell$-vector spaces (including $\ell=p_i$)
$$H^1({G_{K_{v_i}}},\Q_{\ell}(1))\isom\Ext^1_{G_{K_{v_i}}}(\Q_{\ell}(0),\Q_{\ell}(1)) \isom \Ext^1_{G_{K'_{w_i}}}(\Q_{\ell}(0),\Q_\ell(1))\isom H^1(G_{K'_{w_i}},\Q_\ell(1));$$
\item and for ordinary crystalline, two-dimensional $p$-adic representations, an isomorphism
$$H^1_f(G_{K_{v_i}},\Q_{p_i}(1)) \isom H_f^1(G_{K'_{w_i}},\Q_{p_i}(1)),$$
\item and also an isomorphism 
$$H^1_e(G_{K_{v_i}},\Q_{p_i}(1)) \isom H_e^1(G_{K'_{w_i}},\Q_{p_i}(1)).$$
\eenum
\ethm

\bp 
The first assertion is immediate from \Cref{th:ordinary-amphoric} and the rest follows from \Cref{th:anab-galois-h1}.
\ep

\section{Anabelomorphic Density Theorems}\label{se:anab-den1}\nwss
Let us  illustrate  arithmetic  consequences of the anabelomorphic connectivity theorems (Theorem \ref{th:anab-connectivity-thm2}) by proving the following theorems.
\subsection{A Basic Density Theorem}
By \Cref{th:anab-proj-spaces}, one knows that projective spaces over $p$-adic fields or number fields are anabelomorphic varieties. Let us begin with the following elementary result which works out the case of certain open subsets of $\P^n$, but a similar statement can be formulated, using \Cref{th:anab-toric}, for suitable open subsets of a smooth, projective toric variety. 
\bthm[Anabelomorphic Density Theorem]\label{th:anabelomorphic-density} 
Let $n\geq 1$. Let $V_i\subset \P^n_{\Z}$ be the standard open subset defined by non-vanishing of the $i^{th}$-coordinate for $0\leq i\leq n$. Let $U$ be the intersection of some of the $V_0, V_1,\ldots,V_n$. 
Let $\anab{K}{K'}{\vseq{v}{n}}{\vseq{w}{n}}$ be anabelomorphically connected number fields. Then the  inclusion $$U(K')\subset \prod_i U(K'_{w_i})\isom \prod_i U(K_{v_i}) $$ is dense for the $p$-adic topology on the right (the fields $K_{w_i}'$ and $K_{v_i}$ may not be isomorphic).
\ethm
\bp 
From the definition of $U$ one sees that $U=\A^r\times \Gm^s$ for some integers $r,s\geq 0$. Hence from \Cref{th:third-fun-anab}, and \Cref{def:ef:anab-connected}, one sees that the two products in the statement are homeomorphic topological spaces. The density of $U(K')$ in the stated inclusion is immediate from  the fact that  weak Approximation Theorem holds for  $\A^r\times \Gm^s$ because it holds $\A^1$ \cite[Theorem 1.4]{platonov-book} and hence for the open subset $\Gm\subset \A^1$ \cite[Proposition 7.2(4)]{platonov-book} and for their products $\A^r\times\Gm^s$  \cite[Chap 7, Proposition 7.1(1)]{platonov-book}. Hence, the weak approximation holds for $U$ and the density assertion follows.
\ep

\subsection{Anabelomorphic Connectivity Theorem for Elliptic Curves}
In this subsection, we give a simple example illustrating how \Cref{th:anabelomorphic-density} (for $n=1$) can be used to transfer the data of an elliptic curve over  a number field to any  anabelomorphically connected number field while preserving some properties of the elliptic curve at the respective sets of primes of anabelomorphic connectivity. Let $$U=\P^1_{\Q}-\{0,1,\infty\}.$$ Then $U\subset \P^1$ is an open subset of the form considered in \Cref{th:anabelomorphic-density}. Fix an isomorphism of $\Q$-schemes $\P^1-\{0,1,\infty\}\isom \P^1-\{0,1728,\infty\}$. For any field $L$, one has $U(L)=L^*-\{1\}$. If  $L\anabelmap K$ is an anabelomorphism of $p$-adic fields, then one has an isomorphism $L^*\to K^*$ of topological groups and hence an isomorphism topological spaces
$$U(L)=L^*-\{1\}\isom K^*-\{1\}=U(K).$$ The composite mapping $U\to \P^1-\{0,1,\infty\}\isom \P^1-\{0,1728,\infty\}$ allows one to view the open subset $U(L)$ as $j$-invariants of elliptic curves over $L$ except for $j=0,1728$. 

\brem 
Note that these considerations can be applied to $U=(\P^1_\Z-\{0,1,\infty\})^m$ for any integer $m\geq 1$ and hence to moduli of hyperelliptic curves. The next theorem considers the genus one case.
\erem

The next result is motivated by constructions of \cite{taylor02} and the more general \Cref{qu:deform-galois} about deformations of mod-$\ell$ Galois representations, automorphic forms and anabelomorphically connected number fields.

\bthm\label{th:anabelomorphic-connectivty-theorem-elliptic} 
Let $$\anab{K}{K'}{\vseq{v}{n}}{\vseq{w}{n}}$$ be an anabelomorphically connected pair of number fields.  Let $E/K$ be an elliptic curve over a number field $K$  with $j$-invariant $j_E\neq 0,1728$. Assume that $E$ has potential multiplicative reduction at  $\{v_1,\ldots,v_n\}$ and potentially good reduction at all $v\not\in \{v_1,\ldots,v_n\}$.  Then there exists an elliptic curve $E'/K'$ with $j$-invariant $j_{E'}\neq0,1728$ such that
\benumlab
\item $E'/K'$  potentially good reduction at all $w\not\in\{w_1,\ldots,w_n\}$.
\item For $1\leq i\leq n$, one has $\ord_{v_i}(j_E)=\ord_{w_i}(j_{E'})$.
\item Hence $E'/K'$ has potential multiplicative reduction at all $\{w_i\}$.
\eenum
\ethm
\bp 
Let $j=j_E$ be the $j$-invariant of $E/K$. By \cite[Chap VII, Prop 5.5]{silverman-arithmetic} \cite[Chap. 10, Proposition 2.33]{qingliu-book}, at any place $v$ of bad potential multiplicative reduction one has $v(j)<0$. Let $\alpha_i:K_{v_i}\anabelmap K'_{w_i}$ be the given anabelomorphisms; let $j_i=\alpha_i(j)\in K_{w_i}^{'*}$. Then by \Cref{th:anabelomorphic-density} (applied to $U=\A^1_{\Z}-\{ 0, 1728\}$) one sees that $$\prod_i U(K_{v_i})=\prod_i \left(K_{v_i}^*-\{1\}\right) \isom   \prod_i \left(K_{w_i}^{'*}-\{1\}\right)=
\prod_i  U(K'_{w_i}) \supset U(K')$$ and the inclusion on the right is dense. This is not adequate to prove the theorem because weak approximation used in \Cref{th:anabelomorphic-density} is not adequate to control the behavior at primes $w\neq w_i$. However, as the coarse moduli of elliptic curves over $K'$ is the $j$-line $\A^1$ and hence by the Strong Approximation Theorem \cite[Theorem 1.5]{platonov-book}, there exists a $j'\in K'$ which is sufficiently close to $j_i$ for all $i=1,\ldots,n$ and $j'$ is integral at all other non-archimedean primes of $K'$. Since $\abs{j_i}_{w_i}>1$, and $j'$ is sufficiently close to the $j_i$, one sees that $j'\neq 0,1728$.   Hence,   there exists  a $j'\in K'-\{ 0, 1728\}$ which is sufficiently close to each of the $j_i$ and is $w$-integral for all other non-archimedean valuations $w$ of $K'$. 

By \cite[Chap. III, Proposition 1.4(c)]{silverman-arithmetic}, there exists an elliptic curve $E'/K'$ with $j$-invariant $j'$. By construction $j_{E'}=j'$, and as $j'$ is sufficiently close to $j_i$ for each $w_i$  and as $E/K$ has potential multiplicative reduction at each $v_i$, the valuation of $j'$ at each $w_i$ is negative. Hence, $E'/K'$ has potential multiplicative reduction over $K_{w_i'}'$. Moreover, for other non-archimedean valuations $w\neq w_1,\ldots,w_n$ of $K'$, $j'$ is $w$-integral by construction and so $E'$ has potential good reduction at such $w$. This proves all the assertions.
\ep

\subsection{Anabelomorphic version of Moret-Bailly's Theorem}\label{se:anab-den2}
Let $K$ be a $p$-adic field and let $\bK$ be an algebraic closure of $K$, let $G_K=\gal(\bK/K)$ be the absolute Galois  and let $X/K$ be a geometrically connected, smooth quasi-projective variety; let $\bar X=X\times_K\bK$. Let $\pi_1(X/K)$  be the \'etale fundamental groups of $X$ and let $\pi_1(\bar{X}/\bK)$ be the geometric \'etale fundamental group of $X/K$ computed using some geometric base-point $*:\Spec(\bK)\to X$. Then one has an exact sequence of topological groups
$$1\to  \pi_1(\bar{X}/\bK) \to \pi_1(X/K) \mapright{\eta} G_K\to 1.$$   

It is standard that any $K$-rational point of $X$ provides a section of $\eta$ (see \cite[Page xiv]{stix-book}). The \textit{Section Conjecture} of Alexander Grothendieck \cite[Cojecture 2, Page xiv]{stix-book} asserts that  $\pi_1(\bar{X}/\bK)$-conjugacy classes of continuous sections  $s:G_K\to \pi_1(X/K)$ of the surjection $\eta$ are in bijection with the set of rational points $X(K)$ (note that this statement is far broader than the one conjectured by Grothendieck). In this subsection, we write   $ {\rm Sect }(G_K,\pi_1(X/K))$ for the set of $\pi_1(\bar{X}/\bK)$-conjugacy classes of sections of $\eta$. Thus, the Section Conjecture asserts that one has a bijection (of sets) $X(K) \mapright{\isom} {\rm Sect }(G_K,\pi_1(X/K))$. As the set on the right is purely group theoretic, Grothendieck's Section Conjecture thus asserts that the set $X(K)$ is an amphoric set (\Cref{le:groth-conj}). Roughly speaking, \Cref{th:anab-affine-spaces}, \Cref{th:anab-proj-spaces} are examples of this sort of phenomenon.

\blem\label{le:groth-conj}
Suppose that  $X/K$ and $Y/L$ are  two geometrically connected, smooth, quasi-projective anabelomorphic varieties over $p$-adic fields $K,L$ (this hypothesis, together with \Cref{pr:anabelomorphy-of-base}, says that one has an anabelomorphism $K\anabelmap L$). Assume that Grothendieck's Section Conjecture holds for $X/K$ and $Y/L$. Then one has a natural bijection of sets
$$ X(K) \isom Y(L),$$
and in particular,  if $X(K)\neq \emptyset$ then $Y(L)\neq\emptyset$.
\elem

\bp[Proof of \Cref{le:groth-conj}] 
If $X/K$ and $Y/L$ are anabelomorphic  varieties i.e. $$\alpha:\pi_1(X/K)\mapright{\isom} \pi_1(Y/L),$$
then by \cite[Corollary 2.8(ii)]{mochizuki-topics1}, $\alpha$ preserves the corresponding geometric \'etale fundamental groups
$$\alpha(\pi_1(\bar{X}/\bK))=\pi_1(\bar{Y}/\bL),$$
and  hence by \Cref{pr:anabelomorphy-of-base}, the fields $K\anabelmap L$ are anabelomorphic i.e. $G_K\isom G_L$. This together with Grothendieck's section conjecture implies that there is a natural bijection of sets 
$$X(K)\isom{\rm Sect}(G_K,\pi_1(X/K)) \mapright{\isom} {\rm Sect}(G_L,\pi_1(Y/L)) \isom  Y(L).$$
The last assertion is obvious.
\ep

\newcommand{\vseqq}[2]{{#1}_{1,1},\ldots, {#1}_{1,r_1} ;\ldots ; {#1}_{{#2},1},\ldots, {#1}_{{#2},r_{#2} }}
Let us extend the notion of anabelomorphically connected number fields slightly. 
\bdefn
We  write $$\anab{K}{K'}{\vseq{v}{n}}{\vseqq{v'}{n}}$$ and say that $K,K'$ are \emph{anabelomorphically connected} along non-archimedean places $\vseq{v}{n}$ of $K$ and non-archimedean places $\vseqq{v'}{n}$ of $K'$ if one has anabelomorphisms
$$   K_{v_i}\anabelmap K'_{v'_{i,j}}  \text{ for each } i \text{ and for all } 1\leq j\leq r_i.$$
\edefn
Clearly this extends the notion introduced in \Cref{def:ef:anab-connected} by allowing several primes of $K'$ to correspond with each of the primes  $\vseq{v}{n}$ of $K$.

\brem 
A simple, but anabelomorphically trivial example of this definition is the following. Let $K$ be a number field and let $v$ be a non-archimedean prime of $K$. Suppose $K'/K$ is a finite extension such that $v$ splits completely in $K'$, say $w_1,\ldots,w_m$ are all the primes of $K'$ lying over $v$. Then one has an isomorphism of $p$-adic fields $K_v\isom K'_{w_i}$ for all $1\leq i\leq m$. Hence, one has (trivial) anabelomorphisms $K_v\anabelmap K'_{w_i}$ for $i=1,\ldots,m$. Hence, one sees that 
$$(K,\{v\}) \anabelmap (K',\{w_1,\ldots,w_m\})$$
are anabelomorphically connected number fields in the sense of the above definition. The main theorem of \cite{moret-bailly89}, and its application in \cite{taylor02} (and other works on potential automorphy) are via the formulation \cite[Theorem G]{taylor02}, are related to this example and motivates the next theorem.
\erem

\bthm\label{th:general-density-thm}
Let $K$ be a number field and let $v_1,\ldots,v_n$ be a finite set of non-archimedean places of $K$. Let $\anab{K}{K'}{\vseq{v}{n}}{\vseq{v'}{n}}$ be anabelomorphically connected number field. Let $X/K$ (resp. $Y/K'$) be a geometrically connected, smooth, quasi-projective variety over $K$ (resp. $K'$). Suppose the following conditions are met:
\benumlab 
\item  $X/K_{v_i}$ and $Y/K'_{v'_i}$ are anabelomorphic varieties for $1\leq i\leq n$, and
\item $X(K_{v_i})\neq\emptyset$ for all $1\leq i\leq n$, and
\item Grothendieck's section conjecture holds for each $X/K_{v_i}$ and $Y/K'_{v'_i}$, and 
\item suppose that one is given a non-empty open subset (in the $v_i$-adic topology) $ U_i\subseteq X(K_{v_i})$.
\eenum 
Then there exists a finite extension $K''/K'$ and places $\vseqq{v''}{n}$ of $K''$ such that
\benumlab
\item  one has the anabelomorphic connectivity chain $$\anab{K}{K'}{\vseq{v}{n}}{\vseq{v'}{n}}\anabelmap({K''},\{\vseqq{v''}{n}\})$$
\item for all corresponding primes in the above connectivity chain, bijections  $$ Y(K''_{v''_{i,j}})\isom  Y(K'_{v'_i}) \isom  X(K_{v_i})$$
\item a point $y\in Y(K'')$ whose image in $Y(K''_{v_{i,j}})\isom Y(K'_{v_i})\isom  X(K_{v_i})$ (for all $i,j$) is contained in $U_i$.
\eenum
\ethm

\brem
The situation considered in  \cite[Th\'eor\`eme 1.3]{moret-bailly89} (also see \cite[Corollary 1.5]{conrad2005}), the $p$-adic fields $K,L$ are isomorphic so one may take $Y=X$ (and hence $X/K$ and $X/L$ are trivially anabelomorphic) and the section conjecture  is unnecessary in \cite{moret-bailly89}. 
\erem

\bp[Proof of \Cref{th:general-density-thm}] 
The proof will use Lemma~\ref{le:groth-conj}. By the hypothesis that $X/K_{v_i},Y/K'_{v'_i}$ are anabelomorphic, one has by  Lemma~\ref{le:groth-conj}, that  for each $i$, there is a natural bijection of sets
$$X(K_{v_i})\isom Y(K'_{v'_i}),$$
and hence the latter sets are non-empty because of our hypothesis.

Now \cite[Th\'eor\`eme 1.3]{moret-bailly89} (or \cite[Corollary 1.5]{conrad2005}, or \cite[Theorem G]{taylor02}) can be applied to $Y/K'$ with $S=\{v_1',\ldots,v'_n\}$ so there exists a finite extension $K''/K'$ which is totally split at all the primes $v_i'$ into primes $v_{i,j}''$ with $j=1,\ldots, r_i=[K'':K']$ and hence for each $i$ one has isomorphisms $K'_{v'_i}\isom K''_{v_{i,j}}$ (for all $j$) and hence for each $i$ one has  $K'_{v'_i}\anabelmap K''_{v_{i,j}}$ (for all $j$).  Hence, one has the
anabelomorphic connectivity $(K',\{v_1',\ldots,v_n'\})\anabelmap({K''},\{\vseqq{v''}{n}\})$. The remaining conclusions are consequences of \cite[Th\'eor\`eme 1.3]{moret-bailly89} (\cite[Corollary 1.5]{conrad2005}).
\ep

\brem\label{cor:anab-moret-bailly-proj}
Let $K,K'$ be anabelomorphically connected number fields as in \Cref{th:general-density-thm}. Then  \Cref{th:general-density-thm} holds unconditionally (i.e. without assuming Grothendieck's Section Conjecture) for the following two cases: 
\benumlab
\item $X=\P^n_K$ and $Y=\P^n_{K'}$, or
\item  $X=\A^n_K$ and $Y=\A^n_{K'}$.
\eenum
This follows from \Cref{th:anab-proj-spaces} for {\bf(1)} and from \Cref{th:anab-affine-spaces} for {\bf(2)} and the proof of \Cref{th:general-density-thm}. 
\erem

\section{Weak Anabelomorphy}\label{se:weak-anabelomorphy}\nwss
\subsection{Definitions} As noted in \ssep\ref{se:intro}, one may think of anabelomorphy as an anabelian method of base change. In this section we want to elaborate on  this base change aspect. To this effect, let $F$ be a $p$-adic field, let $\bF$ be an algebraic closure of $F$. Let $X/F$ be a geometrically connected, smooth, quasi-projective variety over $F$. For any field extension $F'/F$ contained in $\bF$, write $X_{F'}=X\times_F{F'}$ for the base change of $X$ to $F'$. Consider the set 
$$[X,F]:=\{ X_{F'}: [F':F]<\infty \},$$
of all possible base changes of $X/F$ to finite extensions $F'/F$ (contained in $\bF$). We define an equivalence relation on the set $[X,F]$ as follows.

\bdefn\label{def:weakly-anab}
Let $X_K,X_L\in [X,F]$, then  one says that $X_K, X_L$ are \emph{weakly anabelomorphic} if  $K\anabelmap L$.
\edefn

The following is fundamental in understanding this:

\bpro\label{pr:weakly-anab}
Let $X/F$ be a geometrically connected, smooth, quasi-projective variety. Let $X_K,X_L\in[X,F]$.
\benumlab
\item Weak anabelomorphy is an equivalence relation $\sim$ on $[X,F]$. 
\item If $X_K$ and $X_L$ are anabelomorphic then they are also weakly anabelomorphic. 
\eenum
\epro
\bp 
The first assertion is immediate from the properties of anabelomorphic of $p$-adic fields. The second assertion follows from \Cref{pr:anabelomorphy-of-base}. 
\ep

\bdefn\label{def:weakly-amphoric}
Let $X/F$ be a geometrically connected, smooth, quasi-projective variety over a $p$-adic field $F$. Let $X_K\in [X,F]$. Then a quantity $Q_{X_K}$ or a property  $\sP$ associated to $X_K$ is said to be a \emph{weakly amphoric quantity} (resp.  \emph{weakly amphoric property}) if this quantity (resp. property) depends only on the weak anabelomorphism class of $X_K$ in $[X,F]$. More precisely: if,  $X_K\sim\ {X_L}$ for a pair $X_K, X_L\in [X,F]$, then one has $Q_{X_K}=Q_{X_L}$ (resp. the property $\sP$ holds for $X_K$ if and only if $\sP$ holds for $X_L$); one says that an algebraic structure $A_{X_K}$ is a \emph{weakly amphoric algebraic structure} if there is an isomorphism $A_{X_K}\isom A_{X_L}$ which is functorial in anabelomorphisms $K\anabelmap L$.  \edefn
\brem 
As will be seen in \Cref{th:kodaira-sym-unamphoric}{\bf(1,2)} and its proof, weakly amphoric quantities and properties do exist. However, at the moment, we do not know any nice examples of weakly amphoric algebraic structures. 
\erem

\subsection{Weak anabelomorphy and elliptic curves}\label{ss:weak-anab-elliptic}
\numberwithin{table}{subsection}
 As we have pointed out in \Cref{re:standin}, the upper-numbering ramification filtration is a stand-in for the field structure. From  \Cref{th:discriminant-is-unamphoric} one knows that the ramification filtration also impacts discriminants of anabelomorphic $p$-adic fields and from \Cref{th:artin-swan-unamphoric} one knows that Artin and Swan conductors are not amphoric in general. The ramification filtration enters discriminants and conductors of curves via presence of the wild ramification term i.e. the Swan conductor in \cite{serre-conductors} and the Grothendieck-Ogg-Shafarevich formula \cite[Page 450]{silverman-arithmetic}. Thus, one expects discriminants and conductors of elliptic and higher genus curves are not weakly amphoric. This led one to search for examples and led to \Cref{th:kodaira-sym-unamphoric} (and \Cref{th:discriminant-is-unamphoric2}) given below  (unexpectedly, some other well-known invariants of elliptic curves are not weakly amphoric).  As was shown in \cite{ogg1967}, for genus one, the wild ramification term in the conductor is zero for $p\neq 2,3$. Hence, in these and all other genus one examples presented in the tables,  $p=2$ or $p=3$ (in contrast to \Cref{th:discriminant-is-unamphoric}).   It is still possible that for genus one curves, the Kodaira Symbol i.e. the reduction type of the special fiber does jump around for $p\geq 5$. But one does not have any examples of this phenomenon  because simplest examples require working with $5$-adic fields with sufficiently deep ramification and this forces the calculation beyond the scope of present software.  See \Cref{th:discriminant-is-unamphoric2} for the higher genus case.
 
\bthm\label{th:kodaira-sym-unamphoric} 
Let $E/F$ be an elliptic curve over a $p$-adic field $F$. Let $E_K,E_L\in [E,F]$ be weakly anabelomorphic. Then
\benumlab
\item $E_K$ has potential good reduction if and only if $E_L$ has potential good reduction.
\item $E_K$ has potential multiplicative reduction if and only if $E_L$ has potential multiplicative reduction. 
\item In general, the following quantities are not weakly amphoric. 
\benum
\item The valuation of the discriminant of $E_K$,
\item the exponent of conductor of $E_K$.
\item The Kodaira Symbol of $E_K$, and
\item the Tamagawa number of $E_K$.
\eenum
\item In particular,  the number of irreducible components, counted without multiplicities, of the special fiber of $E_K$ is not weakly amphoric.
\eenum
\ethm

\brem 
In the split multiplicative reduction case, (limited) numerical evidence (see \Cref{table:long-numerical-data-table2}) suggests that the valuation of the discriminant,  exponent of the conductor, Kodaira Symbol, and Tamagawa number of $E_K$ are all weakly amphoric. But one does not know how to prove this at the moment.
\erem

\brem
The first two assertions of \Cref{th:kodaira-sym-unamphoric} are similar to \cite[Theorem~2.14(ii)]{mochizuki-topics1}.
\erem

\bp[Proof of \Cref{th:kodaira-sym-unamphoric}] 
Let $j_E$ be the $j$-invariant of $E$. As $K\supset F\subset L$, one has $j_E=j_{E_K}=j_{E_L}$, so write $j$ for this quantity. 
From \Cref{le:weil-deligne-groups}, one sees that $\ord_K(j)=\ord_L(j)$. By \cite[Chap VII, Prop 5.5]{silverman-arithmetic}, $E_K$ has potential good reduction if and only if $\ord_K(j)\geq 0$. If $j=0$ then $j$-invariant is integral in both $K$ and $L$ (because it is already so in $F$). So assume $j\neq 0$.  Then $\ord_F(j)\geq0$ if and only if $\ord_K(j)\geq0$ and $\ord_F(j)\geq0$ if and only if $\ord_L(j)\geq0$. This proves the first assertion.

Again from \cite[Chap VII, Prop. 5.5]{silverman-arithmetic} one sees that $E_K$ has potential multiplicative reduction if and only if $v_K(j)<0$ and as $v_K(j)<0$ if and only if $v_F(j)<0$ one similarly gets (2). 

So it remains to prove  assertions {\bf(3,4,5)}. The assertion {\bf(5)} is clear from the fact that there exists only finitely many $p$-adic fields (in any chosen algebraic closure of $F$) which are anabelomorphic to $K$.  To prove {\bf(3,4)} it suffices to give examples. This is done in the tables given below. All of these computations were carried out using SageMath \cite{sage}. The assertion {\bf(3)}(d)  is immediate from {\bf(3)}(c) because the Kodaira Symbol of $E_K$ corresponds to the dual graph of the special fiber and hence also encodes the number of irreducible components of the special fiber (of the minimal model of $E_K$) counted without multiplicities. Let  $m$ be the number (counted with multiplicities) of connected components of the special fiber (over the algebraic closure of the residue field). The dual graph of the special fiber and the number $m$ can be read off from \cite[Table 4.1, Page 365]{silverman-advanced}.

Let $F=\Q_3(\zeta_9)$, let $K=F(\sqrt[9]{3})$ and $L=F(\sqrt[9]{2})$. Then $K\anabelmap L$ from \Cref{le:basic-example}. Both of these field have degree
$$[K:\Q_3]=[L:\Q_3]=54.$$

Let $E: y^2=x^3+3x^2+9$ and $E_K$ and $E_L$ be as above.  Let $\Delta$ be the minimal discriminant (over the relevant field), $f$ be the exponent of the conductor, the list of Kodaira Symbols and the definition of the Tamagawa number are in \cite[Chap. IV]{silverman-advanced}. The following table shows the values for $E_K$ and $E_L$.

\begin{center}
	\begin{tabular}{|c|c|c|c|c|c|}
		\hline \rule[-2ex]{0pt}{5.5ex} Curve  & $v(\Delta)$ & $f$ & Kodaira Symbol &  Tamagawa Number & $m$ \\ 
		\hline \rule[-2ex]{0pt}{5.5ex} $E_K$ & $6$ & $4$ &  $IV$ &  $1$ & $3$\\ 
		\hline \rule[-2ex]{0pt}{5.5ex} $E_L$ & $6$ & 2 & $I_0^*$ &  $4$ & $5$\\ 
		\hline 
	\end{tabular} 
\end{center}

Here is another example let $E: y^2=x^3+3x^2+3$ and let $K,L, E_K,E_L$ be as above. Then one has
\begin{center}
	\begin{tabular}{|c|c|c|c|c|c|}
		\hline \rule[-2ex]{0pt}{5.5ex} Curve  & $v(\Delta)$ & $f$ & Kodaira Symbol &  Tamagawa Number & $m$ \\ 
		\hline \rule[-2ex]{0pt}{5.5ex} $E_K$ & $12$ & $6$ &  $IV^*$ &  $3$ & $7$\\ 
		\hline \rule[-2ex]{0pt}{5.5ex} $E_L$ & $12$ & 10 & $IV$ &  $1$ & $3$\\ 
		\hline 
	\end{tabular} 
\end{center}

Here is an example with $p=2$. let $E: y^2=x^3-129784x+17996160$ (LMFDB Label 5888.d1). Let  $K=\Q_2(\zeta_{16},\sqrt[4]{1-\zeta_{16}})$, $L=\Q_2(\zeta_{16},\sqrt[4]{5})$ and $E_K,E_L$ be as above. Then one has
\begin{center}
	\begin{tabular}{|c|c|c|c|c|c|}
		\hline \rule[-2ex]{0pt}{5.5ex} Curve  & $v(\Delta)$ & $f$ & Kodaira Symbol &  Tamagawa Number & $m$ \\ 
		\hline \rule[-2ex]{0pt}{5.5ex} $E_K$ & $12$ & $12$ &  $II$ &  $1$ & $1$\\ 
		\hline \rule[-2ex]{0pt}{5.5ex} $E_L$ & $24$ & 20 & $I_0^*$ &  $4$ & $5$\\ 
		\hline 
	\end{tabular} 
\end{center}
This completes the proof of all the assertions.
\ep

\subsection{Additional numerical examples}\label{ss:add-num-data} Here are two more  examples where all the quantities are simultaneously different.

\newcommand{\zetap}{\zeta_9}
Let $$F=\Q_3(\zeta_9)\quad K=F(\sqrt[9]{3}) \quad L=F(\sqrt[9]{4}),$$
and let
$$E:y^2=x^3+(-\zetap^5 + 8\zetap^4 - \zetap^3 + \zetap^2 - 2\zetap - 11) x+ (-408\zetap^5 - 6\zetap^4 + 201\zetap^3 + 37\zetap^2 - 38\zetap + 1348).$$

\begin{center}
	\begin{tabular}{|c|c|c|c|c|c|}
		\hline \rule[-2ex]{0pt}{5.5ex} Curve  & $v(\Delta)$ & $f$ & Kodaira Symbol &  Tamagawa Number & $m$  \\ 
		\hline \rule[-2ex]{0pt}{5.5ex} $E_K$ & $15$ & $15$ &  $II$ &  $1$ & $1$\\ 
		\hline \rule[-2ex]{0pt}{5.5ex} $E_L$ & $39$ & $37$ & $IV$ &  $3$ & $3$\\ 
		\hline 
	\end{tabular} 
\end{center}

For the same fields $F,K,L$ as in the previous example and for the curve
$$E:y^2=x^3+ (-2\zetap^5 + \zetap^4 + \zetap^3 - \zetap^2 + 2\zetap + 5) x+ ( 869\zetap^5 + 159\zetap^4 - 47\zetap^3 - 125\zetap^2 + 354\zetap + 713).$$
\begin{center}
	\begin{tabular}{|c|c|c|c|c|c|}
		\hline \rule[-2ex]{0pt}{5.5ex} Curve  & $v(\Delta)$ & $f$ & Kodaira Symbol &  Tamagawa Number & $m$  \\ 
		\hline \rule[-2ex]{0pt}{5.5ex} $E_K$ & $15$ & $9$ &  $IV^*$ &  $3$ & $7$\\ 
		\hline \rule[-2ex]{0pt}{5.5ex} $E_L$ & $27$ & $19$ & $II^*$ &  $1$ & $9$ \\ 
		\hline 
	\end{tabular} 
\end{center}

Now let us provide two examples for $p=2$. These examples are taken from our data because at least four of  the five quantities are simultaneously different. Let $F=\Q_2(\zeta_{16})$, $K=F(\sqrt{\zeta_8-1},\sqrt{\zeta_8^3-1})$, $L=F(\sqrt[4]{\zeta_4-1})$. The two fields $K,L$ were shown to be anabelomorphic in \cite{jarden79} and are totally ramified extensions of $\Q_2$ of degree $n=32$.
\newcommand{\zetas}{\zeta_{16}}
\begin{multline*}
E:y^2=x^3+(-2\zetas^7 + 2 \zetas^6 - 2 \zetas^5 + 2 \zetas^4 - 2 \zetas^3 + 4 \zetas^2 + 6 \zetas + 30)x \\
\qquad \qquad \qquad + (32 \zetas^7 - 76 \zetas^6 - 8 \zetas^5 + 32 \zetas^4 - 24 \zetas^3 - 20 \zetas^2 + 16 \zetas - 28).
\end{multline*}

Then
\begin{center}
	\begin{tabular}{|c|c|c|c|c|c|}
		\hline \rule[-2ex]{0pt}{5.5ex} Curve  & $v(\Delta)$ & $f$ & Kodaira Symbol &  Tamagawa Number & $m$ \\ 
		\hline \rule[-2ex]{0pt}{5.5ex} $E_K$ & $64$ & $60$ &  $I_0^*$ &  $2$ & $5$ \\ 
		\hline \rule[-2ex]{0pt}{5.5ex} $E_L$ & $52$ & $52$ & $II$ &  $1$ & $1$\\ 
		\hline 
	\end{tabular} 
\end{center}

\begin{multline*}
E:y^2=x^3+(-2 \zetas^6 - 2 \zetas^4 + 4 \zetas^2 + 2)x \\
\qquad \qquad \qquad + (28 \zetas^6 - 40 \zetas^5 - 24 \zetas^4 + 8 \zetas^3 + 16 \zetas^2 - 40 \zetas + 60).
\end{multline*}

Then
\begin{center}
	\begin{tabular}{|c|c|c|c|c|c|}
		\hline \rule[-2ex]{0pt}{5.5ex} Curve  & $v(\Delta)$ & $f$ & Kodaira Symbol &  Tamagawa Number & $m$  \\ 
		\hline \rule[-2ex]{0pt}{5.5ex} $E_K$ & $68$ & $60$ &  $II^*$ &  $1$ & $9$\\ 
		\hline \rule[-2ex]{0pt}{5.5ex} $E_L$ & $56$ & $52$ & $I_0^*$ &  $2$ & $5$ \\ 
		\hline 
	\end{tabular} 
\end{center}
\brem
Numerical data of  Table~\ref{table:long-numerical-data-table2} suggests that if $E$ has semistable reduction, then the four quantities considered above are weakly amphoric These examples reveal that Tate's algorithm \cite[Chapter IV, 9.4]{silverman-advanced} for determining the special fiber of an elliptic curve over a $p$-adic field is dependent on the intertwining between the additive and multiplicative structure of the $p$-adic fields. \erem

\subsection{Weak anabelomorphy of Artin Conductors, Swan Conductors and Discriminants of curves}\label{se:swan-II}
The results of this section complement the results for genus $1$ of \ssep\ref{ss:weak-anab-elliptic} and  \Cref{th:artin-swan-unamphoric} on Swan Conductors.  Let $F$ be a $p$-adic field, let $\bF$ be an algebraic closure of $F$, let $X/F$ be a geometrically connected, smooth quasi-projective curve over $F$. Let $X_K, X_L\in [X,F]$ with $K\anabelmap L$  anabelomorphic $p$-adic fields containing $F$. Write $X_K=X\times_F K$ and $X_L=X\times_F L$.

\newcommand{\swan}{\mathop{\rm Swan}}

For geometric applications discussed in this section it will be convenient to work with  strictly Henselian rings. As Artin and Swan conductors are unaffected by passage to unramified extensions, this passage to strictly Henselian rings is harmless. In particular, one can work over $\knr$.

Let $X/K$ be a geometrically connected, smooth, proper curve and $X_{\bar{\eta}}$ (resp. $X_s$) is the geometric generic fiber (resp. special fiber) of a regular, proper model  then one has a \emph{discriminant} $\Delta_{X/K}$ defined as in \cite{saito1988}. This discriminant coincides with the usual discriminant if $X/K$ is an elliptic curve.

\bthm\label{th:discriminant-is-unamphoric2} Let $F$ be a $p$-adic field. Then the Swan conductor and the discriminant of a geometrically connected, smooth, projective curve of genus $g\geq2$ over $F$ is not weakly amphoric in general. 
\ethm
\bp 
This will be proved using the conductor formula of \cite{saito1988}. That paper works with strictly Henselian discrete valuation rings. \Cref{le:anab-knr} says that if $K\anabelmap L$ are anabelomorphic $p$-adic fields, then so are $\knr\anabelmap \lnr$. Thus, one can work with the valuation ring of $\knr$ (resp. $\lnr$) which is strictly Henselian and preserve the property that the new base fields $\knr,\lnr$ are also anabelomorphic. Write $S_K=\Spec(\O_{ \knr})$ and similarly define $S_L$. Given a geometrically connected, smooth,  projective curve $X/K$, one can choose a regular, proper and flat model $\mathfrak X\to S_K$ with generic fiber isomorphic to $X_{\knr}=X\times_K \knr$, and geometric generic fiber $\mathfrak{X}_{\bK}$ isomorphic to $X_{\bK}=X\times_K\bK$ and similarly choose $\mathfrak Y\to S_L$. Let $\mathfrak{X}_s$ denote the special fiber of $\mathfrak{X} \to S_K$. Then Swan conductor of $X/K$ is defined by
$$\swan(X_K):=\swan(\mathfrak{X}/S_K)=\sum_{i\geq 0}(-1)^i\swan(H^i_{\text{\'et}}(X_{\bK},\Q_\ell)).$$ Since  $X/K$ is a curve, then  the $H^0$ and $H^2$ terms in the  formula for $\swan(X_K)$ are zero as the action of $\gal(\bK/\knr)$ on these two terms is trivial. Hence, the alternating sum reduces to the $H^1$ term. 

The relationship between the discriminant, Artin conductor, the Swan conductor is established by Saito's conductor formula (\cite[Theorem 1]{saito1988}) which, under the hypothesis on $X/K$ asserts,  that one has
$$
-\ord_K(\Delta_{X/K})={\rm Artin}(X/K)=\chi_{\text{\'et}}(\mathfrak X_{\bK})-\chi_{\text{\'et}}(\mathfrak{X}_s)-\swan(H^1_{\text{\'et}}(X_{\bK},\Q_\ell)),$$
where $\Delta_{X/K}$ is the discriminant defined in \cite[Page 155]{saito1988}.  If  $X/K$ is of genus $1$, then, all the quantities entering Saito's formula above agree with the familiar quantities: discriminant, Swan and Artin conductor which appear in \Cref{th:kodaira-sym-unamphoric}.

Let $K\anabelmap L$ be an anabelomorphism of $p$-adic fields containing a $p$-adic field $F$ and contained in some fixed algebraic closure of $F$. By \Cref{le:anab-knr}, one has the anabelomorphism $\knr\anabelmap \lnr$.  
To prove the theorem one wants to compare $\swan(X_{K})$ and $\swan(X_{L})$. 

In \Cref{th:artin-swan-unamphoric} and \Cref{co:break-decomp-unamphoric}, we have shown that Artin and Swan conductors of Galois representations, and the breaks in the break-decomposition  are not amphoric in general  because these invariants depend on the ramification filtration. Hence, the assertion of the theorem is immediate. 
\ep

\brem We provide some comments on finding examples of genus two curves exhibiting the behavior predicted by \Cref{th:discriminant-is-unamphoric2}. Like the genus one case (\Cref{th:kodaira-sym-unamphoric}),  genus two examples can be found using the algorithm for determination of the special fiber of a genus two curve \cite{qingliu1994} (also see \cite{ogg1966}, \cite{ueno1973}). The phenomenon of $\swan(X_K)\neq 0$ is related to whether or not $X/F$ acquires stable reduction over wildly ramified extensions. For genus $g=1$ this happens only for $p=2,3$. In genus $2$ case, if $p\geq 7$, then $\swan(X_K)=0$ (\cite{ogg1966}, \cite{ueno1973}) and in general, if $X/F$ is hyperelliptic and if $p>2g+1$, then $\swan(X_K)=0$ \cite[Remark 5.7]{dokchitser2022}. For $g=2$ and $p=2,3,5$ examples may be found by reading off the conductor from the algorithm for determining the special fiber. Explicitly \cite[Example 2, Page 78]{qingliu1994} shows that for $p=5$ and  $y^2=x^5+t$ and a suitable $t$, one has a non-zero contribution for the Swan conductor. Since Swan conductor is sensitive to the inertia filtration, one sees that there exist some extensions $K\anabelmap L$ such that two Swan conductors, $\swan(X_K)$ and $\swan(X_L)$ are distinct. More generally, for any odd prime $p$, for a non-singular hyper-elliptic curve $y^2=x^p+t$ over $\Q_p$ of genus $g=\frac{p-1}{2}$ one has wild action of inertia. Other examples of this type also exist  \cite[Theorem 12.3, Remark 12.4, Example 12.6]{dokchitser2022} and also \cite{azon2024}. Numerical examples can be found using the methods used for \Cref{th:kodaira-sym-unamphoric} and \href{https://people.maths.bris.ac.uk/~matyd/redlib/}{Tim Dokchitser's Magma Scripts}.
\erem

\section{Anabelomorphy and perfectoid fields and spaces}\label{se:perfectoids}\nwss
\newcommand{\kf}{K^\flat}
\newcommand{\lf}{L^\flat}
The purpose of this section is to illustrate the unexpected parallels and relationships between perfectoid geometry \cite{scholze12-perfectoid-ihes} and anabelomorphy studied in the present paper. This relationship came to light in the course of writing of this paper and \cite{joshi-formal-groups}. 

\subsection{Anabelomorphy of perfectoid fields} Let $K$ be a perfectoid field of characteristic zero (see \cite[Section 3]{scholze12-perfectoid-ihes}). Let $\kf$ be its tilt (see \cite[Lemma 3.3]{scholze12-perfectoid-ihes}). The following is a formulation of \cite[Theorem 3.7]{scholze12-perfectoid-ihes} from the point of view of anabelomorphy:
\bthm\label{th:perfect-anab}
Let $K,L$ be perfectoid fields with an isometry $K^\flat\isom L^\flat$ between their respective tilts. Then 
\benumlab 
\item one has anabelomorphisms of perfectoid fields $$K\anabelmap K^\flat \anabelmap L^\flat \anabelmap L,$$
\item and one has topological isomorphism of the multiplicative monoids (i.e. an isomorphism of their multiplicative structures): $$\invlim_{x\mapsto x^p}K\isom K^\flat \isom L^\flat\isom  \invlim_{x\mapsto x^p}L,$$
\item 
In particular, if $F$ is a perfectoid field of characteristic $p>0$ and $(K,K^\flat\isom F)$ and $(L,L^\flat\isom F)$ are two untilts of $F$, then one has an anabelomorphism $K\anabelmap L$.\eenum
\ethm
\bp 
Let $G_K$ (resp. $G_{K^\flat}$) be the absolute Galois group of $K$ (resp. $K^\flat$). Then by  \cite[Theorem 3.7]{scholze12-perfectoid-ihes} (also see \cite{fontaine-wintenberger1979b}) one has an isomorphisms 
$$G_K\isom G_{K^\flat}\textrm{ and }G_L\isom G_{L^\flat}.$$ 
Since $K^\flat \isom L^\flat$, one also has an isomorphism $G_{K^\flat}\isom G_{L^\flat}$. Putting both of these together one obtains {\bf(1)}.
The assertion {\bf(2)} is immediate from \cite[Lemma 3.4(iii)]{scholze12-perfectoid-ihes} and {\bf(3)} is immediate from {\bf(1)}.
This proves the assertion.
\ep

\begin{example}\label{ex:perfectoid-fields}
Here is an explicit example of \Cref{th:perfect-anab}. The $p$-adic completions $$K=\widehat{\Q_p(\zeta_p,\zeta_{p^2},\cdots)}\subset \C_p \textrm{ and } L=\widehat{\Q_p(\sqrt[p]{p},\sqrt[p^2]{p},\cdots)}\subset \C_p$$ are both perfectoid fields.  From \cite[Example 2.1.1]{weinstein-aws}, one has an isometry  $\kf\isom \F_p((t^{1/p^\infty}))\isom \lf$ of their tilts, and hence one has an anabelomorphism $K\anabelmap L$ of perfectoid fields.
\end{example}

\subsection{Anabelomorphy of perfectoid spaces}
\Cref{th:perfect-anab} has a higher dimensional analog.
Suppose that $K$ is a perfectoid field. Let $X/K$ be a connected perfectoid space over $K$ \cite[Definition 6.15]{scholze12-perfectoid-ihes}. Let $X^\flat/\kf$ be its tilt (see \cite[Definition 6.16]{scholze12-perfectoid-ihes}). Let $\pi_1(X/K)$ be its \'etale fundamental group for a suitable choice of geometric base point. Then one has the following:
\bthm\label{th:perfect-anab3} 
Let $K,L$ be perfectoid fields with isometric tilts. Let $X/K, Y/L$ be connected perfectoid spaces with an isomorphism of the tilts $$X^\flat/K^\flat\isom Y^\flat/L^\flat.$$
Then one has  anabelomorphisms of perfectoid spaces
$$X/K\anabelmap X^\flat/K^\flat \anabelmap Y^\flat/L^\flat \anabelmap Y/L.$$
\ethm
\bp 
This is a consequence of the stronger assertion  \cite[Theorem 7.12]{scholze12-perfectoid-ihes} which implies that the categories of finite \'etale covers of $X/K$ and $X^\flat/K^\flat$ are naturally equivalent.
\ep

The following theorem is the perfectoid analog of \Cref{th:anab-proj-spaces}.
\bthm\label{th:anab-proj-spaces-perf} 
Let $K,L$ be perfectoid fields with isometric tilts $K^\flat\isom L^\flat$. Let $(\P_K^n)^{\mathit{perf}}$ be the perfectoid projective space. Then one has an anabelomorphism
$$ (\P_K^n)^{\mathit{perf}} \anabelmap (\P_L^n)^{\mathit{perf}}$$
and an isomorphism of topological spaces
$$\abs{(\P_K^n)^{\mathit{perf}}} \isom  \abs{(\P_L^n)^{\mathit{perf}}}.$$
\ethm
\bp 
This is a consequence of \cite[Theorem 7.1]{scholze14-icm2}.
\ep

Other examples of this phenomenon arise in the theory of Diamonds \cite{scholze-diamonds}:
\bthm
Let $F$ be a perfectoid field of characteristic $p>0$. Let $K\anabelmap L$ be anabelomorphic $p$-adic fields (i.e. $G_K\isom G_L$), then the diamonds  $\mathscr{X}_{F,K}^\Diamond$ and $\mathscr{X}_{F,L}^\Diamond$ (\cite[Lecture 15]{scholze-weinstein-book}) associated with complete Fargues-Fontaine curves $\mathscr{X}_{F,K}$ and $\mathscr{X}_{F,L}$ over $K$ and $L$ respectively are anabelomorphic:
$$\pi_1(\mathscr{X}_{F,K}^\Diamond)\isom G_K\isom G_L\isom \pi_1(\mathscr{X}_{F,L}^\Diamond).$$
\ethm
\bp 
This is immediate from the proof of \cite[Theorem 16.3.1]{scholze-weinstein-book}. 
\ep
\brem 
In contrast, in \cite{joshi-gconj}, it is shown that if $K\anabelmap L$ is a strict anabelomorphism of $p$-adic fields, then one has a strict anabelomorphism of schemes $$\mathscr{X}_{F,K}\anabelmap \mathscr{X}_{F,L}$$ i.e. these schemes are not $\Z$-isomorphic. Also see \cite[\ssep 8.7]{joshi-teich-rosetta}. \erem

\section{Anabelomorphy for $p$-adic differential equations}\label{se:anab-diff-eq}\nwss
\newcommand{\xan}{X^{{\rm an}}}
\newcommand{\yan}{Y^{{\rm an}}}
This section is motivated by the results in the archimedean case detailed in \ssep\ref{se:anab-arch}.
\subsection{Anabelomorphy of some $p$-adic differential equations}
Let $X/K$ be a geometrically connected, smooth, quasi-projective variety over a $p$-adic field $K$. Let $\xan/K$ denote the strictly analytic Berkovich space associated by \cite{berkovich-book} to $X/K$.   In this section $\pi_1(X/K)$ will stand for the \'etale fundamental group  of  $\xan/K$ (computed using a choice of a geometric base-point) defined in \cite[Chapter III, 1.4.1]{andre03}. A reference for $p$-adic differential equations on $\xan/K$ considered below is \cite[Chap. III, Section 3]{andre-book}. The purpose of this section is to prove the following (for the archimedean version see \Cref{th:anab-arch-diff-eq}). 
\bthm\label{th:anab-diff-eq1}
Let $X/K$ and $Y/L$ be two geometrically connected, smooth, quasi-projective varieties over $p$-adic fields $K$ and $L$. Assume that $\xan/K$ and $\yan/L$ are anabelomorphic strictly analytic spaces with an anabelomorphism $\alpha: \pi_1(\yan/L)\isom \pi_1(\xan/K)$ and suppose that this induces an anabelomorphism $L\anabmapright{\alpha} K$ of $p$-adic fields. Then $\alpha$ induces a natural bijection $\alpha$ between rank one, \'etale $p$-adic differential equations on $\yan_{et}/L$ and $\xan_{et}/K$ which associates to a rank one, \'etale $p$-adic differential equation $(M_{et},\nabla)$ on  $\yan_{et}/L$, rank one  $p$-adic differential equation $(N_{et},\nabla)$ on $\xan_{et}/K$ such that the associated discrete $K$-representation of $\pi_1(\xan/K)$ (provided by the $p$-adic Riemann-Hilbert correspondence \cite[Chapter III, Theorem 3.4.6]{andre-book}) is given composing with $\alpha^{-1}:\pi_1(\yan/L)\mapright{\isom} \pi_1(\xan/K)$. 
\ethm
\bp 
By \cite[Corollary 2.8(ii)]{mochizuki-topics1}, the tempered anabelomorphism $\alpha$ induces an anabelomorphism $L\anabmapright{\alpha} K$ of $p$-adic fields. The  Riemann-Hilbert Correspondence \cite[Chapter III, Theorem 3.4.6]{andre-book} for rank one differential equations establishes an equivalence between the category of discrete one dimensional representations $$\pi_1(\yan/L) \mapright{\rho}\glo L=L^*$$ and the category pairs $(M_{et},\nabla)$ consisting of a locally free rank one $\O_{\yan_{et}/L}$-module $M_{et}$ and an $L$-linear connection $\nabla$ on $M_{et}$. 

Now the composition 
$$\rho:\pi_1(\xan/K)\mapright{\alpha^{-1}}\pi_1(\yan/L) \mapright{\rho} \glo L=L^*\mapright{\alpha} K^*$$
gives a discrete representation 
$$\rho':\pi_1(\xan/K)\mapright{\alpha^{-1}\circ\rho\circ \alpha } K^*=\glo K$$
and hence by \cite[Chapter III, Theorem 3.4.6]{andre-book}, a rank one differential equation $(N_{et},\nabla)$ on $\xan_{et}/K$ of rank one.
\ep

Based on \Cref{th:automorphic-ordinary-synchronization}, other results of \ssep\ref{se:anab-langlands} and \Cref{th:anab-arch-diff-eq}, one expects that 
\bcon\label{conj:p-adic-diff-eq}
\Cref{th:anab-diff-eq1} holds true for differential equations of all ranks $n\geq 1$.
\econ
\subsection{Weak anabelomorphy and $p$-adic differential equations}
Let $X/F$ be a geometrically connected, smooth, quasi-projective variety over a $p$-adic field $F$. Let $\bF$ be an algebraic closure of $K$, and let $K\anabelmap L$ be anabelomorphic $p$-adic fields containing $F$ and contained in $\bF$. Then given any anabelomorphism $\alpha:K \anabelmap L$ one can consider the given $p$-adic differential equation as giving $p$-adic differential equations on $\xan/K$ and $\xan/L$ respectively.  In particular, it is possible to ask if there are quantities, properties algebraic structures associated to a differential equation on $X/K$ which are weakly amphoric or not weakly amphoric (with respect to  anabelomorphisms $K\anabelmap L$). 

An important invariant of a $p$-adic differential equation is the index of irregularity \cite[Chapter III, 3.1.2]{andre-book} at a singular point. It is well-known that  the Swan conductor of a Galois representation  is  the analog, in theory of differential equations, of the local index of irregularity. Hence,  \Cref{th:artin-swan-unamphoric} suggests that the following conjecture is natural.
\bcon[Index of Irregularity is not weakly amphoric]\label{conj:irreg-unamphoric}
In the above notation, assume that $X/F$ is a curve (i.e. $\dim(X)=1$). Then  the index of irregularity of a $p$-adic differential equation $(M,\nabla)$ on $X/F$ is not weakly amphoric in general. More generally, the irregularity module of the differential equation $(M,\nabla)$ over $X/F$ is not weakly amphoric ($X$ need not be a curve for this).
\econ

\section{Anabelomorphy at Archimedean primes}\label{se:anab-arch}\nwss
As is well-known, many non-isomorphic quasi-projective complex varieties have isomorphic topological and hence also \'etale fundamental groups. So the naive approach to the study of  anabelomorphy via isomorphisms of fundamental groups is not useful in this archimedean setting. During the writing of \cite{joshi-teich}, one came to recognize that an alternate approach to the archimedean theory of anabelomorphisms, via the work of \cite{nakai59, nakai60,nakai-book,nakai1972} exists, and it naturally includes classical Teichm\"uller Theory of Riemann surfaces. Since \cite{joshi-teich} is not yet published, here we provide a self-contained treatment in the archimedean setting. The theory given here is broader than the one used in classic works in algebraic geometry such as \cite{deligne-zeta-values}, \cite{hain88} where isomorphisms between fundamental groups arise from isomorphisms between algebraic varieties. 
\subsection{Nakai quasi-isometries}

For the basics of Riemannian Geometry, the reader should consult standard texts on the subject (for example, \cite{JohnLee-riemannian-book} or \cite{Petersen-riemannian-book}). Let $(M,g)$ be a Riemannian manifold with a Riemannian metric $g$. Let $x,y\in M$ be two points and $\gamma$ be a rectifiable path joining $x,y$ in $M$. Then the line integral $\int_\gamma ds$ of $\gamma$ is well-defined using the metric $g$ of $M$. Let $$d_{(M,g)}(x,y)=\text{inf}_{\gamma} \{ \int_\gamma ds\}$$
where the infimum is over all the rectifiable paths joining $x,y$ in $M$.	We call $d_M(x,y)$ the \emph{natural distance between $x,y\in M$}. The function $M\times M\to \R$ given by $(x,y)\mapsto d_M(x,y)$ is a metric on $M$  which induces the same topology as the given topology of $M$ (see \cite[Chapter 5]{Petersen-riemannian-book} or \cite[Chapter 6]{JohnLee-riemannian-book}).

\bdefn
 Let $(M_1,g_1), (M_2,g_2)$ be two Riemannian manifolds. A \emph{Nakai quasi-isometry} between $f:(M_1,g_1)\to (M_2,g_2)$ is a homeomorphism $f:M_1\to M_2$ satisfying
$$K^{-1}d_{(M_1,g_1)}(x,y)\leq d_{(M_2,g_2)}(f(x),f(y))\leq  K d_{(M_1,g_1)}(x,y)$$
for some constant $K\geq 1$ and all $x,y\in M_1$.
\edefn
\brem\  
\benumlab
\item As is noted in \cite{nakai1972}, a Nakai quasi-isometry is a quasi-conformal mapping of Riemannian manifolds (for quasi-conformal mappings between Riemann surfaces see \cite{ahlfors-quasiconf-book} or \cite[Chapter 4]{imayoshi-book}).
\item If $f:(M_1,g_1)\to (M_2,g_2)$ is Nakai quasi-isometry, then so is $f^{-1}:(M_2,g_2)\to (M_1,g_1)$ \cite{nakai1972}. 
\item It is easily checked from the definition that the composition of two Nakai quasi-isometries is also a Nakai quasi-isometry. 
\item Hence, `Nakai quasi-isometry' gives an equivalence relation on Riemannian manifolds. 
\item Readers should beware that there are many, possibly inequivalent, definitions of quasi-isometries, here we will work with the specific one made in \cite{nakai1972}, \cite[Appendix]{nakai-book} (also see \cite{nakai59}, \cite{nakai60}).
\item The Riemannian manifolds of interest to here are of the following type. Let $X/\C$ be a connected, smooth quasi-projective variety over complex numbers. Then $X\into \P^n$ for some $n$ by the virtue of its quasi-projectivity and the restriction of Fubini-Study metric on $\P^n$ (see \cite[Example 3.1.9]{huybrechts-complex-geometry-book}) to $X\into\P^n$ gives the K\"ahler structure of $\xan$ (see \cite[Chapter 3]{huybrechts-complex-geometry-book}). Hence, $\xan$ is a K\"ahler manifold (and hence a Riemannian manifold). We will consider complex, smooth quasi-projective varieties as Riemannian manifolds via this well-known fact.  
\eenum
\erem

\subsection{Anabelomorphisms of complex varieties}
From now on, we work with  connected, smooth quasi-projective varieties over $\C$. If $X$ is such a variety, let $\xan_\C$ (or $\xan$)  be the associated complex manifold, considered as being equipped with a Riemannian metric (say the Fubini-Study metric given using some embedding of $X\into\P^n_\C$ for some $n\geq 1$). Write $\pi_1(X)$ for the topological  fundamental group, $\pi_1(\xan)$, of $\xan$ computed using some choice of a base-point, and $\pi_1^{et}(X)$ for the \'etale fundamental group of $X$ computed using some (geometric) base-point.

\bdefn\label{def:anab-arch}
Let $X,Y$ be  connected, smooth, quasi-projective varieties over $\C$. Consider the complex manifolds $\xan, \yan$ as equipped with  Riemannian metrics. We say that \textit{$X$, $Y$ are anabelomorphic complex quasi-projective varieties} if there exists a Nakai quasi-isometry $f:\xan \to \yan$ of Riemannian manifolds. We write this as $X\anabmapright{f} Y$. \edefn

\brem\ 
\benumlab 
\item Note that an anabelomorphism $X\anabmapright{f} Y$ defined above need not be a morphism of complex varieties and it may need not even be a holomorphic mapping. In particular $f$ need not be algebraic.
\item `Nakai quasi-isometry' defines an equivalence relation on connected, quasi-projective, smooth complex varieties. 
\item This allows us to extend the notion of amphoric quantities, properties and algebraic structures given in \Cref{def:amphoric-var} to anabelomorphisms of quasi-projective complex varieties. For example a property of a smooth quasi-projective complex variety $X$ is  amphoric if every connected smooth quasi-projective complex variety $Y$ which is anabelomorphic to $X$ has this property etc.
\eenum
\erem

The following lemma explains the significance of this definition for anabelomorphy:
\blem\label{le:anab-arch-prop}
Suppose $X\anabmapright{f} Y$ is an anabelomorphism of complex varieties in the sense of \Cref{def:anab-arch}, then one has the following:
\benumlab
\item $f$ is a homeomorphism of the analytic spaces $\xan\mapright{f}\yan$, and 
\item if $\xan,\yan$  are Riemann surfaces then the Nakai quasi-isometry  $f$ is simply a quasi-conformal mapping $f:\xan \to \yan$ of Riemann surfaces (for a definition see \cite{ahlfors-quasiconf-book}).
\item In all dimensions, following hold:
\begin{enumerate}
\item $\pi_1(\xan) \isom \pi_1(\yan)$; and hence
\item $\pi_1^{et}(X)\isom \pi_1^{et}(Y)$;
\item $\dim(X)=\dim(Y)$;
\item  if $\dim(X/\C)=1$, the topological type $(g,n)$ of $\xan$ is the same for $\yan$.
\end{enumerate}
\eenum
\elem

\bp 
Clearly {\bf(1)} is by the definition. The assertion {\bf(2)} is proved in \cite[Page 398]{nakai1972}.
By definition, $f$ induces a homeomorphism of Riemannian manifolds and hence {\bf3(a)} holds true. As $X,Y$ are complex varieties,  {\bf3(b)} follows by passage to profinite completions. The assertion {\bf3(c)} is the invariance of dimensions under homeomorphisms   \cite[Theorem 13.22]{LeeJohn-topological-manifolds-book}. To prove {\bf3(d)}, note that one has $$\text{rank}(H_1(\xan,\Z))=2g_X+n_Y-1=2g_Y-n_Y-1=\text{rank}(H_1(\yan,\Z))$$ where $(g_X,n_X)$ (resp. $(g_Y,n_Y)$) is the genus and the number of punctures in $X$ (resp. $Y$). Since $f$ is a homeomorphism,  $f$ preserves the number of punctures i.e. $n_X=n_Y$, and this gives $g_X=g_Y$. This completes the proof.
\ep

\brem 
In \ssep\ref{se:anabel-galois-reps} and \ssep\ref{se:p-adic-hodge} we saw that the notion of anabelomorphisms of $p$-adic fields may be applied to the study the respective categories of Galois representations and $p$-adic Hodge theories,
\Cref{le:anab-arch-prop} allows us to compare local systems and Hodge structures on anabelomorphic complex varieties. This is taken up in the next few subsections.
\erem

\brem 
From \cite[Appendix 4E, Theorem, Page 409]{nakai-book}, one sees that if $X\anabelmap Y$ is an anabelomorphism of connected quasi-projective complex varieties, then $X, Y$ both possess a Green's function or neither one does, i.e. the possession or non-possession of a Green's function is an amphoric property. While we do not use this property here, we remark that Green's functions play a central role in Arakelov Theory and hence this property is of interest. 
\erem

\subsection{Anabelomorphy and ordinary linear differential equations}
Following \cite{deligne-diff-eq}, by an \textit{ordinary linear differential equation on $\xan$}, we mean a pair $(M^{an},\nabla^{an})$ consisting of a finite rank vector bundle $M^{an}$  on $\xan$ and  a $\C$-linear, integrable connection $\nabla^{an}$ on $M^{an}$. An \textit{algebraic ordinary linear differential on $X$} is a pair $(M,\nabla)$ on $X$ consisting of a  finite rank vector bundle $M$  on $X$ and  a $\C$-linear, integrable connection $\nabla$ on $M$.

The following result is the complex analytic analog of the results of \ssep\ref{se:anabel-galois-reps} (\Cref{le:irred-amphoric}, \Cref{th:automorphic-ordinary-synchronization}) and is motivated by \Cref{th:anab-diff-eq1}. One could say that the theorem below is gluing differential equations on $X$ (resp. $Y$) by their monodromy representations.

\bthm\label{th:anab-arch-diff-eq} 
Let $X\anabmapright{f}Y$ be an anabelomorphism of connected, smooth, quasi-projective varieties over $\C$. Then 
\benumlab
\item one has an equivalence between the categories of ordinary linear differential equations on $\xan$ and  $\yan$ respectively. 
\item and it takes the local system underlying $(M^{an},\nabla^{an})$ on $\xan$ to  the local system underlying an ordinary linear differential equation $(N^{an},\nabla^{an})$ on $\yan$. 
\item  The correspondence given by {\bf(1)}, takes an algebraic ordinary linear differential equation  with regular singular points on $X$ (for a choice of smooth compactification of $X$ with a normal crossings divisor) is mapped to an algebraic ordinary linear differential equation with regular singular points on $Y$ (for a choice of smooth compactification of $Y$ with a normal crossings divisor).
\eenum
\ethm
\bp 
By \cite[Th\'eor\`eme 2.17]{deligne-diff-eq}, there is an equivalence of categories between ordinary linear differential equations $(M^{an},\nabla^{an})$  of rank $n$ and the category of local systems on $X$ i.e. is the category of finite dimensional representations $\rho:\pi_1(X)\to \gln \C$. By  \Cref{le:anab-arch-prop}, any anabelomorphism $X\anabmapright{f} Y$ induces an isomorphism $f^{-1}:\pi_1(Y)\mapright{\isom} \pi_1(X)$. Hence, the composite $\pi_1(Y)\mapright{\isom} \pi_1(X)\mapright{\rho} \gln \C$ provides a representation $\rho':\pi_1(Y)\to \gln{\C}$. By \cite[Th\'eor\`eme 2.17]{deligne-diff-eq}, this gives rise to an ordinary differential equation $(N^{an},\nabla^{an})$ on $Y^{an}$ with monodromy representation $\rho'$. This proves {\bf(1,2)}.

Now suppose $(M,\nabla)$ is an algebraic, ordinary linear differential equation with regular singular points  on $X$ (see \cite[D\'efinition 4.5]{deligne-diff-eq}--this requires a smooth compactification $\bar{X}$ of $X$ such that $\bar{X}-X$ is a divisor with normal crossings, but from \cite[Proposition 4.4(ii)]{deligne-diff-eq} one obtains independence from the choice of a compactification). Then the analytification $(M,\nabla)\mapsto (M^{an},\nabla^{an})$ provides an ordinary linear differential equation on $\xan$. By {\bf(2)}, this gives rise to a representation $\rho: \pi_1(X)\to \gln{\C}$ of the fundamental group. Then the composite homomorphism  $\pi_1(Y)\mapright{\isom}\pi_1(X)\mapright{\rho} \gln \C$ gives representation of the fundamental group of $Y$. By \cite[Th\'eor\`eme 5.9]{deligne-diff-eq}, one obtains a unique (up to isomorphism) algebraic, ordinary, linear differential equation with regular singular points $(N,\nabla)$ on $Y$.  
\ep

\brem 
It is important to recognize that the anabelomorphism $X\anabmapright{f} Y$ need not be an algebraic or even a holomorphic mapping and hence the correspondence $(M,\nabla)\mapsto (N,\nabla)$ established in {\bf(3)} is highly non-algebraic in general.
\erem

\subsection{Anabelomorphy and Hodge theory}
In this subsection it will be convenient to work with the viewpoint of \cite{hain88}, \cite{deligne-zeta-values}, which requires  one to remember the base-point used to compute fundamental groups (this is also the point of view of \cite{joshi-teich} and \cite[\ssep I3, Page 25]{mochizuki-iut1-4}). The definition of anabelomorphism of complex quasi-projective varieties (\Cref{def:anab-arch}) is broader than has been conventionally used in the literature on mixed Hodge Theory (\cite{hain88} and also in \cite{deligne-zeta-values}, both of which work with morphisms of complex algebraic varieties) and hence while the observations of this section are no doubt elementary (from a certain point of view), the perspective  and emphasis here is Teichm\"uller Theoretic and notably allows non-holomorphic mappings between algebraic varieties. 

Let $X/\C$ be a connected, smooth, quasi-projective variety over $\C$. In this subsection $\pi_1(X,x)=\pi_1(\xan,x)$ will denote the topological fundamental group of the complex manifold $\xan$ computed using a base-point $x\in \xan$. Let $\Z[\pi_1(X,x)]$ be the group ring of $\pi_1(X,x)$ and let $\pi_1(X,x)\to 1$ be the tautological homomorphism to the trivial group. This gives the homomorphism of their respective group rings $\Z[\pi_1(X,x)]\to\Z$. This is the \textit{augmentation homomorphism} and its kernel, $J=\ker(\Z[\pi_1(X,x)]\to\Z)$, is the \textit{augmentation ideal}. One should think of the collection of group rings $\left\{\Z[\pi_1(X,x)]\right\}_{x\in X}$ rather than as a single group ring. All these rings are  all isomorphic to each other because the fundamental group is independent of the choice of the base-point. But remembering $x$ means one has a continuous parameter in play.

\bpro\label{pr:hodge1} Let $X/\C$ be any connected, smooth, quasi-projective variety. Then 
\benumlab
\item the ring $\Z[\pi_1(X,x)]$ and the augmentation ideal $J$ are amphoric, hence
\item for all $n\geq1$ the quotients $\Z[\pi_1(X,x)]/J^n$, and the completion $$\widehat{\Z[\pi_1(X,x)]}=\invlim_{n}\Z[\pi_1(X,x)]/J^n$$ are all amphoric.
\item The lower central series $$\pi_1(X,x)_\mydot=\pi_1(X,x)\supseteq [\pi_1(X,x),\pi_1(X,x)]\supseteq [\pi_1(X,x),[\pi_1(X,x),\pi_1(X,x)]]\supseteq \cdots $$
of $\pi_1(X,x)$
is also amphoric.
\item Write $\pi_1(X,x)_{1}=[\pi_1(X,x),\pi_1(X,x)]$ and  $\pi_1(X,x)_{N}=[\pi_1(X,x),\pi_1(X,x)_{N-1}]$ for $N\geq 2$.
Then  for each $N\geq 1$, the nilpotent group $\pi_1(X,x)^N=\pi_1(X,x)/\pi_1(X,x)_{N}$ and its torsion-free quotient $\pi_1(X,x)^{[N]}=\pi_1(X,x)^N/{\rm Torsion}$ are amphoric.
\item For each $N\geq1$, the Malcev Lie algebra $Lie(\pi_1(X,x)^{[N]})$ attached to $\pi_1(X,x)^{[N]}$ by Malcev's construction and the  unipotent $\Q$-algebraic group whose Lie Algebra is $Lie(\pi_1(X,x)^{[N]})$ are amphoric.
\eenum
\epro

\bp 
Suppose $X\anabmapright{f}Y$ is an anabelomorphism of connected, smooth, quasi-projective complex varieties. Then, since $f$ is a homeomorphism, the mapping $(X,x)\mapsto (Y,f(x))$ induces a natural bijection between the sets $\left\{ (X,x)\right\}_{x\in X}$ and $\left\{ (Y,y)\right\}_{y\in Y}$. Thus, the set of pointed spaces arising from $X$ is amphoric. From \Cref{le:anab-arch-prop} one sees that  an anabelomorphism $X\anabmapright{f}Y$ induces an isomorphism of rings $\Z[\pi_1(X,x)]\isom \Z[\pi_1(Y,y)]$ for all $x\in X$ and all $y\in Y$ and this isomorphism preserves the augmentation ideals on both the sides. Thus, $\Z[\pi_1(X,x)]$ and the augmentation ideal $J\subset \Z[\pi_1(X,x)]$ are amphoric. This proves {\bf(1)}. 
The assertion {\bf(2)} is immediate from {\bf(1)}. \Cref{le:anab-arch-prop} implies {\bf(3)} because of the functorial property of the lower central series of a group. The assertion {\bf(4)} follows from {\bf(3)}. The assertion {\bf(5)} is clear from the constructions \cite[\ssep 9.3--9.7]{deligne-zeta-values} with the nilpotent group $\Gamma=\pi_1(X,x)^N$. 
\ep

\bpro\label{pr:hodge2} 
For each $N\geq 1$, one has a representation of the fundamental group $$\rho_N:\pi_1(X,x)\to {\rm Aut}(\Z[\pi_1(X,x)]/J^N)$$ given by $g\mapsto (U\mapsto g^{-1} U g)$. This representation preserves the filtration by ideals $J^\mydot/J^N\subset \Z[\pi_1(X,x)]/J^N$ and hence is a unipotent representation of $\pi_1(X,x)$ which is manifestly amphoric. 
\epro
\bp 
The amphoricity assertion is self-evident and the properties of $\rho_N$ can be found in \cite{hain_zucker}. 
\ep

The next theorem is the precise archimedean analog of \Cref{th:ordinary-amphoric}.

\bthm\label{th:hodge3}
Let $X\anabmapright{\alpha} Y$ be an anabelomorphism of connected, smooth, quasi-projective varieties over $\C$. Then
\benumlab
\item There is a natural mixed Hodge structure on  $\Z[\pi_1(X,x)]/J^N$ (for $N\geq 1$) but, on the other hand, this mixed Hodge structure is not amphoric. 
\item In particular, the Hodge filtration on this Hodge structure is not amphoric.
\item There is an equivalence between the category of unipotent variation of mixed Hodge structures on $X$ and $Y$ respectively in which underlying (unipotent) monodromy representations of $\pi_1(X)$ are identified with the corresponding unipotent representation of $\pi_1(Y)$.
\item The anabelomorphism provides  a natural quasi-equivalence between the categories of sheaves of commutative, differential graded $\Q$-algebras on $X$ and $Y$ respectively.
\eenum
\ethm
\bp
The existence of mixed Hodge structure on $\Z[\pi_1(X,x)]/J^N$ is given by \cite[Theorem 1]{hain1987}. This mixed Hodge structure is not amphoric in general--for instance this is the case already for Riemann surfaces. To see this,  choose two non-isomorphic Riemann surfaces $X,Y$ of genus $g\geq 2$ and a quasi-conformal mapping i.e. a Nakai quasi-isometry $f:\xan \to \yan$ (so $X,Y$ are anabelomorphic by \Cref{def:anab-arch}). 
If one has an isomorphism of mixed Hodge structures $\Z[\pi_1(X,x)]/J^3\isom \Z[\pi_1(Y,y)]/J^3$, then by the main theorem of \cite{pulte1988}, one has an isomorphism $g:(\xan,x)\isom (\yan,y)$ with $g(x)=y$ for all $x\in\xan$ with at most two exceptions.  In particular $g:\xan \isom \yan$ which contradicts the assumption that $X,Y$ are not isomorphic Riemann surfaces.

Now the weight filtration on $\Z[\pi_1(X,x)]/J^N$ is given by the manifestly amphoric filtration $$\Z[\pi_1(X,x)]/J^n\supset J^\mydot/J^n.$$ Since the Hodge structures are not amphoric, one sees that the Hodge filtration is not amphoric. This proves {\bf(1,2)}.

According to \cite[1.3]{hain88}, a unipotent variation of mixed Hodge structures on $Y$ gives rise to a unipotent representation of $\pi_1(Y)$ and composing with the isomorphism $\alpha:\pi_1(X) \mapright{\isom} \pi_1(Y) \to \gln \C$ one obtains a unipotent representation of $\pi_1(X)$ and hence applying \cite[1.3]{hain88} on $X$,  one obtains a unipotent variation of mixed Hodge structures on $X$. This proves {\bf(3)}.

\newcommand{\cda}[1]{{\rm CDA}_\Q(#1)}
The last assertion is a consequence of \cite{navarro-aznar1987}. Let $\alpha^{-1}:Y\anabelmap X$ be the anabelomorphism inverse to the given anabelomorphism $\alpha:X\to Y$ so one obtains a homeomorphism $\alpha:\xan \to\yan$ and its inverse $\alpha^{-1}:\yan\to \xan$ to which one may apply the formalism of \cite{navarro-aznar1987}. Let $\cda X$ (resp. $\cda Y$) be the category of sheaves of commutative differential graded $\Q$-algebras on $X$ (resp.$Y$). By \cite[Lemma 4.8]{navarro-aznar1987}, one sees that the functor 
$$ \R_{TW}(1_X)_*:\cda X\to \cda X$$
and the functor
$$ \R_{TW}\alpha_*^{-1}\circ \R_{TW}\alpha_*:\cda X \to \cda Y\to \cda X$$
are quasi-equivalent. This proves {\bf(4)} and the theorem.
\ep

\brem 
Already for Riemann surfaces one sees that Teichm\"uller Theory plays a role in the structure of the mixed Hodge structures on $\Z[\pi_1(X,x)]/J^n$. A similar phenomenon occurs in \cite{joshi-teich,joshi-teich-rosetta}.
\erem

\subsection{Theta values,  mixed Hodge structures and the archimedean $L$-invariant}\label{ss:thta-mhs-linv}
In \cite{mochizuki-theta}, one finds a construction of cohomology classes in $H^1(G_{K_v},\Q_p(1))$  (for  each prime $v$ of semi-stable reduction) arising from relationship to theta-values \cite[Proposition 1.4(iii)]{mochizuki-theta} (and also in \cite[Example 3.2, Page 79]{mochizuki-iut1-4}). This group describes 2-dimensional, reducible, semi-stable $p$-adic representations of $G_v$. On the other hand, there are no 1-dimensional $\Q$-Hodge structures of weight one, and the Hodge structure of a general elliptic curve is simple as a $\Q$-Hodge structure. Thus, one can ask if the $p$-adic constructions have an archimedean analog. This question is answered here by \Cref{th:ext-classes-at-inf} which
can be thought of as the archimedean analog of \Cref{th:ord-syn-thm} i.e. as the  \emph{Ordinary Synchronization Theorem at Infinity}. 

\bthm\label{th:ext-classes-at-inf}
Let $E/\C$ be an elliptic curve with Schottky parameter $q=q_E$ such that $0<|q|< 1$. Then
\benumlab
\item there is a mixed Hodge structure $H_E\in\Ext^1(\Z(0),\Z(1))\isom \C^*$ whose extension class corresponds to $q\in\C^*$, and
\item this extension coincides with the class $H_\Theta\in\Ext^1(\Z(0),\Z(1))$  corresponding to a value of a suitably normalized reciprocal of a chosen theta-function $\theta$.
\item There is a continuous period mapping from the Teichm\"uller space in genus one $$\mathcal{T}_1=\mathfrak{H} \to \C^*=\Ext^1(\Z(0),\Z(1))$$ 
which assigns an elliptic curve $E_\tau$ with period lattice $[1,\tau]$, the extension class $H_{E_\tau}$ given by {\bf(1)}.
\eenum
\ethm
\bp The first assertion is proved as follows. From \cite[Section 7.1]{deligne-local}, one has an isomorphism of abelian groups:
\be\label{eq:deligne-isom}\Ext^1_{MHS}(\Z(0),\Z(1))=\C^*.\ee
Schottky  uniformization  of elliptic curves says that one has an isomorphism  $$ \C^*/q_E^\Z\mapright{\isom} E(\C).$$ In particular, the Schottky parameter $q_E\in\C^*$ provides a unique mixed Hodge structure  $$H_E\in \Ext^1_{MHS}(\Z(0),\Z(1))=\C^*$$
(not to be confused with the usual Hodge structure $H^1(E^{an},\Z)$ which is pure of weight $1$ and generally simple). The mixed Hodge structure $H_E$ comes equipped with a weight filtration and unipotent monodromy and is explicitly given using the following formula from \cite[7.2]{deligne-local}:
\be\begin{aligned}
H_\C&=\C e_0\oplus\C e_1,\\
W_{-2}\subset H&= \C e_1,\\
F^0\subset H&=\C e_0,\\
H_\Z&=2\pi i \Z e_0\oplus \Z (e_0+\log(q) e_1)\subset H_\C.
\end{aligned}\ee
The mapping $\Z(1)\to H_\Z$ is given by $2\pi i\mapsto 2\pi ie_1$ and $H_\Z\to \Z(0)$ is given by $e_0\mapsto1$. Then one has an exact sequence of mixed Hodge structures
$$0\to\Z(1)\to H\to \Z(0)\to 0,$$
whose class in $\Ext^1_{MHS}(\Z(0),\Z(1))$ is given by $q\in\C^*$. This proves the first assertion.

Now to prove the second assertion. For this, let  $z\in\C$, $q=e^{i\pi\tau}$, with $\tau$ in the upper half plane, so that $0<|q|<1$.  Let $\vartheta_1(q,z)$ and $\vartheta_3(q,z)$ be the  Jacobi Theta functions  on $E/\C$ given by the formulae in \cite[Chapter XXI, Section 21.1]{whittaker-book}. Let  $$\Theta_E(z,q)=\frac{\vartheta_3(0,q)}{q^{-3/4}\cdot\vartheta_1(z,q)}.$$ Then  using the formulae in \cite[Chapter 21, 21.11 and Example 3]{whittaker-book} one checks easily that $$H_\Theta:=\Theta_E(\frac{\pi+\pi\tau}{2},q)=q\in\C^*.$$
Thus $H_\Theta$ arises as  value of a suitably normalized reciprocal of a theta function on $E$ and 
thus provides us a mixed Hodge structure $H_{\Theta}\in \Ext^1_{MHS}(\Z(0),\Z(1))$. 
This proves {\bf(2)}.

The Teichm\"uller space $\mathcal{T}_1$, in genus one, is identified with the upper half-plane $\mathfrak{H}\subset \C$ by \cite[Chapter 1, Theorem 1.2]{imayoshi-book} and the construction of the mapping is clear from {\bf(2)}. The rest of the assertion {\bf(3)} is clear.
\ep

\brem 
The translation between $\Theta(U,z)$ of \cite[Proposition 1.4]{mochizuki-theta} and classical theta function $\vartheta_1(z,q)$ of \cite[Chapter XXI]{whittaker-book} is tedious but not difficult and the  relationship between the two  is 
$$\Theta(e^{i\cdot z},q)=q^{-\frac{1}{8}}\vartheta_1(z,q^{1/2}),$$
(the series on the left should be viewed as a power series in $U=e^{iz}$ for the equality to hold)
but normalization rules are different. 
\erem

Comparing the definition above of $H_E$ with the formula of Fontaine for $\linv$-invariant \cite{colmez10}, we define the archimedean  $\linv$-invariant as follows:
\bdefn\label{def:l-invariant-inf}
Let  ${\rm Log}$  be the principal branch of the complex logarithm.  
Let $$H\in Ext^1_{MHS}(\Z(0),\Z(1))$$ be a mixed Hodge structure. Let $q_H\in\C^*$ be the extension class of $H$ under the isomorphism given by \eqref{eq:deligne-isom}. Then the \emph{archimedean $L$-invariant}, denoted $\linv_\infty(H)$ is defined as
 $$\linv_\infty(H)=\frac{{\rm Log}(q_H)}{2\pi i}.$$
If  $E/\C$ is  given with its Schottky parametrization, then its $\linv_\infty(E)=\linv_\infty(H_E)$ where $H_E$ is the mixed Hodge structure constructed above.
\edefn
The following is immediate from this definition and \Cref{le:anab-arch-prop}:
\bpro\label{pr:linv-inf-unamph}
Let $E/\C$ be an elliptic curve over $\C$. Let $$E(\C)\isom \C^*/q_E^\Z$$ be a Schottky parametrization of $E$, with Schottky parameter $q_E=e^{2\pi i\tau}\in\C^*$ (with $\tau\in\mathfrak{H}$). Then $$\linv_\infty(E)=\linv_\infty(H_E)=\tau.$$ In particular, $\linv_\infty(E)$ is not amphoric (just as in the non-archimedean case \Cref{th:l-inv-unamphoric}).
\epro
\bp 
From \Cref{le:anab-arch-prop} one knows that if $E'$ is a complex quasi-projective variety which is anabelomorphic to $E$, then $\dim(E')=1$ and $E'$ has the same genus and number of punctures as $E$. Thus $E'$ is of genus one with zero punctures and hence is a complex elliptic curve. By \Cref{le:anab-arch-prop}, an anabelomorphism of complex curves $f:E'\to E$ is a quasi-conformal mapping between 1-dimensional complex torii $E',E$, and it is a basic consequence of the fact that the  Teichm\"uller Space in genus one is the upper half plane \cite[Chapter 1]{imayoshi-book} that a quasi-conformal mapping always exists between any pair of complex elliptic curves $E',E$.  Hence, every complex elliptic curve $E'$ is anabelomorphic to the chosen $E$. Now  say, $E'$ has a  period lattice $[1,\tau']$. Then $\linv_\infty(H_{E'})=\tau'$, and this does not coincide with $\tau=\linv_\infty(H_E)$ in general. Thus, $\linv_\infty(H_E)$ is not amphoric.
\ep

\section{Open questions}\label{se:open-question}\nws
We highlight some open questions which may be of wider interest. \begin{question}\label{qu:hdim1}
The notion of anabelomorphy obviously extends to higher dimensional fields considered in \cite{kato1977,kato1978} (this paper considers the $d=1$ case). So the natural question is, to what extent do the results of this paper generalize to the case of higher dimensional local fields?
\end{question}
Here are some specific versions of this question.
\begin{question}
Is there a criterion analogous to \Cref{th:fourth-fun-anab} for deciding if two $d$-dimensional local fields $K,L$ are anabelomorphic? 
\end{question}

\newcommand{\sK}{\mathscr{K}}
Note that Kato's Reciprocity Law for higher dimensional local fields established in \cite{kato1977,kato1978} immediately implies the following
\bpro 
Let $\sK_d(M)$ denote the Milnor $K$-group (in degree $d$) of the field $M$. Suppose $K,L$ are two anabelomorphic $d$-dimensional local fields with (all) successive residue characteristics equal to $p>0$. Then one has an isomorphism  of topological groups
$$ \invlim_{K'} \sK_d(K)/N_{K'/K}\sK_d(K') \isom G_K^{ab}  \isom G_{L}^{ab} \isom \invlim_{L'} \sK_d(L)/N_{L'/L}\sK_d(L')$$
where inverse limits are over all finite abelian extensions $K'/K$ and $L'/L$ respectively and $N_{K'/K}$ (resp. $N_{L'/L}$) is the norm homomorphism and the topology on the respective Milnor $K$-groups is the one defined in \cite{kato1977,kato1978}.
\epro

\begin{question}
This leads to the following question: suppose one has an anabelomorphism of two $d$-dimensional local fields $K,L$ of the same characteristic, then does the anabelomorphism $G_K\isom G_L$ imply that one has a topological isomorphism
$$\sK_d(K) \isom \sK_d(L)?$$ i.e. is $\sK_d(K)$ amphoric? [Here the topology is as defined in \cite{kato1977,kato1978}.] 
\end{question}

\begin{question}
	In the context of the previous question: does the anabelomorphism $G_K\isom G_L$ imply that one has a topological isomorphism
	$$(K^*)^d \isom (L^*)^d$$ compatibly with the isomorphism of their respective quotients $\sK_d(K)\isom \sK_d(L)$?  [Here the topology on $\sK_d(K),\sK_d(L)$ is as defined in \cite{kato1977,kato1978}.] 
\end{question}

\begin{question}
A simpler question is this: suppose $K,L$ are $d$-dimensional local fields. Then is it true that the topological group $K^*$ is amphoric?
\end{question}

In fact, one does not know the answer to the following basic question:
\begin{question}\label{qu:hdim7}
Suppose $K,L$ are two anabelomorphic $d$-dimensional local fields. Then do $K,L$ necessarily have the same successive residue characteristics?
\end{question}

In the context of \Cref{re:finiteness} let us make the following definition
\bdefn\label{def:zero-dim-teich}
Fix an algebraic closure $\bQ_p$ of $\Q_p$ and let $\C_p$ be the completion of $\bQ_p$. For the next few questions, let $\Q_p\subseteq F\subset \bQ_p$ be a $p$-adic field. Let $n\geq 1$ be an integer.  Let $\fT(F)$ be an anabelomorphism class  of all field extensions $F\subset K\subset \bF$ with $[K:F]=n\geq1$. This means each set $\fT(F)$  consists of finite extensions of $K/F$ with $[K:F]=n$ and for all  $K_1,K_2\in \fT(F)$ one has $$G_{K_1}\isom G_{K_2}.$$ In particular each $\fT(F)$ is a finite set because in $\bF$, there are only finitely many field extensions of $F$ of any given degree.
\edefn

\begin{question}
	Is there a nice answer to the question raised in \Cref{re:finiteness}?
\end{question}

\brem\label{re:finite-teich-space}
As discussed in \Cref{ss:intro1.1} and \Cref{re:finiteness}, each discrete set $\fT(F)$ can be considered to be a zero dimensional Arithmetic Teichm\"uller Space because for each $K_1,K_2\in \fT(F)$ one has $$\pi_1^{et}(\Spec(K_1))=G_{K_1}\isom G_{K_2}=\pi_1^{et}(\Spec(K_2))$$ i.e. $\fT(F)$ has properties similar to the classical Teichm\"uller space \cite{imayoshi-book}. Thus, the sums in  Questions \ref{qu:sum1}, \ref{qu:sum2} and \ref{qu:sum3} given below can be considered as sums  or averages over the zero dimensional Arithmetic Teichm\"uller Space $\fT(F)$ in the spirit of similar averages in classical Teichm\"uller Theory (\cite{wright19}). We have provided these questions as prototypes and readers are welcome to formulate and investigate variants of these questions.
\erem

The motivation for \Cref{qu:sum1} and \Cref{qu:sum2} lies in the amphoricity of the Iwasawa cohomology  $H^i_{Iw}(G_K,\Z_p(1))$ given by  \Cref{pr:iwasawa-coh} (the relationship between this cohomology and cyclotomic $p$-adic $L$-functions is studied in \cite{colmez1999} and elsewhere).

\begin{question}\label{qu:sum1}
	With the notation of \Cref{re:finite-teich-space},
suppose $E/F$ is an elliptic curve over a $p$-adic field $F$ and let $L_p(E/K,1)\in\C_p$ be the value of the $p$-adic $L$-function of $E$ considered as an elliptic curve over $K$.  Suppose $n\geq 1$ is an integer. Then do  there exist nice formulae for the sums (one for each $n\geq 1$)
$$\sum_{K\in \fT(F)} L_p(E/K,1) \in \C_p.$$
\end{question}

More generally:
\begin{question}\label{qu:sum2}
	With the notation of \Cref{re:finite-teich-space},
	suppose $E/F$ is an elliptic curve over a $p$-adic field $F$ and let $L_p(E/K,s)$ be the $p$-adic $L$-function of $E/K$ i.e. of $E$ considered as an elliptic curve over $K$.  Suppose $n\geq 1$ is an integer. Let $\fT(F)$ be as defined above. Then what can one say about  the function (one for each $n\geq 1$)
	$$\sum_{K\in \fT(F)} L_p(E/K,s).$$
	For example what can one say about its zeros?
\end{question}

\begin{question}\label{qu:sum3}
	With the notation of \Cref{re:finite-teich-space}, suppose $E/F$ is an elliptic curve over a $p$-adic field $F$. Suppose $\bF$ is an algebraic closure. For finite extension $F\subset K\subset \bF$, let $Tam(E/K)$ be the Tamagawa number of  $E$ considered as an elliptic curve over $K$.  Suppose $n\geq 1$ is an integer. Let $\fT(F)$ be as defined above. Then do  there exist  nice formulae for the sums
	$$\sum_{K\in \fT(F)} Tam(E/K).$$
\end{question}

\begin{question}
	Fix an integer $n\geq 1$. Then is \Cref{th:syn-super-cusp-n} true for  $\GL_n$ and all primes $p$?
\end{question}

\begin{question}\label{qu:langlands-quest}
Can one generalize the results of \ssep\ref{se:anab-langlands} from $\GL_n$ to an arbitrary reductive group $G$?
\end{question}

\begin{question}
The context for this question is the numerical Langlands Correspondence established in \cite{henniart1988}. Let $K,L$ be anabelomorphic $p$-adic fields and $n\geq 1$ be an integer. \Cref{th:ord-syn-thm}, \Cref{th:syn-super-cusp-n}, and \Cref{th:artin-swan-unamphoric} imply that the decomposition, considered in \cite[1.3, 2.6]{henniart1988}, of the sets of irreducible admissible representations of $\gln K$ (resp. $\gln L$) and $n$-dimensional Galois representations using Swan conductors of representations, is not preserved by anabelomorphy. Can one give more precise description of this phenomenon?
\end{question}

\begin{question}\label{qu:p-adic-langlands}
This question is in the context of the $p$-adic Langlands correspondence (see \ssep\ref{ss:anab-p-adic-langlands}). 
\benumlab
\item What are the amphoric (and unamphoric) quantities, structures and properties in the $p$-adic Langlands Program? 
\item Prove (or disprove) \Cref{con:p-adic-langlands1} and \Cref{con:p-adic-langlands2} for $n\geq2$.
\eenum
\end{question}

\begin{question}\label{qu:galois-rep} This question arises from considerations of \ssep\ref{se:anab-connect-thm} and \Cref{th:anab-deform-thry}.
Suppose that $F$ is a finite field or an $\ell$-adic field for some prime $\ell$ and $\anab{K}{K'}{v_1,\ldots,v_n}{w_1,\ldots,w_n}$ are anabelomorphically connected number fields. Let $\rho_K:G_K\to \gln F$ be an irreducible representation of $G_K$. Let $\rho_{K,v_i}=\rho\big\vert_{G_{K_{v_i}}}$ be the restriction of $\rho_K$ to the decomposition group of $v_i$ in $G_K$. Then under what conditions does there exist an (irreducible) representation $\rho_{K'}:G_{K'}\to \gln F$ with $\rho_{K,v_i}=\rho_{K',w_i}$? for $1\leq i\leq n$. 
\end{question}

\begin{question}\label{qu:deform-galois}
In the notation of \Cref{qu:galois-rep}, suppose $\rho_K$ arises from an automorphic representation $\pi_K$. Then, under what circumstances does $\rho_{K'}$ arise from an automorphic representation $\pi_{K'}$? 
\end{question}

\begin{question}
Is \Cref{conj:p-adic-diff-eq} true?
\end{question}

\begin{question}
Prove or disprove \Cref{conj:irreg-unamphoric}.
\end{question}

\section{Appendix: Relationship between Anabelomorphy and Mochizuki's `Indeterminacy Ind1'}\label{ss:anabelomorphy-and-iut}
The idea of Anabelomorphy (\Cref{def:anabelomorphy-fields}) is motivated by the Author's study of Mochizuki's notion of `Indeterminacy Ind1' as discussed  in \cite[Page 416]{mochizuki-iut1-4}, and as used in the main theorems of that paper. As the setting of \iut\ is global, let us fix a number field $L$. Then the simplest way of stating the relationship between the two notions is as follows:
\vskip0.25cm
\begin{tabular}{cc}
	`Indeterminacy Ind1' & = Anabelomorphy at all non-archimedean completions $L_\wp$ of $L$.
\end{tabular}
\vskip0.25cm
\noindent  `Indeterminacy Ind1' is equivalent to considering (and applying) anabelomorphy for all the non-archimedean completions $L_\wp$ of $L$ simultaneously. [The term `Anabelomorphy' was coined by the Author.]

This was explicitly confirmed by Mochizuki (email written to the author on April 22, 2020). Mochizuki provided his  view of how Anabelomorphy  appears (and is used) in \iut, and his explanation is quoted below:
\vskip0.25cm
{\small
In the parlance of \iut, anabelomorphy, in the case of absolute Galois
groups of p-adic local fields, is closely related to
Mochizuki's indeterminacy (Ind1), i.e., to the
${\rm Aut}(G)$-indeterminacy, where $G$ denotes the absolute
Galois group of a p-adic local field, which, in \iut,
occurs \textit{simultaneously} at all non-archimedean primes \textit{of a number field}. In particular, the following results of the present paper:
Theorem~\ref{th:discriminant-is-unamphoric} (and the table following it), Theorem~\ref{th:artin-swan-unamphoric},
and Theorem \ref{th:kodaira-sym-unamphoric}; (and the data tables  after Theorem
\ref{th:kodaira-sym-unamphoric}--see \ssep\ref{ss:add-num-data}) provide explicit numerical insight concerning
how automorphisms of $G$ that do not arise from field
automorphisms, i.e., concerning automorphisms of the
sort that arise in the (Ind1) indeterminacy of \iut,
can act in a fashion that fails to preserve differents,
discriminants, and the Swan and Artin conductors, as
well as several other quantities associated to elliptic
curves and Galois representations that depend, in an
essential way, on the additive structure of the
p-adic field. 
}
\vskip0.25cm

[In \cite[\ssep8.11]{joshi-teich-rosetta}, the Author demonstrates  that all Indeterminacies (Ind1, Ind2, Ind3) considered and used in \iut\ have a natural arithmetic-geometric origin and may be considered as the analogs of the classical fact that the holomorphic structure, and other quantities that go with the holomorphic structure, of a hyperbolic Riemann surface are indeterminate if one is only given the fundamental group of the Riemann surface.]
\clearpage
\thispagestyle{empty}
\begin{sidewaystable}
	\addtocounter{subsection}{1}
	\addtocounter{table}{0} \captionof{table}{Fragment of data on weak amphoricity of invariants of semistable elliptic curves}\label{table:long-numerical-data-table2}
	\begingroup\scriptsize	
	\setlength{\tabcolsep}{10pt} \renewcommand{\arraystretch}{2} \begin{longtable}{|c|c|c|}
		\hline $E/\Q(\zetap)$ & $E/\Q(\zetap,\sqrt[9]{3})$ & $E/\Q(\zetap,\sqrt[9]{4})$ \\
		\hline   $[a_1,a_2,a_3,a_4,a_6]$ & $[v_K(\Delta),f, \text{ K. Sym}, \text{T. num.} ]$ & $[v_L(\Delta),f, \text{ K. Sym}, \text{T. num.} ]$ \\
		\hline $[0, \zetap^5 + \zetap^4 - 6 \zetap^3 - \zetap - 9, 0, \zetap^5 - \zetap^4 + 8 \zetap^2 - \zetap + 12, \zetap^5 + \zetap^2 + 1]$ & $[9, 1, I_{9}, 9]$ & $[9, 1, I_{9}, 9]$\\
		\hline $[0, 2 \zetap^5 - 2 \zetap^4 - \zetap^3 + \zetap - 5, 0, -\zetap^4 + \zetap^3 - 3 \zetap^2 + 8 \zetap + 11, \zetap^5 + \zetap^4 - 2 \zetap^3 + 3 \zetap^2 - \zetap + 1]$ & $[18, 1, I_{18}, 18]$ & $[18, 1, I_{18}, 18]$\\
		\hline $[0, \zetap^5 + \zetap^4 + 24 \zetap^3 + 11 \zetap^2 + 75, 0, -\zetap^5 + 3 \zetap^4 - \zetap^2 + \zetap + 8, \zetap^5 - 3 \zetap^4 + \zetap^3 + \zetap^2 - 2 \zetap - 1]$ & $[18, 1, I_{18}, 18]$ & $[18, 1, I_{18}, 18]$\\
		\hline $[0, \zetap^5 + 2 \zetap^4 + \zetap^3 + 10 \zetap^2 + \zetap + 31, 0, -\zetap^5 + 3 \zetap^4 - \zetap^2 - \zetap - 2, \zetap^5 - 4 \zetap^3 - 7 \zetap - 23]$ & $[18, 1, I_{18}, 18]$ & $[18, 1, I_{18}, 18]$\\
		\hline $[0, -8 \zetap^5 + 8 \zetap^4 - \zetap^2 + \zetap + 4, 0, 2 \zetap^5 + \zetap^3 - 5 \zetap^2 - 2 \zetap - 10, -3 \zetap^5 + \zetap^4 - \zetap^3 - \zetap^2 + 5 \zetap - 22]$ & $[9, 1, I_{9}, 9]$ & $[9, 1, I_{9}, 9]$\\
		\hline $[0, 3 \zetap^4 + 7 \zetap^2 - 4 \zetap + 16, 0, 2 \zetap^5 + \zetap^4 + 8 \zetap^3 - \zetap^2 + 21, \zetap^5 + 3 \zetap^2 - \zetap + 3]$ & $[9, 1, I_{9}, 9]$ & $[9, 1, I_{9}, 9]$\\
		\hline $[0, -\zetap^5 - 7 \zetap^4 + 2 \zetap^2 - 2 \zetap - 12, 0, \zetap^5 - \zetap^4 + \zetap^3 - \zetap + 4, -\zetap^4 - 3 \zetap^2 + \zetap + 3]$ & $[9, 1, I_{9}, 9]$ & $[9, 1, I_{9}, 9]$\\
		\hline $[0, \zetap^5 - \zetap^4 - 6 \zetap^3 - \zetap^2 + 17, 0, 3 \zetap^4 + \zetap^3 + \zetap^2 + 11, 2 \zetap^5 + \zetap^3 - \zetap^2 + 3 \zetap + 1]$ & $[18, 1, I_{18}, 18]$ & $[18, 1, I_{18}, 18]$\\
		\hline $[0, \zetap^4 + 2 \zetap^3 - \zetap^2 - 10 \zetap - 9, 0, \zetap^4 + 2 \zetap^2 + 4, \zetap^5 - 17 \zetap^4 - \zetap^3 + \zetap^2 + 2 \zetap - 34]$ & $[9, 1, I_{9}, 9]$ & $[9, 1, I_{9}, 9]$\\
		\hline $[0, \zetap^5 + 9 \zetap^4 - 6 \zetap^3 + 3 \zetap^2 + \zetap + 17, 0, -\zetap^5 - 274 \zetap^4 + \zetap^3 + \zetap^2 + 2 \zetap - 553, 2 \zetap^5 + \zetap^4 + 6 \zetap^3 - 4 \zetap^2 + 22]$ & $[9, 1, I_{9}, 9]$ & $[9, 1, I_{9}, 9]$\\
		\hline $[0, -2 \zetap^5 + \zetap^3 + 2 \zetap + 45, 0, 3 \zetap^5 - \zetap^4 + 3 \zetap + 11, 2 \zetap^5 - \zetap^4 - 2 \zetap^3 - 8 \zetap^2 + 8 \zetap + 4]$ & $[9, 1, I_{9}, 9]$ & $[9, 1, I_{9}, 9]$\\
		\hline $[0, 2 \zetap^5 + 7 \zetap^3 + \zetap^2 + 27, 0, \zetap^5 - \zetap^3 - 6 \zetap + 1, 11 \zetap^5 + 2 \zetap^4 + 2 \zetap^3 - 8 \zetap^2 + 17]$ & $[9, 1, I_{9}, 9]$ & $[9, 1, I_{9}, 9]$\\
		\hline $[0, -\zetap^5 - \zetap^4 - \zetap^3 - 2 \zetap^2 - \zetap - 14, 0, \zetap^5 - \zetap^4 + \zetap^3 + 7 \zetap^2 - \zetap + 6, 2 \zetap^5 + \zetap^4 + 2 \zetap^2 - 11 \zetap]$ & $[9, 1, I_{9}, 9]$ & $[9, 1, I_{9}, 9]$\\
		\hline $[0, -\zetap^5 - \zetap^4 - \zetap^3 + \zetap^2 - 3, 0, -\zetap^4 - 2 \zetap^3 - 3 \zetap^2 - \zetap - 16, 31 \zetap^5 - 3 \zetap^4 - \zetap^3 + \zetap + 53]$ & $[27, 1, I_{27}, 27]$ & $[27, 1, I_{27}, 27]$\\
		\hline $[0, -4 \zetap^5 - 2 \zetap^4 + \zetap^3 + \zetap^2 - \zetap - 3, 0, 3 \zetap^3 - 10 \zetap^2 - \zetap - 12, \zetap^5 + \zetap^4 + 2 \zetap^3 - \zetap^2 - 10 \zetap - 14]$ & $[18, 1, I_{18}, 18]$ & $[18, 1, I_{18}, 18]$\\
		\hline $[0, \zetap^4 + 2 \zetap^3 - \zetap^2 + 2 \zetap + 12, 0, -\zetap^5 - \zetap^4 - 70 \zetap^2 + \zetap - 129, \zetap^5 - 3 \zetap^4 - 3 \zetap^3 - 13]$ & $[9, 1, I_{9}, 9]$ & $[9, 1, I_{9}, 9]$\\
		\hline $[0, \zetap^5 - 2 \zetap^3 - \zetap^2 + \zetap + 8, 0, -\zetap^4 + 4 \zetap^3 + \zetap^2 + \zetap + 8, 11 \zetap^5 - \zetap^4 + 84 \zetap^3 - 4 \zetap^2 + 183]$ & $[9, 1, I_{9}, 9]$ & $[9, 1, I_{9}, 9]$\\
		\hline $[0, -4 \zetap^5 + 10 \zetap^4 - 8 \zetap^3 - 4 \zetap - 23, 0, -9 \zetap^5 + \zetap^4 - \zetap^3 + \zetap^2 - \zetap - 20, -\zetap^5 + \zetap^4 - \zetap^3 - \zetap^2 - 7]$ & $[27, 1, I_{27}, 27]$ & $[27, 1, I_{27}, 27]$\\
		\hline $[0, 4 \zetap^5 + 3 \zetap^4 - 2 \zetap^2 + 10 \zetap + 40, 0, \zetap^5 - \zetap^4 + 41 \zetap^3 + 86, \zetap^2 + \zetap + 10]$ & $[27, 1, I_{27}, 27]$ & $[27, 1, I_{27}, 27]$\\
		\hline $[0, \zetap^5 + 4 \zetap^4 - 3 \zetap^3 + 3 \zetap^2 + \zetap + 7, 0, -191 \zetap^5 - 3 \zetap^4 + \zetap^3 + \zetap^2 - \zetap - 379, \zetap^5 + 7 \zetap^4 + \zetap^3 + 21]$ & $[9, 1, I_{9}, 9]$ & $[9, 1, I_{9}, 9]$\\
		\hline $[0, -\zetap^4 - 141 \zetap^3 - \zetap^2 + \zetap - 283, 0, -6 \zetap^4 - \zetap^3 - 4 \zetap^2 + \zetap - 16, -\zetap^5 - \zetap^4 + \zetap^2 - \zetap + 11]$ & $[9, 1, I_{9}, 9]$ & $[9, 1, I_{9}, 9]$\\
		\hline $[0, -6 \zetap^5 - \zetap^4 - 4 \zetap^3 + \zetap^2 - 13, 0, 403 \zetap^5 + \zetap^3 - 11 \zetap^2 + 778, 3 \zetap^5 - \zetap^4 - \zetap^2 - \zetap - 75]$ & $[9, 1, I_{9}, 9]$ & $[9, 1, I_{9}, 9]$\\
		\hline $[0, 6 \zetap^5 + 83 \zetap^4 + 8 \zetap^3 - \zetap^2 - \zetap + 194, 0, -9 \zetap^4 + 2 \zetap^3 + \zetap^2 + \zetap - 6, -\zetap^5 + \zetap^4 + 2 \zetap^3 - 2 \zetap^2 - 4 \zetap - 5]$ & $[9, 1, I_{9}, 9]$ & $[9, 1, I_{9}, 9]$\\
		\hline $[0, 24 \zetap^5 + \zetap^4 - 14 \zetap^3 - \zetap^2 + \zetap + 17, 0, -2 \zetap^5 - 2 \zetap^4 + \zetap^3 + 2 \zetap^2 + \zetap + 1, -\zetap^5 + 2 \zetap^4 + 2 \zetap^3 - \zetap + 4]$ & $[9, 1, I_{9}, 9]$ & $[9, 1, I_{9}, 9]$\\
		\hline $[0, -\zetap^5 - \zetap^3 - 5 \zetap^2 - 7, 0, -2 \zetap^5 + \zetap^4 - \zetap^2 - 54 \zetap - 114, 3 \zetap^5 - 4 \zetap^4 - \zetap^2 - 1]$ & $[9, 1, I_{9}, 9]$ & $[9, 1, I_{9}, 9]$\\
		\hline
	\end{longtable} 
	\vskip2mm 
	\endgroup
\end{sidewaystable}
\clearpage

\bibliography{../../master/masterofallbibs.bib}
\end{document}